\documentclass[final]{siamltex}
\usepackage{amsmath} 
\usepackage{amssymb}  
\usepackage{graphics}
\usepackage{color}
\usepackage{epsfig}
\usepackage[ruled,vlined]{algorithm2e}
\usepackage{color}

\newtheorem{remark}{Remark}[section]
\newtheorem{example}{Example}[section]

\newcommand{\beq}{\begin{equation}}
\newcommand{\eeq}{\end{equation}}

\newcommand{\re}{{\mathbb R}}

\title{A general class of iterative splitting methods for solving linear systems}

\author{P. Novati\thanks{Dipartimento di Matematica, Informatica e Geoscienze,
Universit\`{a} di Trieste, I-34100 Trieste (Italy), {\em e-mail:} novati@units.it}
\and
F. Tagliaferro\thanks{Via delle Fiamme Gialle 8/1, I-34123 Trieste (Italy),
{\em e-mail:} fulvio.tagliaferro@gmail.com}
\and
M. Zennaro\thanks{Dipartimento di Matematica, Informatica e Geoscienze, Universit\`{a}
di Trieste, I-34100 Trieste (Italy), {\em e-mail:} zennaro@units.it}}

\begin{document}

\maketitle

\renewcommand{\thefootnote}{\fnsymbol{footnote}}
\renewcommand{\thefootnote}{\arabic{footnote}}

\begin{abstract}
Recently Ahmadi et al. (2021) and Tagliaferro (2022) proposed some iterative methods
for the numerical solution of linear systems which, under the classical hypothesis of
strict diagonal dominance, typically converge faster than the Jacobi method, but
slower than the forward/backward Gauss-Seidel one.
In this paper we introduce a general class of iterative methods, based on
suitable splittings of the matrix that defines the system, which include all of
the methods mentioned above and have the same cost per iteration in
a sequential computation environment.
We also introduce a partial order relation in the set of the splittings and, partly
theoretically and partly on the basis
of a certain number of examples, we show that such partial order is typically
connected to the speed of convergence of the corresponding methods.
We pay particular attention to the case of linear systems for which the Jacobi iteration
matrix is nonnegative, in which case we give a rigorous
proof of the correspondence between the partial order relation and the magnitude of
the spectral radius of the iteration matrices.
Within the considered general class, some new specific promising methods are proposed
as well.
\end{abstract}



\pagestyle{myheadings}
\thispagestyle{plain}
\markboth{}{}

%
%

\section{Introduction}

\label{introduction}

Iterative methods for the numerical solution of a linear system
\begin{equation}
Ax=b,
\label{lin-syst}
\end{equation}
where $A\in \re^{n\times n}$ is the square matrix of coefficients, $b\in \re^n$ is the
known term and $x\in \re^n$ is the solution to be computed, constitute a widely investigated
field of research.
Various efficient methods have been developed, taking into account the order of magnitude
of the dimension $n$ and the various possible characteristics of the matrix $A$.
Moreover, also the type of architecture of the employed computers is carefully taken into
consideration, with particular attention to parallel implementations.

Some of the most popular iterative methods are the {\em Jacobi method} and the
({\em forward} or {\em backward} or {\em symmetric}) {\em Gauss-Seidel methods}
(see, e.g., Stoer and Bulirsch \cite{SB} or Golub and Van Loan \cite{GV}),
designed on a splitting of the coefficient matrix $A$ in three parts as follows:
\begin{equation}
A = C+D+E,
\label{initial-splitting}
\end{equation}
where $C$ and $E$ are the strictly lower and the strictly upper triangular
part of $A$, respectively, and $D$ is its diagonal, whose elements are required to be 
all $\ne 0$.
Of course, it is understood that these three matrices are filled in with additional
zero elements so as to get the $n\times n$ square form.

Then, starting from an initial approximation $x^{(0)}$ to the solution $x$, a sequence
$\{x^{(k)}\}_{k\ge 1}$ is computed such that
\begin{equation}
x^{(k+1)} = -D^{-1}(C+E)x^{(k)}+D^{-1}b
\label{Jacobi}
\end{equation}
in Jacobi's case,
\begin{equation}
x^{(k+1)} = -(D+C)^{-1}Ex^{(k)}+(D+C)^{-1}b
\label{f_Gauss-Seidel}
\end{equation}
and
$$
x^{(k+1)} = -(D+E)^{-1}Cx^{(k)}+(D+E)^{-1}b
$$
in forward and backward Gauss-Seidel's case, respectively, and
\begin{equation}
x^{(k+1)} = (D+E)^{-1}C(D+C)^{-1}Ex^{(k)}+(D+E)^{-1}[I-C(D+C)^{-1}]b
\label{sym_Gauss-Seidel}
\end{equation}
or
\begin{equation}
x^{(k+1)} = (D+C)^{-1}E(D+E)^{-1}Cx^{(k)}+(D+C)^{-1}[I-E(D+E)^{-1}]b
\label{sym_Gauss-Seidel-2}
\end{equation}
in symmetric Gauss-Seidel's case.

So they show the typical iterative form
\begin{equation}
x^{(k+1)} = Bx^{(k)}+c,
\label{general}
\end{equation}
with
$$
B=B_{J}:=-D^{-1}(C+E) \ \ \ {\rm and} \ \ \ c=c_{J}:=D^{-1}b
$$
for Jacobi,
$$
B=B_{fGS}:=-(D+C)^{-1}E \ \ \ {\rm and} \ \ \ c=c_{fGS}:=(D+C)^{-1}b
$$
for forward Gauss-Seidel,
$$
B=B_{bGS}:=-(D+E)^{-1}C \ \ \ {\rm and} \ \ \ c=c_{bGS}:=(D+E)^{-1}b
$$
for backward Gauss-Seidel and
$$
B=B_{sGS}:=(D+E)^{-1}C(D+C)^{-1}E \ \ \ {\rm and}
 \ \ \ c=c_{sGS}:=(D+E)^{-1}[I-C(D+C)^{-1}]b
$$
or
$$
B=B'_{sGS}:=(D+C)^{-1}E(D+E)^{-1}C \ \ \ {\rm and}
 \ \ \ c=c'_{sGS}:=(D+C)^{-1}[I-E(D+E)^{-1}]b
$$
for symmetric Gauss-Seidel.

We remark that, even if the forward and backward Gauss-Seidel iterations look more
elaborate (as they require the solution of a triangular system, which may be solved by
forward or backward substitution, respectively),
their computational complexity is the same as that of the Jacobi iteration.

Moreover, also the symmetric Gauss-Seidel iteration, which just gathers one iteration
of forward and one of backward Gauss-Seidel in sequence, can be suitably rearranged in
order to reach the same computational cost
(see the related discussion in Section~\ref{TC-TR methods}).

Therefore, in a standard sequential computation environment the efficiency comparison
among all the above mentioned methods is based on their respective speed of convergence.

It is well known that an iterative method of the type (\ref{general}) is convergent
if and only if $\rho(B)<1$, where $\rho(\cdot)$ denotes the spectral radius, and
that the smaller is $\rho(B)$, the faster is the convergence.

As for the four particular methods above, whereas the equality $\rho(B_{sGS}) =
\rho(B'_{sGS})$ always holds since $B'_{sGS}$ is just a cyclic permutation of $B_{sGS}$,
it turns out that all the inequalities $\rho(B_{J})<1$, $\rho(B_{fGS})<1$,
$\rho(B_{bGS})<1$ and $\rho(B_{sGS})<1$ are true if the
coefficient matrix $A$ is {\em strictly diagonally dominant} either {\em by rows}
(i.e., if $\|D^{-1}(C+E)\|_\infty = \|B_{J}\|_\infty <1$) or {\em by columns} (i.e., if
$\|(C+E)D^{-1}\|_1 = \|DB_{J}D^{-1}\|_1 <1$).
As is customary, $\| \cdot\|_\infty$ and $\| \cdot\|_1$ stand for the well-known
matrix {\em infinity} and $1$ norm, respectively.

Furthermore, although not being a general rule, in most cases of practical interest and
when all of such methods converge, symmetric Gauss-Seidel is faster than both
forward and backward Gauss-Seidel and, in turn, these last two methods are faster than Jacobi
(i.e., $\rho(B_{sGS}) < \rho(B_{fGS}) < \rho(B_{J})$ and
$\rho(B_{sGS}) < \rho(B_{bGS}) < \rho(B_{J})$).

One important practical case, typically produced by a finite difference method applied to
a linear system of differential equations, is that of a system characterized
by a matrix $A$, which some authors call {\em L-matrix} (see, for example, Zhang, Huang and Liu~\cite{ZHL}), for which the diagonal elements (those of $D$) are strictly
positive and $C+E\le O$ elementwise (see the well-known Stein-Rosenberg theorem and its
generalizations in, e.g., Young~\cite{Y}).

Nevertheless, the Jacobi method is much more suitable to exploit the potential speedup
of a parallel computer environment than the Gauss-Seidel methods, and this fact may
often change the balance of efficiency into its favour.

Recently, Ahmadi et al.~\cite{AMKS} and, subsequently, Tagliaferro~\cite{Ta} proposed
some iterative methods which are a kind of halfway house between Jacobi and forward and
backward Gauss-Seidel.

More precisely, in \cite{AMKS} the Jacobi iteration matrix $B_{J}$ is splitted into
two or more sets of rows and one passes from the previous iterate to the next one by
applying one after the other the various selected sets of rows in a certain order
fixed a priori, each time working on the already modified components.

Differently, in \cite{Ta} the Jacobi iteration matrix $B_{J}$ is splitted into the lower
and the upper triangular parts and one passes from the previous iterate to the next one
by applying first the upper triangle and then, on the so modified vector, the lower
triangle.

These two strategies are different from each other in the way of choosing the splittings
of $B_{J}$, but both of them may be embedded in the more general methodology that we
want to present in this paper.

The methods introduced in \cite{AMKS}, which we shall call {\em AMKS-methods}, will
be carefully analyzed in Section~\ref{IEEE},
whereas the method introduced in \cite{Ta}, which we shall call {\em T$_U$-method},
will be considered soon as an inspiring model to define our general splitting technique.

With a slightly different notation than that used in \cite{Ta}, it looks like a usual
iterative method (\ref{general}) in the doubled dimension space $\re^{2n}$.
Starting from a pair of initial approximations $\{x^{(0)}_1, x^{(0)}_2\}$,
a sequence of pairs of vectors $\{x^{(k)}_1,x^{(k)}_2\}_{k\ge 1}$ is defined by setting
\begin{equation}
\left\{
\begin{array}{l}
x^{(k+1)}_1 = Ux^{(k)}_1+Lx^{(k)}_2+c, \\
x^{(k+1)}_2 = Ux^{(k+1)}_1+Lx^{(k)}_2+c,
\end{array}
\right.
\label{Tagliaf1}
\end{equation}
where
\begin{equation}
L=-D^{-1}C, \ \ \ \ \ U=-D^{-1}E, \ \ \ \ \ c=D^{-1}b.
\label{Tagliaf2}
\end{equation}

It is immediate to see that the iterative scheme (\ref{Tagliaf1}) may be reformulated in
the more compact form
\begin{equation}
X^{(k+1)} = {\cal B}_{T_U} X^{(k)} + \Gamma_{T_U},
\label{Tagliaf3}
\end{equation}
where
\begin{equation}
{\cal B}_{T_U} = \left[ \begin{array}{cc}I & O \\ -U & I \end{array} \right]^{-1} \cdot
\left[ \begin{array}{cc}U & L \\ O & L \end{array} \right]
= \left[ \begin{array}{cc}U & L \\ U^2 & UL+L \end{array} \right]
\label{Tagliaf4}
\end{equation}
is a $(2n\times 2n)$-matrix and
\begin{equation}
X^{(k)}=\left[ \begin{array}{c}x^{(k)}_1 \\ x^{(k)}_2 \end{array} \right],
\ k\ge 0, \ \ \ {\rm and} \ \ \
\Gamma_{T_U}=\left[ \begin{array}{c}c \\ (U+I)c \end{array} \right]
\label{Tagliaf5}
\end{equation}
are $2n$-vectors, $I$ and $O$ being the identity and the zero $(n\times n)$-matrix,
respectively.

It is obvious that the fixed point of (\ref{Tagliaf3}) (i.e., of (\ref{Tagliaf1}))
is given by $X=[x^T, x^T]^T$, where $x$ is the fixed point
of (\ref{Jacobi}), that is the solution of (\ref{lin-syst}), and that the method is
convergent if and only if $\rho({\cal B}_{T_U})<1$.

Looking carefully at (\ref{Tagliaf1}) reveals that it is possible to perform the
iterations just by proceeding with the pair of $n$-vectors $\{Ux^{(k)}_1, x^{(k)}_2\}$
and that the computational complexity of one iteration still is the
same as for the Jacobi or the Gauss-Seidel methods, at least in a sequential computer
environment.

In a parallel computer environment the situation might instead change to the detriment
of the T$_U$-method (\ref{Tagliaf3}) only if the number of parallel processors available
becomes big enough.
Just to give a rough, and also pessimistic, idea of the situation, if we assume that each
processor can handle the product of a matrix row times a vector in a certain time
$\tau$, with $n$ processors or more available the Jacobi method can perform an iteration
by computing a matrix-vector product all together in the same amount of time $\tau$,
whereas the T$_U$-method needs the computation of two matrix-vector products
in sequence, taking the total time of about $2 \tau$.
However, such a number of available processors is often unrealistic when dealing with
very large systems of equations arising from practical applications.

However, in this paper we do not consider any parallel computation issues any longer and
confine ourselves to evaluate the sequential complexity only.

Again without being a general rule, in many cases of practical interest it has been experimentally observed that, when all the three methods converge, the asymptotic rate of
convergence of the T$_U$-method lies in the between of Jacobi and Gauss-Seidel
(i.e., $\rho(B_{fGS}) < \rho({\cal B}_{T_U}) < \rho(B_{J})$ and
$\rho(B_{bGS}) < \rho({\cal B}_{T_U}) < \rho(B_{J})$).

In this paper we introduce a general class of iterative methods by arbitrarily splitting
the matrix $B_{J}$ and dividing
each iteration into successive updating steps.

\begin{definition}
We say that the $d$-tuple of matrices
${\cal M}_{(d)}=\{M_{1},\ldots ,M_{d}\}$, $M_{p} \in \re^{n\times n}$, is a
{\em splitting mask} for $\re^{n\times n}$ if the following three conditions hold:
\begin{itemize}
\item
$M_{p} \ne O$ \ {\rm for all} \ $p=1,\ldots ,d$;
\item
the $M_{p}$'s are {\em logical matrices}, i.e.,
$(M_{p})_{i,j} \in \{0,1\}$ \ {\rm for all} \ $p=1,\ldots ,d$ {\rm and} $i,j=1,\ldots ,n$;
\item
$\sum_{p=1}^d M_{p} = E$, where $E_{i,j}=1$ \ {\rm for all} \ $i,j=1,\ldots ,n$.
\end{itemize}
\label{splitting_mask}
\end{definition}

\vskip 0.2cm

It is immediate to realize that the {\em Hadamard product}
\begin{equation}
M_{p} \circ M_{q} =O \ \  \ {\rm for \ all} \ p,q=1,\ldots ,d \ \ {\rm with} \ \ p\ne q,
\label{hadamard}
\end{equation}
where (as is wellknown)
$$
(M_{p} \circ M_{q})_{i,j}:=(M_{p})_{i,j}(M_{q})_{i,j} \ \ \ {\rm for \ all} \ i,j=1,\ldots ,n.
$$

By applying the mask ${\cal M}_{(d)}$ to a given matrix $B\in \re^{n\times n}$ we get its
{\em corresponding splitting} in two successive steps.

\begin{definition}
Given a matrix $B\in \re^{n\times n}$, we say that the $d$-tuple of matrices
${\cal B}_{(d)}=\{B_{1},\ldots ,B_{d}\}$, $B_{p} \in \re^{n\times n}$,
is the {\em pre-splitting of $B$ corresponding to the mask ${\cal M}_{(d)}$} of order $d$ if
\begin{itemize}
\item
$B_{p}=M_{p} \circ B$ \ {\rm for all} \ $p=1,\ldots ,d$.
\end{itemize}
\label{pre-splitting}
\end{definition}

\vskip 0.2cm

Indeed, it might well be that one or more elements $B_p$ vanish in ${\cal B}_{(d)}$.
On the other hand, as it will be clear in Section~\ref{subsect_refinement}, the possible
presence of an element $B_p=O$ does not modify the behaviour of the iterative method
associated to the splitting ${\cal B}_{(d)} \setminus \{B_p\}$, that does not include it,
and only results in a more complicated and redundant formulation of the method itself.
Therefore, we give following definition.

\begin{definition}
Given a matrix $B\in \re^{n\times n}$, we say that the $d'$-tuple of matrices
${\cal B}'_{(d')}=\{B_{1},\ldots ,B_{d'}\}$, $B_{p} \in \re^{n\times n}$, $d' \le d$,
is the {\em splitting of $B$ corresponding to the mask ${\cal M}_{(d)}$} of order $d$ if
it is obtained by eliminating all the elements $B_p=O$ from the pre-splitting
${\cal B}_{(d)}$, if any, keeping unchanged the sequential order of all the other
nonzero ones.
\label{splitting}
\end{definition}

\vskip 0.2cm

\begin{remark}
It is immediate to see that, by suitably summing up a certain number of
consecutive elements $M_i$ in ${\cal M}_{(d)}$, we can find another splitting mask
${\cal M}'_{(d')}$ for which ${\cal B}'_{(d')}$ turns out to be the
corresponding pre-splitting of $B$ (which, obviously, coincides with the corresponding
splitting of $B$ as well).
\label{splitting-remark}
\end{remark}

\vskip 0.2cm

We can resume the above discussion in the following result, which substantially
allows us to outline our theory without involving again and again
the splitting masks ${\cal M}_{(d)}$.

\begin{proposition}
Given a matrix $B\in \re^{n\times n}$, a splitting (of $B$) of order $d$ is a $d$-tuple
of matrices ${\cal B}_{(d)}=\{B_{1},\ldots ,B_{d}\}$, $B_{p} \in \re^{n\times n}$,
such that:
\begin{itemize}
\item
$B_{p} \ne O$ for all $p=1,\ldots ,d$;
\item
$B=\sum_{p=1}^d B_{p}$;
\item
the Hadamard product $B_{p} \circ B_{q} =O$ for all $p,q=1,\ldots ,d$ with $p\ne q$.
\end{itemize}
\label{proposition_splitting}
\end{proposition}

\vskip 0.2cm

As a first significant example, we consider the splitting mask associated to the T$_U$-method.
It is clear that ${\cal M}_{(2)}={\cal M}_{T_U}:=\{\tilde{U},\tilde{L}\}$, where
$\tilde{u}_{i,j}=1$ and
$\tilde{l}_{i,j}=0$ if $i\le j$ and $\tilde{u}_{i,j}=0$ and $\tilde{l}_{i,j}=1$ if $i>j$,
and that the corresponding pre-splitting of the Jacobi iteration matrix $B_{J}$ is
${\cal B}_{(2)}={\cal B}_{T_U}:=\{U,L\}$.

\begin{remark}
Note that the values of the diagonal elements of $\tilde{U}$ and $\tilde{L}$ are
indeed negligible since the mask ${\cal M}_{(2)}$ is applied to some Jacobi iteration
matrix $B_J$, whose diagonal elements are always equal to zero.

This clearly applies to all splitting masks ${\cal M}_{(d)}$, at least as long as they
are used within the context of the iterative methods and are applied to some Jacobi
matrix $B_J$.
\label{diagonal}
\end{remark}

We also observe that, if $B_J$ is (lower or upper) triangular, then the splitting
of $B_J$ corresponding to the mask ${\cal M}_{(2)}$ reduces to ${\cal B}_{(1)}=\{B_J\}$
(see Definition~\ref{splitting}) and, clearly, the method (\ref{Tagliaf3}) reduces
to the Jacobi one.

The splitting masks associated to the Gauss-Seidel methods are
not straightforward to detect and will be investigated in Section~\ref{T-method}.

The paper is organized as follows.

In Section~\ref{general_class} we define the iterative method corresponding to a
given splitting ${\cal B}_{(d)}$, show the important property of {\em cyclicity} of
the splittings and, on the set of all of them, we introduce the partial order
relation of {\em refinement}. We also give the notion of {\em essentiality} and define
the concept of {\em potentially optimal} splitting.

In Section~\ref{sect_L-matrices} we treat the important case of linear systems for
which the Jacobi iteration matrix $B_J$ is nonnegative.
In this setting we prove that, the more refined is the splitting, the smaller is the spectral radius of the corresponding iteration matrix.
We also briefly treat the opposite case of linear systems for
which the Jacobi iteration matrix $B_J$ is nonpositive.

In Section~\ref{sect_diag_dom} we extend the convergence theorem under the hypothesis
of strict diagonal dominance, well-known for the Jacobi and the Gauss-Seidel methods,
to the whole class of splitting methods that we propose.

In Section~\ref{IEEE} we prove the inclusion of the methods defined by Ahmadi et
al.~\cite{AMKS} in our general class.

In Section~\ref{T-method} we consider again the T$_U$-method and, after proving a couple
of its basic properties, we consider some specific types of refinements which
include, in particular, all the Gauss-Seidel methods (forward, backward and symmetric).

In Section~\ref{section_altrimetodi} we propose a new family of splitting methods,
namely the {\em alternate triangular column/row methods}, which seem to have a good
potential for a fast convergence, in competition with the symmetric Gauss-Seidel method.

Finally, in Section~\ref{numerics} some numerical illustrative examples are given
and in Section~\ref{conclusions} some conclusions are drawn.

%
%

\section{A general class of iterative methods}

\label{general_class}

To begin with, we observe that, possibly at the cost of an initial preconditioning
and a final diagonal-matrix/vector multiplication, it is not
restrictive to assume that the diagonal of the coefficient matrix $A$ be equal to
the identity matrix $I$, so that
$$
L=-C, \ \ \ \ \ U=-E, \ \ \ \ \ c=b
$$
in (\ref{Tagliaf2}) and thus, by (\ref{initial-splitting}),
\begin{equation}
A=I-(L+U)= I-B_{J}.
\label{normalised}
\end{equation}
In fact, many possible important properties of the coefficient matrix $A$, if present,
are not compromised by this assumption.

For example, diagonal dominance by rows is preserved from $A$ to $D^{-1}A$, in which
case we have to consider the modified system, equivalent to (\ref{lin-syst}),
$$
(D^{-1}A)x=D^{-1}b,
$$
and dominance by columns is preserved from $A$ to $AD^{-1}$, in which case
we have to consider the modified equivalent system
$$
(AD^{-1})y=b,
$$
along with the final transformation $x=D^{-1}y$.
Furthermore, also symmetry together with positive (semi-)definiteness are preserved
from $A$ to $D^{-1/2}AD^{-1/2}$, in which case we have to consider the modified
equivalent system
$$
(D^{-1/2}AD^{-1/2})y=D^{-1/2}b,
$$
along with the final transformation $x=D^{-1/2}y$.

With reference to Proposition~\ref{proposition_splitting}, given a splitting
${\cal B}_{(d)}=\{B_{1},\ldots ,B_{d}\}$ of the Jacobi iteration matrix $B_{J}=L+U$,
we define the corresponding iterative scheme
\begin{equation}
\left\{
\begin{array}{ccl}
x^{(k+1)}_1 & = & B_{1}x^{(k)}_1 + B_{2}x^{(k)}_2 + \cdots + B_{d-1}x^{(k)}_{d-1} +
 B_{d}x^{(k)}_d+ c, \\
x^{(k+1)}_2 & = & B_{1}x^{(k+1)}_1 + B_{2}x^{(k)}_2 + \cdots + B_{d-1}x^{(k)}_{d-1} +
 B_{d}x^{(k)}_d + c, \\
\cdot & \cdot & \ \ \ \ \ \ \ \ \ \ \ \ \ \ \ \ \ \ \ \ \ \ \ \ \ \ \ \ \ \cdots \\
\cdot & \cdot & \ \ \ \ \ \ \ \ \ \ \ \ \ \ \ \ \ \ \ \ \ \ \ \ \ \ \ \ \ \cdots \\
\cdot & \cdot & \ \ \ \ \ \ \ \ \ \ \ \ \ \ \ \ \ \ \ \ \ \ \ \ \ \ \ \ \ \cdots \\
x^{(k+1)}_d & = & B_{1}x^{(k+1)}_1 + B_{2}x^{(k+1)}_2 + \cdots + B_{d-1}x^{(k+1)}_{d-1} +
 B_{d}x^{(k)}_d + c,
\end{array}
\right.
\label{general_method}
\end{equation}
which starts from a given $d$-tuple $\{x^{(0)}_1, \cdots , x^{(0)}_d\}$ of $n$-vectors.

It is straightforward to see that
such a scheme may be reformulated in the more compact form
\begin{equation}
X^{(k+1)} = {\cal B}_{(d)} X^{(k)} + \Gamma_{(d)}
\label{general_compact1}
\end{equation}
with
\begin{equation}
{\cal B}_{(d)} = ({\cal I}_{(d)}-{\cal L}_{(d)})^{-1} {\cal U}_{(d)} =
\left( {\cal I}_{(d)} + \sum_{h=1}^{d-1}{\cal L}_{(d)})^{h} \right) {\cal U}_{(d)}
\label{general_compact2}
\end{equation}
(which, without any risk of misunderstanding, we denote
by the same symbol used for the corresponding splitting),
where
\begin{eqnarray*}
{\cal I}_{(d)} = {\rm blockdiag}(I,\cdots,I), \\ \\
{\cal L}_{(d)} = \left[ \begin{array}{ccccc}
O & O & \cdots & O & O \\ 
B_{1} & O & \cdots & O & O \\
\cdot & \cdot & \cdots & \cdot & \cdot \\
\cdot & \cdot & \cdots & \cdot & \cdot \\
\cdot & \cdot & \cdots & \cdot & \cdot \\
B_{1} & B_{2} & \cdots & O & O \\
B_{1} & B_{2} & \cdots & B_{d-1} & O
\end{array} \right]
\end{eqnarray*}
is a nihilpotent block lower-triangular $(dn\times dn)$-matrix,
\begin{eqnarray*}
{\cal U}_{(d)} = \left[ \begin{array}{ccccc}
B_{1} & B_{2} & \cdots & B_{d-1} &  B_{d} \\ 
O & B_{2} & \cdots & B_{d-1} &  B_{d} \\
\cdot & \cdot & \cdots & \cdot & \cdot \\
\cdot & \cdot & \cdots & \cdot & \cdot \\
\cdot & \cdot & \cdots & \cdot & \cdot \\
O & O & \cdots & B_{d-1} & B_{d} \\
O & O & \cdots & O & B_{d}
\end{array} \right]
\end{eqnarray*}
is a block upper-triangular $(dn\times dn)$-matrix and
\begin{eqnarray*}
X^{(k)}=\left[ \begin{array}{c}x^{(k)}_1 \\ x^{(k)}_2 \\ x^{(k)}_3 \\ \cdot \\ \cdot \\
\cdot \\ x^{(k)}_d \end{array} \right], \ \ \
\Gamma_{(d)} = ({\cal I}_{(d)}-{\cal L}_{(d)})^{-1}
\left[ \begin{array}{c}c \\ c \\ c \\ \cdot \\ \cdot \\
\cdot \\ c \end{array} \right]
= \left[ \begin{array}{c}c \\ (B_{1}+I)c \\ (B_{2}+I)(B_{1}+I)c \\\cdot \\
\cdot \\ \cdot \\ (B_{d-1}+I) \cdots (B_{1}+I)c
\end{array} \right]
\end{eqnarray*}
are $dn$-vectors.

Like for the T$_U$-method (\ref{Tagliaf3}), the fixed point of (\ref{general_compact1})
(i.e., of (\ref{general_method})) is given by the $dn$-vector $X=[x^T, \cdots ,x^T]^T$,
where $x$ is the solution of (\ref{lin-syst}), and convergence obviously holds
if and only if $\rho({\cal B}_{(d)})<1$.

Moreover, the computational complexity of one iteration
is the same as for the Jacobi or the Gauss-Seidel methods, at least for the number of
multiplications.
This fact might be realized better by looking at the two-step iterative scheme
$$
\left\{
\begin{array}{ccl}
x^{(k+1)}_1 & = & x^{(k)}_d + B_{d} (x^{(k)}_d-x^{(k-1)}_d), \\
x^{(k+1)}_2 & = & x^{(k+1)}_1 + B_{1} (x^{(k+1)}_1-x^{(k)}_1), \\
\cdot & \cdot & \ \ \ \ \ \ \ \ \ \ \ \ \ \ \ \ \cdots \\
\cdot & \cdot & \ \ \ \ \ \ \ \ \ \ \ \ \ \ \ \ \cdots \\
\cdot & \cdot & \ \ \ \ \ \ \ \ \ \ \ \ \ \ \ \ \cdots \\
x^{(k+1)}_d & = & x^{(k+1)}_{d-1} + B_{d-1} (x^{(k+1)}_{d-1}-x^{(k)}_{d-1}),
\end{array}
\right.
$$
which is equivalent to (\ref{general_method}) for $k\ge 1$.

Anyway, at least from a formal point of view, the necessary amount of memory involved
in the computations seems to strongly depend on the order $d$ of the
splitting ${\cal B}_{(d)}$. In fact, the dimension of the problem grows
up to $dn$ from the initial value $n$.
Nevertheless, the increase of dimension is mostly due to the adopted formalism and
a smart organization of the data together with a smart formulation of the implementation
algorithm may well reduce significantly the inconvenience.
In any case, in this paper we do not consider this practical aspect any more.

%
%

\subsection{Cyclicity of the splittings}

\label{subsect_cyclicity}

The results of this section are particularly important for analyzing the convergence
properties of the proposed class of methods.
\begin{definition}
Given a matrix $B\in \re^{n\times n}$, we define the {\em cyclic shifting map}
$$
{\cal S}: {\rm Spl}(B) \longrightarrow {\rm Spl}(B)
$$
on the set ${\rm Spl}(B)$ of its splittings by setting
$$
{\cal S}({\cal B}_{(d)}) := \{B_{d},B_{1},\ldots ,B_{d-2},B_{d-1}\}
$$
for any splitting ${\cal B}_{(d)}=\{B_{1},B_{2},\ldots ,B_{d-1},B_{d}\}$
and for any order $d\ge 1$.

Moreover, we say that the splitting ${\cal S}({\cal B}_{(d)})$ is the {\em cyclic shift}
of ${\cal B}_{(d)}$.
\label{def_shift}
\end{definition}

\vskip 0.2cm

Note that, for the limit case $d=1$, we only have ${\cal S}(\{B)\})=\{B\}$.

\begin{lemma}[Cyclicity]
Let ${\cal B}_{(d)}=\{B_{1},B_{2},\ldots ,B_{d-1},B_{d}\}$ be a splitting of the
Jacobi iteration matrix $B_{J}$ of order $d$ and let
${\cal B}_{(d)}'=\{B_{d},B_{1},\ldots ,B_{d-2},B_{d-1}\}$ be its cyclic shift.
If $\lambda \ne 0$ is an eigenvalue of (the iteration matrix corresponding to)
${\cal B}_{(d)}$ and $\big[\alpha_1^T, \ldots, \alpha_d^T\big]^T$ is one of its
related eigenvectors, then $\lambda$ and 
$\big[\alpha_d^T, \lambda \alpha_1^T, \ldots, \lambda \alpha_{d-1}^T\big]^T$
are an eigenvalue and the related eigenvector, respectively, of (the iteration matrix
 corresponding to) ${\cal B}_{(d)}'$.
\label{cyclicity}
\end{lemma}
\begin{proof}
By looking at the iterative scheme (\ref{general_method}) (see also (\ref{general_compact1})
and (\ref{general_compact2})), it is immediate to realize that the assumed hypothesis
is equivalent to
\begin{equation}
(\lambda {\cal L}_{(d)}+{\cal U}_{(d)})\alpha = \lambda \alpha
\label{eigen0}
\end{equation}
which, in expanded form, reads
\begin{equation}
\left\{
\begin{array}{ccc}
B_{1}\alpha_1 + B_{2}\alpha_2 + \cdots + B_{d-1}\alpha_{d-1} +
 B_{d}\alpha_d & = & \lambda \alpha_1, \\
\lambda B_{1}\alpha_1 + B_{2}\alpha_2 + \cdots + B_{d-1}\alpha_{d-1} +
 B_{d}\alpha_d & = & \lambda \alpha_2, \\
\cdots & \cdot & \cdot \\
\cdots & \cdot & \cdot \\
\cdots & \cdot & \cdot \\
\lambda B_{1}\alpha_1 + \lambda B_{2}\alpha_2 + \cdots + B_{d-1}\alpha_{d-1} +
 B_{d}\alpha_d & = & \lambda \alpha_{d-1}, \\
\lambda B_{1}\alpha_1 + \lambda B_{2}\alpha_2 + \cdots + \lambda B_{d-1}\alpha_{d-1} +
 B_{d}\alpha_d & = & \lambda \alpha_d.
\end{array}
\right.
\label{eigen1}
\end{equation}
Now, if we operate a cyclic shift both to the set of equations and to the addends
in the left-hand side and, subsequently, multiply the last $d-1$ equations by
$\lambda$, we obtain
$$
\left\{
\begin{array}{ccc}
B_{d}\alpha_d + \lambda B_{1}\alpha_1 + \lambda B_{2}\alpha_2 + \cdots +
\lambda B_{d-2}\alpha_{d-2} + \lambda B_{d-1}\alpha_{d-1} & = & \lambda \alpha_d, \\
\lambda B_{d}\alpha_d + \lambda B_{1}\alpha_1 + \lambda B_{2}\alpha_2 + \cdots +
\lambda B_{d-2}\alpha_{d-2} + \lambda B_{d-1}\alpha_{d-1} & = & \lambda^2 \alpha_1, \\
\lambda B_{d}\alpha_d + \lambda^2 B_{1}\alpha_1 + \lambda B_{2}\alpha_2 + \cdots +
\lambda B_{d-2}\alpha_{d-2} + \lambda B_{d-1}\alpha_{d-1} & = & \lambda^2 \alpha_2, \\
\cdots & \cdot & \cdot \\
\cdots & \cdot & \cdot \\
\cdots & \cdot & \cdot \\
\lambda B_{d}\alpha_d + \lambda^2 B_{1}\alpha_1 + \lambda^2 B_{2}\alpha_2 + \cdots +
\lambda^2 B_{d-2}\alpha_{d-2} + \lambda B_{d-1}\alpha_{d-1} & = & \lambda^2 \alpha_{d-1},
\end{array}
\right.
$$
which, similarly as before, is easily seen to be equivalent to the desired result.
\end{proof}

\vskip 0.2cm

\begin{remark}
For each $d\ge 1$ it obviously holds that ${\cal S}^d({\cal B}_{(d)})={\cal B}_{(d)}$
for all the splittings ${\cal B}_{(d)}$ of order $d$.
\label{d_cyclic}
\end{remark}

\vskip 0.2cm

In view of the forthcoming statement and similar ones that will be made throughout the
paper, we need to give the following definition.

\begin{definition}
We say that two iterative methods are {\em spectrum-equivalent} if their respective
iteration matrices ${\cal B}$ and ${\cal B}'$, possibly of different dimensions
$n$ and $n'$, have the same nonzero eigenvalues.
\label{asympt_equiv}
\end{definition}

\vskip 0.2cm

\begin{remark}
Two spectrum-equivalent iterative methods have the same asymptotic speed of convergence
since their respective iteration matrices ${\cal B}$ and ${\cal B}'$ are such that
$\rho({\cal B}) = \rho({\cal B}')$.
\label{rem_equiv}
\end{remark}

\vskip 0.2cm

In the light of the Remark~\ref{d_cyclic}, it is clear that the iterated application
of the Cyclicity Lemma for $d$ times implies the following corollary.
\begin{corollary}
Let ${\cal B}_{(d)}=\{B_{1},\ldots ,B_{d}\}$ be a splitting of the Jacobi iteration matrix
$B_{J}$ of order $d$.
Then the iterative methods corresponding to the iteration matrices ${\cal B}_{(d)}$ and
${\cal S}^r({\cal B}_{(d)})$ are spectrum-equivalent.
\label{cyclic_spectrum}
\end{corollary}

\vskip 0.2cm

It is interesting to observe that the (iterated) use of the cyclic shifting map ${\cal S}$
not only leads to splittings with the same nonzero eigenvalues of the corresponding iteration matrices, but even leave the original iterative methods substantially unchanged.
In fact, by writing twice the equations (\ref{general_method}) for $k$ and $k+1$
in sequence and then by neglecting the first $k-1$ and the last of the resulting
$2k$ equations, we obtain
\begin{equation}
\left\{
\begin{array}{ccl}
x^{(k+1)}_d & = & B_{d}x^{(k)}_d + B_{1}x^{(k+1)}_1 + \cdots + B_{d-2}x^{(k+1)}_{d-2}
+ B_{d-1}x^{(k+1)}_{d-1} + c, \\
x^{(k+2)}_1 & = & B_{d}x^{(k+1)}_d + B_{1}x^{(k+1)}_1 + \cdots + B_{d-2}x^{(k+1)}_{d-2}
+ B_{d-1}x^{(k+1)}_{d-1} + c, \\
x^{(k+2)}_2 & = & B_{d}x^{(k+1)}_d + B_{1}x^{(k+2)}_1 + \cdots + B_{d-2}x^{(k+1)}_{d-2}
+ B_{d-1}x^{(k+1)}_{d-1} + c, \\
\cdot & \cdot & \ \ \ \ \ \ \ \ \ \ \ \ \ \ \ \ \ \ \ \ \ \ \ \ \ \ \ \ \ \ \ \cdots \\
\cdot & \cdot & \ \ \ \ \ \ \ \ \ \ \ \ \ \ \ \ \ \ \ \ \ \ \ \ \ \ \ \ \ \ \ \cdots \\
\cdot & \cdot & \ \ \ \ \ \ \ \ \ \ \ \ \ \ \ \ \ \ \ \ \ \ \ \ \ \ \ \ \ \ \ \cdots \\
x^{(k+2)}_{d-1} & = & B_{d}x^{(k+1)}_d + B_{1}x^{(k+2)}_1 + \cdots + B_{d-2}x^{(k+2)}_{d-2} +
 B_{d-1}x^{(k+1)}_{d-1} + c.
\end{array}
\right.
\label{cycled_method}
\end{equation}
Now, by setting
$$
y^{(k)}_1:=x^{(k)}_d \ \ \ {\rm and} \ \ \ y^{(k)}_i:=x^{(k+1)}_{i-1}, \ \ \ \ \ k\ge 0,
$$
we get the iterative scheme
\begin{equation}
\left\{
\begin{array}{ccl}
y^{(k+1)}_1 & = & B_{d}y^{(k)}_1 + B_{1}y^{(k)}_2 + \cdots + B_{d-2}y^{(k)}_{d-1}
+ B_{d-1}y^{(k)}_{d} + c, \\
y^{(k+1)}_2 & = & B_{d}y^{(k+1)}_1 + B_{1}y^{(k)}_2 + \cdots + B_{d-2}y^{(k)}_{d-1}
+ B_{d-1}y^{(k)}_{d} + c, \\
y^{(k+1)}_3 & = & B_{d}y^{(k+1)}_1 + B_{1}y^{(k+1)}_2 + \cdots + B_{d-2}y^{(k)}_{d-1}
+ B_{d-1}y^{(k)}_{d} + c, \\
\cdot & \cdot & \ \ \ \ \ \ \ \ \ \ \ \ \ \ \ \ \ \ \ \ \ \ \ \ \ \ \ \ \ \ \ \cdots \\
\cdot & \cdot & \ \ \ \ \ \ \ \ \ \ \ \ \ \ \ \ \ \ \ \ \ \ \ \ \ \ \ \ \ \ \ \cdots \\
\cdot & \cdot & \ \ \ \ \ \ \ \ \ \ \ \ \ \ \ \ \ \ \ \ \ \ \ \ \ \ \ \ \ \ \ \cdots \\
y^{(k+1)}_{d} & = & B_{d}y^{(k+1)}_1 + B_{1}y^{(k+1)}_2 + \cdots + B_{d-2}y^{(k+1)}_{d-1} +
 B_{d-1}y^{(k)}_{d} + c,
\end{array}
\right.
\label{cycled_method_2}
\end{equation}
which corresponds to the cyclic shift ${\cal S}({\cal B}_{(d)})$ and, apart from some
first initial values, produces the same approximations to the
solution $x$ of (\ref{lin-syst}) that are given by the scheme (\ref{general_method}).

The Cyclicity Lemma also simplifies significantly the proof of the next result, which
regards the eigenvectors of the iteration matrices.
\begin{proposition}
Let ${\cal B}_{(d)}=\{B_{1},\ldots ,B_{d}\}$ be a splitting of the Jacobi iteration
matrix $B_{J}$ of order $d$, let $\lambda \ne 0$ be an eigenvalue of (the iteration matrix)
${\cal B}_{(d)}$ and $\big[\alpha_1^T, \ldots, \alpha_d^T\big]^T$ be one of its related
eigenvectors.
Then it holds that
\begin{equation}
\alpha_p \ne 0 \ \ \ {\rm for \ all} \ \ \ p=1,\ldots ,d.
\label{nonzero2}
\end{equation}
\label{nonzero1}
\end{proposition}
\begin{proof}
Possibly by applying the Cyclicity Lemma~\ref{cyclicity} the necessary number of times,
we can confine ourselves to prove that assuming $\alpha_1=0$ makes an absurde.

Indeed, if it were so, the first equation of (\ref{eigen1}) would imply
$$
B_{2}\alpha_2 + \cdots + B_{d-1}\alpha_{d-1} + B_{d}\alpha_d = 0.
$$
Then substituting into the second equation would yield $\lambda \alpha_2 =0$ and thus,
being $\lambda \ne 0$, $\alpha_2 = 0$ as well.
Consequently, from the previous equality we would also get
$$
B_{3}\alpha_3 + \cdots + B_{d-1}\alpha_{d-1} + B_{d}\alpha_d = 0.
$$

Now it is clear that, by recursively repeating analogous substitutions into all the
remaining equations of (\ref{eigen1}),
we could prove that $\alpha_p=0$ for all $p=1,\ldots ,d$, which is impossible.
\end{proof}

%
%

\subsection{Refinement of splittings and essentiality}

\label{subsect_refinement}

In this section we want to suggest a procedure that, given a splitting ${\cal B}_{(d)}$
of a convergent Jacobi iteration matrix $B_J$, lets us modify it in such a way that, in
most cases, the resulting splitting ${\cal B}'_{(d')}$ determines an iterative scheme with
a higher rate of asymptotic convergence than the one associated to ${\cal B}_{(d)}$.

To this aim, given a matrix $B\in \re^{n\times n}$, we define a strict partial ordering
on the set ${\rm Spl}(B)$ of its splittings.
This partial ordering is defined in two consecutive steps, the former regarding two
splittings of orders $d$ and $d+1$ and the latter two splittings of arbitrary orders $d$
and $d'$ with $d'\ge d+2$.
\begin{definition}
Given two splittings ${\cal B}_{(d)}=\{B_{1},\ldots ,B_{d}\}$ and
${\cal B}'_{(d+1)}= \{B'_{1}, \ldots ,\\ B'_{d+1}\}$ of a matrix $B\in \re^{n\times n}$
of orders $d$ and $d+1$, respectively,
we say that ${\cal B}'_{(d+1)}$ is {\em more refined} than (or, equivalently, is a
{\em refinement} of) ${\cal B}_{(d)}$, and write
${\cal B}'_{(d+1)} \succ {\cal B}_{(d)}$, if there exist two integers $r,s$ with
$0\le r\le d-1$ and $0\le s\le d$ such that, with
$$
{\cal C}_{(d)}=\{C_{1},\ldots ,C_{d}\} := {\cal S}^r({\cal B}_{(d)})
$$
and
$$
{\cal C}'_{(d+1)}=\{C'_{1},\ldots ,C'_{d},C'_{d+1}\} := {\cal S}^s({\cal B}'_{(d+1)}),
$$
the following conditions are satisfied:
\begin{itemize}
\item[{\rm (I)}]
$C'_{p}=C_{p}$ for all $p=1,\ldots ,d-1$;
\item[{\rm (II)}]
the pair of matrices
$\{C'_{d}, C'_{d+1}\}\subseteq {\cal C}'_{(d+1)}$ is a splitting of $C_{d}\in {\cal C}_{(d)}$;
\item[{\rm (III)}]
$C'_{d+1} C'_{d}\ne O$.
\end{itemize}
\label{order1}
\end{definition}

\vskip 0.2cm

\begin{definition}
Given two splittings ${\cal B}_{(d)}$ and ${\cal B}'_{(d')}$ of a matrix
$B\in \re^{n\times n}$ of orders $d$ and $d'$, respectively, with $d'\ge d+2$, we say
that ${\cal B}'_{(d')}$ is {\em more refined} than (or, equivalently, is a {\em refinement}
of) ${\cal B}_{(d)}$, and write ${\cal B}'_{(d')} \succ {\cal B}_{(d)}$, if there exists
a chain of $d'-d-1$ splittings ${\cal B}'_{(d+1)}$, ${\cal B}'_{(d+2)}$,
$\cdots $ , ${\cal B}'_{(d'-1)}$ of increasing orders $d+1,d+2, \ldots, d'-1$,
respectively, such that ${\cal B}'_{(d')} \succ {\cal B}'_{(d'-1)} \
\succ \cdots \succ \ {\cal B}'_{(d+1)} \succ {\cal B}_{(d)}$ in the sense of
Definition~\ref{order1}.
\label{order2}
\end{definition}

\vskip 0.2cm

It is immediate to realize that the strict partial ordering introduced by
Definitions~\ref{order1} and \ref{order2} on ${\rm Spl}(B)$ is {\em compatible} with the
cyclic shifting map ${\cal S}$ in the sense that
\begin{equation}
{\cal B}'_{(d')} \succ {\cal B}_{(d)} \ \iff \ 
{\cal S}^r({\cal B}'_{(d')}) \succ {\cal S}^s({\cal B}_{(d)})
\label{compatible}
\end{equation}
for all integers $r,s$ with $0\le r\le d'-1$ and $0\le s\le d-1$.

It is obvious that the unique element of order $d=1$, i.e., ${\cal B}_{(1)}=\{B\}$, is
always a minimal element in such a strict partial ordering.

Now we show that condition (III) in Definition~\ref{order1} is essential in order to hope
for an effective improvement of a given splitting ${\cal B}_{(d)}$.

\begin{proposition}
With reference to Definition~\ref{order1}, we assume that two splittings
${\cal B}_{(d)}=\{B_{1},\ldots ,B_{d}\}$ and
${\cal B}'_{(d+1)}= \{B'_{1}, \ldots ,B'_{d+1}\}$ of the Jacobi iteration matrix $B_{J}$
of orders $d$ and $d+1$, respectively, satisfy the conditions (I) and (II) but do not
satisfy (III), i.e.,
\begin{equation}
C'_{d+1} C'_{d}=O.
\label{uguale_zero}
\end{equation}
Then the corresponding iterative methods are spectrum-equivalent.
\label{prodotto_zero_raggio}
\end{proposition}
\begin{proof}
First assume that $\lambda \ne 0$ is an eigenvalue of ${\cal B}_{(d)}$.
Thus Corollary~\ref{cyclic_spectrum} implies that $\lambda$ is an eigenvalue of
${\cal C}_{(d)}$ as well.

Now let $\big[\alpha_1^T, \ldots, \alpha_d^T\big]^T$
be the related eigenvector. Then we set
\begin{equation}
\alpha_p' := \alpha_p \ \ \ {\rm for \ all} \ \ \ p=1,\ldots ,d
\label{eq_x}
\end{equation}
and, additionally, define 
\begin{equation}
\alpha_{d+1}' := \alpha_d-\frac{1-\lambda}{\lambda}C'_d \alpha_d,
\label{eq_xx}
\end{equation}
so that (\ref{uguale_zero}) and (\ref{eq_x}) for $p=d$ yield
\begin{equation}
C_{d+1}'\alpha_{d+1}' = C_{d+1}'\alpha_d'.
\label{eq_xxx}
\end{equation}

On the other hand, our hypothesis is equivalent to the set of equalities (\ref{eigen1})
for the splitting ${\cal C}_{(d)}$ and hence, taking into account (I) and (II) of
Definition~\ref{order1} along with (\ref{eq_x}), (\ref{eq_xx}) and (\ref{eq_xxx}), we get
\begin{equation}
\left\{
\begin{array}{ccc}
C'_{1}\alpha'_1 + C'_{2}\alpha'_2 + \cdots + C'_{d}\alpha'_{d} +
 C'_{d+1}\alpha'_{d+1} & = & \lambda \alpha'_1, \\
\lambda C'_{1}\alpha'_1 + C'_{2}\alpha'_2 + \cdots + C'_{d}\alpha'_{d} +
 C'_{d+1}\alpha'_{d+1} & = & \lambda \alpha'_2, \\
\cdots & \cdot & \cdot \\
\cdots & \cdot & \cdot \\
\cdots & \cdot & \cdot \\
\lambda C'_{1}\alpha'_1 + \lambda' C'_{2}\alpha'_2 + \cdots + C'_{d}\alpha'_{d} +
 C'_{d+1}\alpha'_{d+1} & = & \lambda \alpha'_{d}, \\
\lambda C'_{1}\alpha'_1 + \lambda C'_{2}\alpha'_2 + \cdots + \lambda C'_{d}\alpha'_{d} +
 C'_{d+1}\alpha'_{d+1} & = & \lambda \alpha'_{d+1},
\end{array}
\right.
\label{eigen1x}
\end{equation}
which, since $\big[{\alpha'_1}^T, \ldots, {\alpha'_{d+1}}^T\big]^T \ne
\big[0^T, \ldots, 0^T\big]^T$, is equivalent to $\lambda$ being an eigenvalue of
${\cal C}'_{(d+1)}$ and, consequently, of ${\cal B}'_{(d+1)}$ as well (see again
Corollary~\ref{cyclic_spectrum}).

Vice versa, now we assume that $\lambda \ne 0$ is an eigenvalue of ${\cal B}'_{(d+1)}$,
and hence of ${\cal C}'_{(d+1)}$ as well,
and that $\big[{\alpha'_1}^T, \ldots, {\alpha'_{d+1}}^T\big]^T$ is the related
eigenvector of ${\cal C}'_{(d+1)}$,
i.e., that the set of equalities (\ref{eigen1x}) holds true.
Therefore, subtracting the last two equalities from one another yields
\begin{equation}
\lambda (\alpha'_{d+1}-\alpha'_d) = (\lambda -1)C'_{d}\alpha'_{d},
\label{eq_xyz}
\end{equation}
so that, by setting
$$
\alpha_p := \alpha'_p \ \ \ {\rm for \ all} \ \ \ p=1,\ldots ,d,
$$
we get (\ref{eq_xx}).
Moreover, using (\ref{uguale_zero}) and multiplying both sides of (\ref{eq_xx})
by $C_{d+1}'$ lead us to (\ref{eq_xxx}) (with $\alpha'_d=\alpha_d$).

Finally, conditions (I) and (II) of Definition~\ref{order1} allow us to state that
the first $d$ equalities of (\ref{eigen1x}) coincide with the set of equalities
(\ref{eigen1}) for the splitting ${\cal C}_{(d)}$.
In addition, since Proposition~\ref{nonzero1} assures that $\alpha'_p \ne 0$ for
all $p=1,\ldots ,d$, we can claim that 
$\big[\alpha_1^T, \ldots, \alpha_d^T\big]^T \ne \big[0^T, \ldots, 0^T\big]^T$.

In conclusion, $\lambda$ is an eigenvalue of ${\cal C}_{(d)}$,
and hence of ${\cal B}_{(d)}$ as well.
\end{proof}

\vskip 0.2cm

Similarly to what happens when applying the cyclic shifting map ${\cal S}$, now we
shall see that, under the above hypotheses on the splittings ${\cal B}_{(d)}$ and
${\cal B}'_{(d+1)}$, not only the corresponding iteration matrices share the same
nonzero eigenvalues but also, from a practical point of view, the
corresponding iterative schemes coincide to one another.
\begin{proposition}
Under the hypotheses of Proposition~\ref{prodotto_zero_raggio} (where, without loss of
generality and for the sake of simplicity, we assume ${\cal C}'_{(d)}={\cal B}'_{(d)}$
and ${\cal C}'_{(d+1)}={\cal B}'_{(d+1)}$, i.e., $r=s=0$) and with the initial choices
\begin{equation}
y^{(0)}_p = x^{(0)}_p \ \ \ for \ all \ \ \ p=1,\ldots ,d
\label{starting_value_1}
\end{equation}
and
\begin{equation}
y^{(0)}_{d+1} = x^{(0)}_{d} + {\cal B}'_{(d)}u = y^{(0)}_{d} + {\cal B}'_{(d)}u
\ \ \ with \ \ \ u\in \re^n,
\label{starting_value_2}
\end{equation}
the iterative schemes corresponding to the splittings
${\cal B}_{(d)}=\{B_{1},\ldots ,B_{d}\}$ and
${\cal B}'_{(d+1)}=\{B'_{1},\ldots ,B'_{d+1}\}$ give rise to the sequences of approximations
$$
X^{(k)}=\big[ {x^{(k)}_1}^T, \ldots, {x^{(k)}_d}^T \big]^T \ \ \ and \ \ \
Y^{(k)}=\big[ {y^{(k)}_1}^T, \ldots, {y^{(k)}_d}^T, {y^{(k)}_{d+1}}^T \big]^T, \ \ \ k\ge 0,
$$
to the solution $x$ of (\ref{lin-syst}), respectively, where the equalities
\begin{equation}
y^{(k)}_p = x^{(k)}_p \ \ \ for \ all \ \ \ p=1,\ldots ,d
\label{y=x_1}
\end{equation}
and
\begin{equation}
y^{(k)}_{d+1} = x^{(k)}_{d} + {\cal B}'_{(d)} \big( x^{(k)}_{d}-x^{(k-1)}_{d} \big)
= y^{(k)}_{d} + {\cal B}'_{(d)} \big( y^{(k)}_{d}-y^{(k-1)}_{d} \big)
\label{y=x_2}
\end{equation}
hold for all $k\ge 1$.
\label{prodotto_zero_metodo}
\end{proposition}
\begin{proof}
The iterative scheme associated to the splitting ${\cal B}'_{(d+1)}$ is given by
\begin{equation}
\left\{
\begin{array}{ccl}
y^{(k+1)}_1 & = & B'_{1}y^{(k)}_1 + B'_{2}y^{(k)}_2 + \cdots + B'_{d}y^{(k)}_d +
 B'_{d+1}y^{(k)}_{d+1} + c, \\
y^{(k+1)}_2 & = & B'_{1}y^{(k+1)}_1 + B'_{2}y^{(k)}_2 + \cdots + B'_{d}y^{(k)}_{d} +
 B'_{d+1}y^{(k)}_{d+1} + c, \\
\cdot & \cdot & \ \ \ \ \ \ \ \ \ \ \ \ \ \ \ \ \ \ \ \ \ \ \ \ \ \ \ \ \cdots \\
\cdot & \cdot & \ \ \ \ \ \ \ \ \ \ \ \ \ \ \ \ \ \ \ \ \ \ \ \ \ \ \ \ \cdots \\
\cdot & \cdot & \ \ \ \ \ \ \ \ \ \ \ \ \ \ \ \ \ \ \ \ \ \ \ \ \ \ \ \ \cdots \\
y^{(k+1)}_d & = & B'_{1}y^{(k+1)}_1 + B'_{2}y^{(k+1)}_2 + \cdots + B'_{d}y^{(k)}_{d} +
 B'_{d+1}y^{(k)}_{d+1} + c, \\
y^{(k+1)}_{d+1} & = & B'_{1}y^{(k+1)}_1 + B'_{2}y^{(k+1)}_2 + \cdots + B'_{d}y^{(k+1)}_{d} +
 B'_{d+1}y^{(k)}_{d+1} + c.
\end{array}
\right.
\label{general_method_(d+1)}
\end{equation}
Therefore, since (\ref{starting_value_2}) and (\ref{uguale_zero}) imply
$$
B'_{d+1}y^{(0)}_{d+1} = B'_{d+1}y^{(0)}_{d},
$$
substituting (\ref{starting_value_1}) in (\ref{general_method_(d+1)}) for $k=0$ and using
the conditions (I) and (II) of Definition~\ref{order1} make the first $d$ equalities
to coincide with the iterative scheme (\ref{general_method}) for $k=0$.
Consequently, (\ref{y=x_1}) is proved for $k=1$.
Furthermore, computing the difference between the last two equalities clearly yields 
(\ref{y=x_2}) for $k=1$ as well.

Finally, it is also evident that the same computations may be done to perform the
general induction step from $k$ to $k+1$, concluding the proof.
\end{proof}

\vskip 0.2cm

\begin{remark}
The foregoing Propositions~\ref{prodotto_zero_raggio} and \ref{prodotto_zero_metodo}
clearly prove that, as anticipated in Section~\ref{introduction},
a possible element $B_p=O$ may well be eliminated from a splitting ${\cal B}_{(d)}$
without remotely modifying the behaviour of the iterative method associated to it.
\label{zero-remark}
\end{remark}

\vskip 0.2cm

In the light of Propositions~\ref{prodotto_zero_raggio} and \ref{prodotto_zero_metodo}
it makes sense to concentrate on splittings that enjoy the following
{\em essentiality} property.

\begin{definition}
We say that a splitting ${\cal B}_{(d)}=\{B_{1},\ldots ,B_{d}\}$ of a matrix $B$
of order $d$ is {\em nonessential} if there exists an integer $r$ with $0\le r\le d-1$
such that, with
$$
{\cal C}_{(d)}=\{C_{1},\ldots ,C_{d}\} := {\cal S}^r({\cal B}_{(d)}),
$$
we have 
\begin{equation}
C_{d} \, C_{d-1}=O.
\label{uguale_zero_bis}
\end{equation}
Otherwise ${\cal B}_{(d)}$ is said to be {\em essential}.
\label{essentiality}
\end{definition}

\vskip 0.2cm

It is obvious that the splitting of order $d=1$, i.e., ${\cal B}_{(1)}=\{B\}$,
is essential.

\begin{remark}
It is quite evident that, for a given splitting ${\cal B}_{(d)}$, both the properties of
being essential and being nonessential are preserved when operating a (iterated)
cyclic shifting.
\label{rem_essentiality_shift}
\end{remark}

\vskip 0.2cm

Indeed, a nonessential splitting determines an iterative scheme which, substantially,
produces the same approximations to the solution $x$ of (\ref{lin-syst}) as those
produced by another splitting possessing one element less.
Therefore, the former splitting is computationally less efficient than the latter since
it formally involves iterates $X^{(k)}$ of bigger dimension and, also, shows a smaller
potential for parallelism.

Furthermore, without being necessarily a general rule, later in this paper (see,
in particular, the forthcoming Section~\ref{sect_L-matrices}) we shall see that,
in case of convergence of the Jacobi method (i.e., $\rho(B_J)<1$),
two splittings ${\cal B}'_{(d')}$ and ${\cal B}_{(d)}$ of the Jacobi iteration matrix
$B_{J}$ of orders $d'$ and $d$ (with $d'\ge d+1 \ge 2$), respectively,
that satisfy the condition ${\cal B}'_{(d')} \succ {\cal B}_{(d)} \succ \{B_J\}$
are very often such that $\rho({\cal B}'_{(d')}) < \rho({\cal B}_{(d)}) < \rho(B_J)$.

Therefore, we think it worth giving the following qualitative definition.
\begin{definition}
We say that a splitting ${\cal B}_{(d)}=\{B_{1},\ldots ,B_{d}\}$ of order $d$ of the
Jacobi iteration matrix $B_{J}$ is {\em potentially optimal} if it is both essential and
maximal in the strict partial ordering introduced on ${\rm Spl}(B_{J})$
by Definitions~\ref{order1} and \ref{order2}.
\label{optimality}
\end{definition}

%
%

\subsection{A general necessary condition for convergence}

\label{subsect_convergence_condition}

Now we look for a necessary and sufficient condition for the presence of the eigenvalue
$\lambda = 1$ in the spectrum of the iteration matrix of a given splitting
${\cal B}_{(d)} \in {\rm Spl}(B_{J})$.

It is clear that such an occurence prevents the corresponding iterative scheme
(\ref{general_method}) from converging.

\begin{lemma}
With reference to Definition~\ref{order1}, we assume that two splittings
${\cal B}_{(d)}=\{B_{1},\ldots ,B_{d}\}$ and
${\cal B}'_{(d+1)}= \{B'_{1}, \ldots ,B'_{d+1}\}$ of orders $d$ and $d+1$, respectively,
of the Jacobi iteration matrix $B_{J}$ satisfy conditions (I) and (II).
Then, no matter whether they also satisfy condition (III) or not, it holds that
$\lambda = 1$ is an eigenvalue of (the iteration matrix) ${\cal B}_{(d)}$ if and only if
$\lambda = 1$ is an eigenvalue of (the iteration matrix) ${\cal B}'_{(d+1)}$.
\label{lambda=1}
\end{lemma}
\begin{proof}
If condition (III) does not hold, then the statement is just a particular case of
Proposition~\ref{prodotto_zero_raggio}.

On the other hand, if condition (III) is satisfied, assuming $\lambda = 1$
allows us to repeat the same proof of Proposition~\ref{prodotto_zero_raggio}, the only
difference being that (\ref{eq_xx}) and (\ref{eq_xyz}) directly become
$\alpha_{d+1}' := \alpha_d$ and $\alpha_{d+1}' = \alpha_d'$, respectively.
\end{proof}

\vskip 0.2cm

\begin{proposition}
Let ${\cal B}_{(d)}=\{B_{1},\ldots ,B_{d}\}$ be a splitting of order $d$ of the
Jacobi iteration matrix $B_{J}$.
Then $\lambda = 1$ is an eigenvalue of (the iteration matrix) ${\cal B}_{(d)}$ if and
only if $\lambda = 1$ is an eigenvalue of $B_{J}$.
\label{nec_cond_convergence}
\end{proposition}
\begin{proof}
We consider the {\em chain of splittings} of increasing order $p=1, \ldots ,d$
\begin{equation}
\begin{array}{l}
{\cal B}_{(1)}=\{B_J\}=\{B_{1}+ \ldots +B_{d} \}, \\ 
{\cal B}_{(2)}=\{B_{1}, B_{2}+ \ldots +B_{d} \}, \\ 
{\cal B}_{(3)}=\{B_{1}, B_{2}, B_{3}+ \ldots +B_{d} \}, \\
\ \ \ \ \ \ \ \vdots \\
{\cal B}_{(d-1)}=\{B_{1}, \ldots , B_{d-2}, B_{d-1}+B_{d}\}, \\ 
{\cal B}_{(d)}=\{B_{1},\ldots ,B_{d}\},
\end{array}
\label{chain}
\end{equation}
{\em connecting} $B_J$ and ${\cal B}_{(d)}$.

Since condition (III) may well fail to hold, it is not necessarily ordered in
the sense of Definition~\ref{order1}.
On the other hand, each pair of two consecutive splittings ${\cal B}_{(p)}$ and
${\cal B}_{(p+1)}$, $p=1, \ldots , d-1$, satisfies all the other properties involved by
Definition~\ref{order1} and, therefore, the application of Lemma~\ref{lambda=1} for $d-1$
times along the above chain concludes the proof.
\end{proof}

\vskip 0.2cm

The above result clearly allows us to conclude as follows.

\begin{corollary}
Let ${\cal B}_{(d)}=\{B_{1}, \ldots ,B_{d}\}$ be a splitting of order $d$ of the
Jacobi iteration matrix $B_{J}$.
Then a necessary condition for the convergence of the corresponding iterative scheme
(\ref{general_method}) is that $\lambda = 1$ be not an eigenvalue of $B_{J}$.
\label{lambda=1_bis}
\end{corollary}

\vskip 0.2cm

It is worth stressing that the foregoing results regard the sole case of the eigenvalue
$\lambda = 1$ and do not generalize to the other eigenvalues of unitary modulus.
We illustrate this fact by means of the following counterexample.

\begin{example}
Consider the $3 \times 3$ Jacobi iteration matrix
$$
B_J= \gamma \left[ \begin{array}{ccc}
0 & -1 & -1 \\ 
0.5 & 0 & 0 \\
0 & 0.5 & 0
\end{array} \right] \ \ \ {\rm with} \ \ \ \gamma=1.241706082017...
$$
having three different eigenvalues, two of which being complex conjugate with unitary
modulus (namely, $0.23931... \pm 0.97094... \ i$) and the third one being real (namely,
$-0.47862...$), so that $\rho(B_J)=1$.

A few elementary calculations show that the corresponding T$_U$-method (\ref{Tagliaf3})
has a $6 \times 6$ iteration matrix ${\cal B}_{T_U}$ with three different eigenvalues as
well, one of which being equal to $0$ with multiplicity $4$ and the remaining two ones
being complex conjugate (namely, $0.38545... \pm 0.57449... \ i$) with a modulus
strictly less than $1$, so that $\rho({\cal B}_{T_U})<~1$.
\label{esempio-1}
\end{example}

%
%

\section{The case of irreducible nonnegative Jacobi iteration matrices}

\label{sect_L-matrices}

In this section we consider the particular case of matrices $A$ for
which the Jacobi iteration matrix $B_J=-D^{-1}(C+E)=L+U$ satisfies the
{\em nonnegativity condition}
\begin{equation}
B_J \ge O \ \ \ \ \ {\rm elementwise}.
\label{L-matrix}
\end{equation}

They include an important class of matrices, which some authors call {\em L-matrices} 
(see, for example, Zhang, Huang and Liu~\cite{ZHL}).

We shall often assume that the matrix $A$ is {\em irreducible}, too.
As is wellknown, this means that there do not exist any permutation matrix $P$
such that $PAP^{-1}$ is block upper triangular.

Such a requirement is not particularly restrictive in the linear system solving
environment. In fact, if the matrix $A$ were reducible, it would be possible to split
the system into two smaller subsystems to solve one after the other.

It is immediate to see that, if $A$ is irreducible, the Jacobi iteration matrix $B_J$
is irreducible as well (and vice-versa).
Therefore, the well-known Perron-Frobenius Theorem assures that the spectral radius
$\rho(B_J)$ is a simple leading eigenvalue of the Jacobi iteration matrix
$B_J$, whose corresponding eigenvector $\alpha$ may be chosen to be {\em positive}, i.e.,
$$
\alpha > 0 \ \ \ \ \ {\rm elementwise}.
$$

%
%

\subsection{Splittings of nonnegative Jacobi iteration matrices}

\label{L-matrices_splittings}

Let us consider a splitting ${\cal B}_{(d)}=\{B_{1},\ldots ,B_{d}\}$ of
the Jacobi iteration matrix $B_{J}$.
It is immediate to see that condition (\ref{L-matrix}) implies
\begin{equation}
{\cal B}_{(d)} \ge O \ \ \ \ \ {\rm elementwise}
\label{B_nonneg}
\end{equation}
also for the iteration matrix (\ref{general_compact2}) of the related iterative scheme
(\ref{general_method}).

Unfortunately, the possible irreducibility of $B_{J}$ is not inherited by
${\cal B}_{(d)}$ so that, in principle, the general Perron-Frobenius theory only assures
that $\rho({\cal B}_{(d)})$ is an eigenvalue of ${\cal B}_{(d)}$, not
necessarily simple, and that there exists a nonnegative eigenvector $\alpha_{(d)}$
associated to it, not necessarily positive.
Nevertheless, we will see that, if $B_{J}$ is irreducible, then the positivity
of $\alpha_{(d)}$ is assured all the same.

\begin{lemma}
Let the Jacobi iteration matrix $B_{J}$ satisfy condition (\ref{L-matrix})
and let ${\cal B}_{(d)}=\{B_{1},\ldots ,B_{d}\}$ be a splitting of $B_{J}$.
Then the iteration matrix (\ref{general_compact2}) of the related iterative scheme
(\ref{general_method}) has a nonnegative leading eigenvector $\alpha_{(d)}=
\big[\alpha_1^T, \ldots, \alpha_d^T\big]^T$ corresponding
to the eigenvalue $\lambda := \rho({\cal B}_{(d)})$ (assumed to be $>0$) such that
\begin{equation}
0 \le \lambda \alpha_{1} \le \alpha_{d} \le \alpha_{d-1} \le \ldots \le \alpha_{2}
\le \alpha_{1} \ \ \ \ \ {\rm elementwise}
\label{chain-1}
\end{equation}
if $\lambda \le 1$ and
\begin{equation}
0 \le \alpha_{1} \le \alpha_{2} \le \ldots \le \alpha_{d-1} \le \alpha_{d}
\le \lambda \alpha_{1} \ \ \ \ \ {\rm elementwise}
\label{chain-2}
\end{equation}
if $\lambda \ge 1$.
\label{nonneg_eigen}
\end{lemma}
\begin{proof}
Using the expression (\ref{general_compact2}) of the iteration matrix ${\cal B}_{(d)}$,
as already done in the proof of the Cyclicity Lemma~\ref{cyclicity},
the equation for $\lambda$ and $\alpha_{(d)}$ may be written in the expanded form
(\ref{eigen1}).

Taking the difference of each pair of consecutive equations in (\ref{eigen1}) we get
$$
\lambda (\alpha_{p+1} - \alpha_p) = (\lambda - 1)B_p \alpha_p, \ \ \ p=1,\ldots ,d-1,
$$
that is
\begin{equation}
\alpha_{p+1} = (I-\sigma B_p)\alpha_p, \ \ \ p=1,\ldots ,d-1,
\label{differences}
\end{equation}
where
\begin{equation}
\sigma := (1-\lambda)/ \lambda.
\label{sigma}
\end{equation}
Moreover, subtracting the last equation to the first one yields
\begin{equation}
\lambda \alpha_1 = (I-\sigma B_d)\alpha_d
\label{differences-2}
\end{equation}
and, because of (\ref{B_nonneg}), we can assume
\begin{equation}
\alpha_p \ge 0, \ \ \ p=1,\ldots ,d, \ \ \ \ \ {\rm elementwise}.
\label{alpha>=0}
\end{equation}

In conclusion, if $\lambda \le 1$, then (\ref{sigma}) yields $\sigma \ge 0$ and, consequently,
since $B_p \ge O$ for all $p=1,\ldots ,d$ elementwise, the equalities
(\ref{differences}) and (\ref{differences-2}) and the inequality (\ref{alpha>=0}) imply
(\ref{chain-1}).

Analogously, if $\lambda \ge 1$, then (\ref{sigma}) yields $\sigma \le 0$ and, consequently,
(\ref{differences}), (\ref{differences-2}) and (\ref{alpha>=0}) imply (\ref{chain-2}).
\end{proof}

\vskip 0.2cm

The next result is an obvious consequence of (\ref{chain-1}) and (\ref{chain-2}).
\begin{corollary}
Let the hypotheses of Lemma~\ref{nonneg_eigen} hold.
Then the nonnegative leading eigenvector
$\alpha_{(d)}=\big[\alpha_1^T, \ldots, \alpha_d^T\big]^T$ corresponding
to the eigenvalue $\rho({\cal B}_{(d)})$ satisfies the following {\rm zero components}
property:
\begin{itemize}
\item[{\rm (z.c.)}]
if for some $i \in \{1,\ldots ,n\}$ and some $q \in \{1,\ldots ,d\}$ we have that
the $i$-th component $(\alpha_{q})_i=0$, then it holds that the $i$-th components
$(\alpha_{p})_i=0$ for all $p=1,\ldots ,d$.
\end{itemize}
\label{zero_components}
\end{corollary}

\vskip 0.2cm

Now we are in a position to state the mentioned positivity result for the leading
eigenvector $\alpha_{(d)}$ of the iteration matrix ${\cal B}_{(d)}$ under the hypothesis
of irreducibility of the Jacobi matrix $B_{J}$.

\begin{proposition}
Let the Jacobi iteration matrix $B_{J}$ be irreducible and satisfy condition (\ref{L-matrix})
and let ${\cal B}_{(d)}=\{B_{1},\ldots ,B_{d}\}$ be a splitting of $B_{J}$.
Then the iteration matrix (\ref{general_compact2}) of the related iterative scheme
(\ref{general_method}) has a positive leading eigenvector $\alpha_{(d)}=
\big[\alpha_1^T, \ldots, \alpha_d^T\big]^T$ satisfying the following stronger form of
(\ref{chain-1}) and (\ref{chain-2}):
\begin{equation}
0 < \lambda \alpha_{1} \le \alpha_{d} \le \alpha_{d-1} \le \ldots \le \alpha_{2}
\le \alpha_{1} \ \ \ \ \ {\rm elementwise}
\label{chain-3}
\end{equation}
if $\lambda \le 1$ and
\begin{equation}
0 < \alpha_{1} \le \alpha_{2} \le \ldots \le \alpha_{d-1} \le \alpha_{d}
\le \lambda \alpha_{1} \ \ \ \ \ {\rm elementwise}
\label{chain-4}
\end{equation}
if $\lambda \ge 1$.
\label{positive_eigen}
\end{proposition}
\begin{proof}
By contradiction, let us assume that $\alpha_{(d)}$ is not positive, but just
nonnegative (see Lemma~\ref{nonneg_eigen}). Then there exist two indices
$i\in \{1,\ldots ,n\}$ and $q\in \{1,\ldots ,d\}$
such that $(\alpha_{q})_i=0$ and, thus, Corollary~\ref{zero_components} implies
$$
(\alpha_{p})_i=0 \ \ \ {\rm for \ all} \ \ \ p=1,\ldots ,d.
$$
Consequently, the first equality in (\ref{eigen1}) yields
\begin{equation}
(B_{1}\alpha_1 + B_{2}\alpha_2 + \cdots + B_{d-1}\alpha_{d-1} + B_{d}\alpha_d)_i = 0.
\label{yyy}
\end{equation}

On the other hand, the assumed irreducibility of $B_{J}$ assures that for any
$k \in \{1,\ldots ,n\}$, $k \ne i$, there exists a chain of indices
$h_1, \ldots ,h_r \in \{1,\ldots ,n\}$, all different from one another and also
from $i$ and $k$, such that
\begin{equation}
b_{i h_1}>0, \ b_{h_1 h_2}>0, \ \ldots \ , b_{h_{r-1} h_r}>0, \ b_{h_r k}>0.
\label{zzz}
\end{equation}
Hence, in view of Definition~\ref{splitting}, using the fact
that $B_p \ge O$ elementwise, $p=1,\ldots ,d$, using the nonnegativity of $\alpha_{(d)}$
and, again, Corollary~\ref{zero_components}, by
(\ref{yyy}) we get
$$
(\alpha_{p})_{h_1}=0, \ \ \ p=1,\ldots ,d.
$$
Therefore, using (\ref{zzz}) we can repeat the same reasoning for another $r$ times
and arrive at
$$
(\alpha_{p})_{k}=0, \ \ \ p=1,\ldots ,d.
$$
The arbitrariness of $k$ let us conclude that
$$
\alpha_{p}=0, \ \ \ p=1,\ldots ,d,
$$
that is $\alpha_{(d)}=0$, which makes the absurde.

Finally, we just observe that the positivity of $\alpha_{(d)}$,
along with (\ref{chain-1}) and (\ref{chain-2}), obviously leads to
(\ref{chain-3}) and (\ref{chain-4}).
\end{proof}

%
%

\subsection{Monotonicity of the splittings}

\label{monotonicity}

Let $B_J$ and $\tilde{B}_J$ be two Jacobi iteration matrices (of the same size $n \times n$)
such that
\begin{equation}
O \le B_J \le \tilde{B}_J \ \ \ \ \ {\rm elementwise}.
\label{ordered_J}
\end{equation}
Then let ${\cal B}_{(d)}=\{B_{1},\ldots ,B_{d}\}$ be a splitting of $B_{J}$ of order $d$
corresponding to a certain mask ${\cal M}_{(d)}$ (see Definition~\ref{splitting} and
Remark~\ref{splitting-remark}) and
$\tilde{{\cal B}}_{(d)}=\{\tilde{B}_{1},\ldots ,\tilde{B}_{d}\}$ be the splitting of
$\tilde{B}_{J}$ corresponding to the same mask ${\cal M}_{(d)}$.

Because of (\ref{ordered_J}), we have that
$$
B_{p} \le \tilde{B}_{p}, \ \ \ p=1,\ldots ,d,  \ \ \ \ \ {\rm elementwise}.
$$
Therefore, by using (\ref{general_compact2}) for both the related
iteration matrices ${\cal B}_{(d)}$ and $\tilde{{\cal B}}_{(d)}$, 
since it clearly holds that ${\cal L}_{(d)} \le \tilde{{\cal L}}_{(d)}$ and
${\cal U}_{(d)} \le \tilde{{\cal U}}_{(d)}$ elementwise, we easily get
\begin{equation}
O \le {\cal B}_{(d)} \le \tilde{{\cal B}}_{(d)} \ \ \ \ \ {\rm elementwise}.
\label{ordered_splittings}
\end{equation}

By applying the Gelfand spectral radius theorem using the matrix infinity norm
$\| \cdot \|_\infty$ (which is monotone with respect to the ``elementwise
$\le$'' order relation) to the inequalities (\ref{ordered_J}) and
(\ref{ordered_splittings}), we easily obtain
\begin{equation}
\rho(B_J) \le \rho(\tilde{B}_J) \ \ \ {\rm and}
\ \ \ \rho({\cal B}_{(d)}) \le \rho(\tilde{{\cal B}}_{(d)}).
\label{ordered_radii}
\end{equation}

%
%

\subsection{Acceleration of convergence}

\label{acceleration}

We have collected all the necessary tools to establish the relationships between
the partial order introduced by Definition~\ref{order1} on ${\rm Spl}(B_{J})$ and
the possible acceleration of convergence of the related iterative schemes for the
cases in which the matrix $A$ defining the system (\ref{lin-syst}) determines a
(irreducible) nonnegative Jacobi iteration matrix $B_{J}$.

We first analyze the more general case in which the matrix $A$ may be reducible.
\begin{lemma}
Let the Jacobi iteration matrix $B_{J}$ with $\rho(B_{J})>0$ satisfy
condition (\ref{L-matrix}) and let  ${\cal B}_{(d)}=\{B_{1},\ldots ,B_{d}\}$ a splitting
of $B_J$ of order $d \ge 2$.
Then the following implications hold:
\begin{itemize}
\item[{\rm (a1)}]
$\rho(B_{J}) < 1 \ \iff \ \rho({\cal B}_{(d)}) < 1$;
\item[{\rm (a2)}]
$\rho(B_{J}) = 1 \ \iff \ \rho({\cal B}_{(d)}) = 1$;
\item[{\rm (a3)}]
$\rho(B_{J}) > 1 \ \iff \ \rho({\cal B}_{(d)}) > 1$.
\end{itemize}
\label{not_slower-1}
\end{lemma}
\begin{proof}
For $\gamma \in [0, + \infty)$ we define the function
$$
f(\gamma) := \rho(\tilde{\cal B}_{(d)}(\gamma)),
$$
where $\tilde{\cal B}_{(d)}(\gamma)$ is the iteration matrix of the splitting corresponding
to the Jacobi matrix $\tilde{B}_{J}(\gamma) := \gamma B_J$
(in particular, $f(1)=\rho({\cal B}_{(d)})$).

The function $f(\gamma)$ is clearly nondecreasing (see (\ref{ordered_splittings}) and
(\ref{ordered_radii})), continuous and such that $f(0)=0$.

Now assume that $\rho(B_{J}) \le 1$ and, by contradiction, that
$\rho({\cal B}_{(d)}) > 1$.
Therefore, since $f(1)=\rho({\cal B}_{(d)}) > 1$, there exists $\bar{\gamma} < 1$ such
that $f(\bar{\gamma}) = \rho(\tilde{\cal B}_{(d)}(\bar{\gamma}))= 1$.
On the other hand, the nonnegativity of $\tilde{\cal B}_{(d)}(\bar{\gamma})$ implies that
$f(\bar{\gamma}) = 1$ is one of its eigenvalues and,
consequently, Proposition~\ref{nec_cond_convergence} implies the absurde chain of
inequalities
$$
1 \ge \rho(B_{J}) > \bar{\gamma}\rho(B_{J}) = \rho(\bar{\gamma}B_{J})
= \rho(\tilde{B}_{J}(\bar{\gamma})) \ge 1.
$$
Summarizing, we have proved that
\begin{equation}
\rho(B_{J}) \le 1 \ \ \ \Longrightarrow \ \ \ \rho({\cal B}_{(d)}) \le 1.
\label{implication}
\end{equation}

Now assume instead that $\rho(B_{J}) \ge 1$ and, by contradiction, that
$\rho({\cal B}_{(d)}) < 1$.
It obviously holds that $\rho(B_{J}/\rho(B_{J})) = 1$ and thus, like before, by the
nonnegativity of $B_{J}$ (which implies that $\rho(B_{J})$ is one of its eigenvalues)
and by using Proposition~\ref{nec_cond_convergence} in the opposite direction, we
obtain the absurde chain of inequalities
$$
1 > \rho({\cal B}_{(d)}) = f(1) \ge f(\bar{1/\rho(B_{J})}) =
\rho(\tilde{\cal B}_{(d)}(1/\rho(B_{J}))) \ge 1.
$$
Summarizing, we have also proved that
$$
\rho(B_{J}) \ge 1 \ \ \ \Longrightarrow \ \ \ \rho({\cal B}_{(d)}) \ge 1
$$
which, along with (\ref{implication}), the nonnegativity of $B_{J}$ and ${\cal B}_{(d)}$
and again Proposition~\ref{nec_cond_convergence}, implies (a1), (a2), (a3).
\end{proof}

\vskip 0.2cm

\begin{theorem}
Let the Jacobi iteration matrix $B_{J}$ with $\rho(B_{J})>0$ satisfy
condition (\ref{L-matrix}) and let  ${\cal B}_{(d)}=\{B_{1},\ldots ,B_{d}\}$ and
${\cal B}'_{(d+1)}= \{B'_{1}, \ldots , B'_{d+1}\}$ be two splittings of $B_J$ of orders
$d$ and $d+1$, respectively, possibly being $d=1$ (i.e., ${\cal B}_{(1)}=\{B_{J})\}$).
Moreover, with reference to Definition~\ref{order1}, let ${\cal B}_{(d)}$
and ${\cal B}'_{(d+1)}$ satisfy conditions (I) and (II) (but not necessarily (III)).

Then the following implications hold:
\begin{itemize}
\item[{\rm (b1)}]
if $\rho(B_{J}) < 1$, then $\rho({\cal B}'_{(d+1)}) \le \rho({\cal B}_{(d)})
\le \rho(B_{J})$;
\item[{\rm (b2)}]
if $\rho(B_{J}) = 1$, then $\rho({\cal B}'_{(d+1)}) = \rho({\cal B}_{(d)}) = 1$;
\item[{\rm (b3)}]
if $\rho(B_{J}) > 1$, then $\rho({\cal B}'_{(d+1)}) \ge \rho({\cal B}_{(d)})
\ge \rho(B_{J})$.
\end{itemize}
\label{not_slower-2}
\end{theorem}

\vskip 0.2cm

\begin{proof}
Without loss of generality and for the sake of simplicity, we assume that
${\cal C}_{(d)}={\cal B}_{(d)}$ and ${\cal C}'_{(d+1)}={\cal B}'_{(d+1)}$ (i.e., $r=s=0$).
Therefore, we have:
\begin{itemize}
\item[(I')]
$B'_{p}=B_{p}$ for all $p=1,\ldots ,d-1$;
\item[(II')]
the pair of matrices
$\{B'_{d}, B'_{d+1}\}\subseteq {\cal B}'_{(d+1)}$ is a splitting of $B_{d}\in {\cal B}_{(d)}$.
\end{itemize}

Since $B_{J}$ satisfies condition (\ref{L-matrix}), we have that also
the iteration matrix
$$
{\cal B}'_{(d+1)} \ge O \ \ \ \ \ {\rm elementwise}.
$$

If condition (III) is not satisfied, i.e., if $B'_{d+1} B'_{d}=O$, by
Proposition~\ref{prodotto_zero_raggio} we immediately get
\begin{equation}
\rho({\cal B}'_{(d+1)}) = \rho({\cal B}_{(d)}).
\label{B'=B}
\end{equation}

Otherwise, if condition (III) is satisfied, Proposition~\ref{prodotto_zero_raggio} can
not be used and we must proceed in a different way.

To this aim, we observe that, by Lemma~\ref{nonneg_eigen}, $\lambda:=\rho({\cal B}'_{(d+1)})$
(assumed to be $>0$ without any restriction) is one of the eigenvalues of
${\cal B}'_{(d+1)}$ and that a corresponding eigenvector has the form
$$
\alpha_{(d+1)}'=\big[{\alpha_1'}^T, \ldots, {\alpha_{d+1}'}^T\big]^T \ge
\big[0^T, \ldots, 0^T\big]^T \ \ \ \ \ {\rm elementwise}.
$$
So, similarly to (\ref{eigen1}), using (I') we get
\begin{equation}
\left\{
\begin{array}{ccc}
B_{1}\alpha_1' + B_{2}\alpha_2' + \cdots + B_{d-1}\alpha_{d-1}' +
 B_{d}'\alpha_d' + B_{d+1}'\alpha_{d+1}' & = & \lambda \alpha_1', \\
\lambda B_{1}\alpha_1' + B_{2}\alpha_2' + \cdots + B_{d-1}\alpha_{d-1}' +
 B_{d}'\alpha_d' + B_{d+1}'\alpha_{d+1}' & = & \lambda \alpha_2', \\
\cdots & \cdot & \cdot \\
\cdots & \cdot & \cdot \\
\cdots & \cdot & \cdot \\
\lambda B_{1}\alpha_1' + \lambda B_{2}\alpha_2' + \cdots + B_{d-1}\alpha_{d-1}' +
 B_{d}'\alpha_d' + B_{d+1}'\alpha_{d+1}' & = & \lambda \alpha_{d-1}', \\
\lambda B_{1}\alpha_1' + \lambda B_{2}\alpha_2' + \cdots + \lambda B_{d-1}\alpha_{d-1}' +
 B_{d}'\alpha_d' + B_{d+1}'\alpha_{d+1}' & = & \lambda \alpha_d', \\
\lambda B_{1}\alpha_1' + \lambda B_{2}\alpha_2' + \cdots + \lambda B_{d-1}\alpha_{d-1}' +
 \lambda B_{d}'\alpha_d' + B_{d+1}'\alpha_{d+1}' & = & \lambda \alpha_{d+1}'.
\end{array}
\right.
\label{eigen1'}
\end{equation}

Now we assume that $\rho(B_{J}) \le 1$, so that Lemma~\ref{not_slower-1} yields
$\lambda \le 1$ as well. Thus, using (II') and (\ref{chain-1}), we arrive at
\begin{equation}
\left\{
\begin{array}{ccc}
B_{1}\alpha_1' + B_{2}\alpha_2' + \cdots + B_{d-1}\alpha_{d-1}' +
 B_{d}\alpha_d' & \ge & \lambda \alpha_1', \\
\lambda B_{1}\alpha_1' + B_{2}\alpha_2' + \cdots + B_{d-1}\alpha_{d-1}' +
 B_{d}\alpha_d' & \ge & \lambda \alpha_2', \\
\cdots & \cdot & \cdot \\
\cdots & \cdot & \cdot \\
\cdots & \cdot & \cdot \\
\lambda B_{1}\alpha_1' + \lambda B_{2}\alpha_2' + \cdots + B_{d-1}\alpha_{d-1}' +
 B_{d}\alpha_d' & \ge & \lambda \alpha_{d-1}', \\
\lambda B_{1}\alpha_1' + \lambda B_{2}\alpha_2' + \cdots + \lambda B_{d-1}\alpha_{d-1}' +
 B_{d}\alpha_d' & \ge & \lambda \alpha_d',
\end{array}
\right.
\label{eigen1''}
\end{equation}
elementwise, which is the expanded form for
$$
(\lambda {\cal L}_{(d)}+{\cal U}_{(d)})\big[{\alpha_1'}^T, \ldots, {\alpha_{d}'}^T\big]^T
\ge \lambda \big[{\alpha_1'}^T, \ldots, {\alpha_{d}'}^T\big]^T
\ \ \ \ \ {\rm elementwise}.
$$
By using (\ref{general_compact2}) and the fact that condition (\ref{L-matrix}) yields
$$
({\cal I}_{(d)}-{\cal L}_{(d)})^{-1} \ge O \ \ \ {\rm and} \ \ \
{\cal U}_{(d)}\ge O \ \ \ \ \ {\rm elementwise},
$$
we finally get
$$
\lambda \big[{\alpha_1'}^T, \ldots, {\alpha_{d}'}^T\big]^T \le
{\cal B}_{(d)} \big[{\alpha_1'}^T, \ldots, {\alpha_{d}'}^T\big]^T 
\ \ \ \ \ {\rm elementwise}.
$$
In turn, we easily obtain
$$
\lambda^k \big[{\alpha_1'}^T, \ldots, {\alpha_{d}'}^T\big]^T \le
{{\cal B}_{(d)}}^k \big[{\alpha_1'}^T, \ldots, {\alpha_{d}'}^T\big]^T
 \ \ \ \ \ {\rm elementwise} \ {\rm for \ all} \ k\ge 1 
$$
and thus, since
$\big[{\alpha_1'}^T, \ldots, {\alpha_{d}'}^T\big]^T \ne \big[0^T, \ldots, 0^T\big]^T$,
by computing the infinity norm $\|\cdot \|_\infty$ of both sides and by applying the
Gelfand spectral radius theorem we can conclude that
$$
\lambda = \rho({\cal B}'_{(d+1)}) \le \rho({\cal B}_{(d)}),
$$
which is obviously implied by (\ref{B'=B}).

Summarizing, using again Lemma~\ref{not_slower-1} also for ${\cal B}_{(d)}$, we have
proved that
\begin{equation}
\rho(B_J) \le 1 \ \ \ \Longrightarrow \ \ \ 
\rho({\cal B}'_{(d+1)}) \le \rho({\cal B}_{(d)}) \le 1,
\label{(d+1)-(d)}
\end{equation}
no matter whether condition (III) is satisfied or not.

Now we consider again the chain of splittings (\ref{chain})
introduced in the proof of Proposition~\ref{nec_cond_convergence}, where each pair of two
consecutive splittings ${\cal B}_{(p)}$ and ${\cal B}_{(p+1)}$, $p=1, \ldots , d-1$,
satisfies all the properties involved by Definition~\ref{order1}, possibly except
condition (III).

Therefore, starting with the pair of splittings $({\cal B}_{(1)} = \{B_J\}, {\cal B}_{(2)})$,
we have two possibilities:
\begin{itemize}
\item
either condition (III) is not satisfied, in which case we apply
Proposition~\ref{prodotto_zero_raggio} and get
$$
\rho({\cal B}_{(2)}) = \rho(B_J),
$$
\item
or condition (III) is satisfied, in which case we can repeat
the previous procedure for such a pair and get
\begin{equation}
\rho({\cal B}_{(2)}) \le \rho(B_J).
\label{B2-B1}
\end{equation} 
\end{itemize}
However, in any case we have that (\ref{B2-B1}) holds.

By repeating in sequence the procedure on all the pairs $({\cal B}_{(p)}, {\cal B}_{(p+1)})$,
$p=1, \ldots , d-1$, we arrive at
$$
\rho({\cal B}_{(d)}) \le \rho(B_J)
$$
and hence, together with (\ref{(d+1)-(d)}), at
\begin{equation}
\rho(B_J) \le 1 \ \ \ \Longrightarrow \ \ \ 
\rho({\cal B}'_{(d+1)}) \le \rho({\cal B}_{(d)}) \le \rho(B_J).
\label{final}
\end{equation}

Conversely, if $\rho(B_J) \ge 1$, by using (\ref{chain-2}) we similarly arrive at
$$
\rho({\cal B}'_{(d+1)}) \ge \rho({\cal B}_{(d)}) \ge 1
$$
and then, by employing once again the chain of splittings (\ref{chain}), we obtain
the symmetric implication of (\ref{final})
$$
\rho(B_J) \ge 1 \ \ \ \Longrightarrow \ \ \ 
\rho({\cal B}'_{(d+1)}) \ge \rho({\cal B}_{(d)}) \ge \rho(B_J).
$$

Therefore, (b1), (b2) and (b3) are proved.
\end{proof}

\vskip 0.2cm

The next corollary is straightforward (its proof is included in the previous one).

\begin{corollary}
Let the Jacobi iteration matrix $B_{J}$ with $\rho(B_{J})>0$ satisfy
condition (\ref{L-matrix}) and let  ${\cal B}_{(d)}=\{B_{1},\ldots ,B_{d}\}$
a splitting of $B_J$ of order $d\ge 2$.

Then the following implications hold:
\begin{itemize}
\item[{\rm (c1)}]
if $\rho(B_{J}) \le 1$, then $\rho({\cal B}_{(d)}) \le \rho(B_{J})$;
\item[{\rm (c2)}]
if $\rho(B_{J}) \ge 1$, then $\rho({\cal B}_{(d)}) \ge \rho(B_{J})$.
\end{itemize}
\label{not_slower-3}
\end{corollary}

\vskip 0.2cm

Now we refine our analysis by assuming that the matrix $A$, and thus also $B_J$, is
irreducible and obtain, as an obvious corollary, the generalization of the classical
Stein-Rosenberg theorem (regarding the Gauss-Seidel methods) to any generic splitting
${\cal B}'_{(d+1)}$ that is more refined than $B_J$.

\begin{theorem}
Let the Jacobi iteration matrix $B_{J}$ with $\rho(B_{J})>0$ be irreducible and satisfy
condition (\ref{L-matrix}) and let  ${\cal B}_{(d)}=\{B_{1},\ldots ,B_{d}\}$ and
${\cal B}'_{(d+1)}= \{B'_{1}, \ldots , B'_{d+1}\}$ be two splittings of $B_J$ of orders
$d$ and $d+1$, respectively, possibly being $d=1$ (i.e., ${\cal B}_{(1)}=\{B_{J})\}$).
Then, if ${\cal B}'_{(d+1)}$ is more refined than ${\cal B}_{(d)}$ (i.e., if
${\cal B}'_{(d+1)} \succ {\cal B}_{(d)}$), the following implications hold:
\begin{itemize}
\item[{\rm (d1)}]
if $\rho(B_{J}) < 1$, then $\rho({\cal B}'_{(d+1)}) < \rho({\cal B}_{(d)})
\le \rho(B_{J})$;
\item[{\rm (d2)}]
if $\rho(B_{J}) = 1$, then $\rho({\cal B}'_{(d+1)}) = \rho({\cal B}_{(d)}) = 1$;
\item[{\rm (d3)}]
if $\rho(B_{J}) > 1$, then $\rho({\cal B}'_{(d+1)}) > \rho({\cal B}_{(d)})
\ge \rho(B_{J})$.
\end{itemize}
\label{faster}
\end{theorem}

\vskip 0.2cm

\begin{proof}
Similarly as in the proof of Theorem~\ref{not_slower-2}, without loss of generality
we assume that (I'), (II') and, furthermore, that
\begin{itemize}
\item[(III')]
$B'_{d+1} B'_{d}\ne O$
\end{itemize}
are satisfied.

It is clear that (d2) holds as it is nothing but (b2).

Then, in order to prove (d1), we assume that $\rho(B_{J}) < 1$ and, by contradiction,
that
\begin{equation}
\lambda := \rho({\cal B}'_{(d+1)}) = \rho({\cal B}_{(d)}) > 0.
\label{equality}
\end{equation}

By Proposition~\ref{positive_eigen} both the iteration matrices ${\cal B}_{(d)}$ and
${\cal B}'_{(d+1)}$ have a positive leading eigenvector, namely
$$
\alpha_{(d)}=\big[{\alpha_1}^T, \ldots, {\alpha_{d}}^T\big]^T >
\big[0^T, \ldots, 0^T\big]^T \ \ \ \ \ {\rm elementwise}
$$
and
$$
\alpha_{(d+1)}'=\big[{\alpha_1'}^T, \ldots, {\alpha_{d}'}^T, {\alpha_{d+1}'}^T\big]^T >
\big[0^T, \ldots, 0^T\big]^T \ \ \ \ \ {\rm elementwise},
$$
respectively, which satisfy the recursive relations
\begin{equation}
\alpha_{p+1} = (I-\sigma B_p)\alpha_p, \ \ \ p=1,\ldots ,d-1,
\label{differences-B}
\end{equation}
and, in view of (I'),
\begin{equation}
\alpha_{p+1}' = (I-\sigma B_p)\alpha_p', \ \ \ p=1,\ldots ,d-1, \ \ \
{\rm and} \ \ \ \alpha_{d+1}' = (I-\sigma B_d')\alpha_d',
\label{differences-B'}
\end{equation}
respectively, where
\begin{equation}
\sigma = (1-\lambda)/ \lambda > 0
\label{sigma-2}
\end{equation}
(see the proof of Lemma~\ref{nonneg_eigen}).

Since the eigenvectors are defined short of a multiplicative constant $\ne 0$, the
positivity of $\alpha_{(d)}$ allows us to assume that
$$
\alpha_p' \le \alpha_p, \ \ \ p=1,\ldots ,d, \ \ \ \ \ {\rm elementwise},
$$
which implies
\begin{equation}
\delta_p := \alpha_p - \alpha_p' \ge 0, \ \ \ p=1,\ldots ,d, \ \ \ \ \ {\rm elementwise}.
\label{delta}
\end{equation}
Therefore, (\ref{differences-B}), (\ref{differences-B'}), (\ref{sigma-2})
and the nonnegativity of the $B_p$'s yield
$$
0 \le \delta_{p+1} = (I-\sigma B_p)\delta_p \le \delta_p, \ \ \ p=1,\ldots ,d-1,
 \ \ \ \ \ {\rm elementwise}
$$
and, in turn,
\begin{equation}
0 \le \delta_d \le \ldots \le \delta_1 \ \ \ \ \ {\rm elementwise}.
\label{delta_chain}
\end{equation}
Finally, by substituting the last equality of (\ref{differences-B'}) into the last but
first equality of (\ref{eigen1'}), by using (II') and (\ref{delta}) and by
subtracting to the last equality of (\ref{eigen1}), we get
\begin{equation}
0 \le \lambda B_{1}\delta_1 + \cdots + \lambda B_{d-1}\delta_{d-1} +
 B_{d}\delta_d + \sigma B_{d+1}' B_d' \alpha_d' = \lambda \delta_d
 \ \ \ \ \ {\rm elementwise}.
\label{main}
\end{equation}

Now, again by the arbitrariness of the multiplicative constant of $\alpha_{(d)}$, the
chain of inequalities (\ref{delta_chain}) allows us to suppose that
\begin{equation}
(\delta_d)_i=0 \ \ \ {\rm for \ some} \ \ \ i \in \{1, \ldots ,n\}
\label{delta_i_zero}
\end{equation}
and, consequently, by (\ref{main}) we get
$$
(\lambda B_{1}\delta_1 + \cdots + \lambda B_{d-1}\delta_{d-1} + B_{d}\delta_d)_i = 0
$$
and
$$
(B_{d+1}' B_d' \alpha_d')_i = 0.
$$

On the other hand, being $\lambda \ne 0$, the irreducibility of the Jacobi iteration
matrix $B_J=B_{1} + \cdots + B_{d}$ is clearly transmitted to the matrix
$\lambda B_{1} + \cdots + \lambda B_{d-1} + B_{d}$.
Therefore, since (\ref{delta_chain}) yields
$$
(\delta_p)_h=0 \ \ \ {\rm for \ some} \ \ \ p \in \{1, \ldots ,d-1\}
\ \ \ \Longrightarrow \ \ \ (\delta_d)_h=0,
$$
following the same flow of the proof of Proposition~\ref{positive_eigen},
we can show that
$$
(\delta_d)_k=0
$$
for any $k \in \{1,\ldots ,n\}$, $k \ne i$, and consequently that
$$
(B_{d+1}' B_d' \alpha_d')_k = 0.
$$

In conclusion, since $\alpha_d' > 0$ and $B_{d+1}' B_d' \ge O$ elementwise, we get
$$
B_{d+1}' B_d' = O,
$$
which contradicts (III'). Hence, (\ref{equality}) can not be true and so, in view of (b1),
the stronger implication (d1) is proved.

The proof of (d3) is analogous.
\end{proof}

%
%

\subsection{Some remarks on the case of nonpositive Jacobi iteration matrices}

\label{subsect_nonpos}

In this section we briefly consider the opposite particular case of matrices $A$ for
which the Jacobi iteration matrix $B_J=-D^{-1}(C+E)=L+U$ satisfies the
{\em nonpositivity condition}
\begin{equation}
B_J \le O \ \ \ \ \ {\rm elementwise}.
\label{<=0-matrix}
\end{equation}

As is usually done, if
$B=[b_{i,j}]_{i,j=1}^n$ is an $(n \times n)$-matrix, we shall denote by
\begin{equation}
|B|:=[\, |b_{i,j}|\, ]_{i,j=1}^n
\label{absval}
\end{equation}
the matrix made up with the absolute values of the entries of $B$.

It is clear that
$$
B_J \le O \ \ \ {\rm elementwise} \ \ \ \ \ \iff \ \ \ \ \ |B_J| = -B_J.
$$

Consequently, we have the following obvious result for the spectrum of the
Jacobi iteration matrix.
\begin{proposition}
Let the Jacobi iteration matrix $B_J$ satisfy (\ref{<=0-matrix}).
Then $\lambda$ is an eigenvalue of $B_J$ if and only if $-\lambda$ is an eigenvalue
of $|B_J|$ and the related eigenvectors coincide.

In particular, it holds that
$$
\rho(B_J) = \rho(|B_J|).
$$
\label{Jacobi-nonpos}
\end{proposition}

From a practical point of view, this means that passing from the nonnegative to the
nonpositive case for the iteration matrix $B_J$ does not change the convergence
properties of the Jacobi method.

On the contrary, in general, this is not the case of any splitting method more refined
than Jacobi.

From a heuristic point of view, this fact seems to be clear enough.
In fact, if ${\cal B}_{(d)}=\{B_{1},\ldots ,B_{d}\}$ is a splitting of order $d$ of the
Jacobi iteration matrix $B_{J} \le O$ elementwise, we consider the splitting
\begin{equation}
{\cal B}_{(d)}':=\{|B_{1}|,\ldots ,|B_{d}|\},
\label{B_d'}
\end{equation}
that is the splitting of $|B_J|$ corresponding to the same mask ${\cal M}_{(d)}$ of
${\cal B}_{(d)}$ (see Definition~\ref{splitting}).
By (\ref{general_compact2}) applied to both the related iteration matrices
${\cal B}_{(d)}'$ and ${\cal B}_{(d)}$ we immediately obtain
$$
|{\cal B}_{(d)}| \le {\cal B}_{(d)}' \ \ \ \ \ {\rm elementwise}
$$
and, consequently,
\begin{equation}
\rho({\cal B}_{(d)}) \le \rho({\cal B}_{(d)}')
\label{rhoB_d<=rhoB_d'}
\end{equation}
(see the analogous derivation of (\ref{ordered_radii}) from (\ref{ordered_splittings})).

Now, due to possible cancellation of positive and negative terms, a splitting that is
more refined than Jacobi, starting already from those with 
just two elements (i.e., of the type ${\cal B}_{(2)}=\{B_{1},B_{2}\}$),
gives rise to an iteration matrix in which some elements of $|{\cal B}_{(d)}|$ are very
likely strictly less than the corresponding elements in ${\cal B}_{(d)}'$.
Therefore, the week inequality in (\ref{rhoB_d<=rhoB_d'}) is very likely to become a
strict inequality.

Indeed it often happens that, the more refined is the splitting, the bigger
is its hope of success in attaining a strict inequality in (\ref{rhoB_d<=rhoB_d'}).

This fact becomes particularly interesting when the Jacobi method does not converge
(i.e., when $\rho(B_J)=\rho(|B_J|) \ge 1$) and, at the same time, although (obviously)
not converging if applied to the system related to the opposite Jacobi iteration matrix
$|B_J|=-B_J$, some more refined splitting methods instead do converge.

Nevertheless, there exist also interesting examples, like that of the symmetric tridiagonal
matrices $A$, in which the particular position of the few nonzero elements in the
corresponding Jacobi iteration matrix $B_J$ cause the failure of the above heuristic
reasoning.

However, it is not a purpose of this paper to state rigorous sufficient conditions for
the success of the strict inequality in (\ref{rhoB_d<=rhoB_d'}).

Some illustrative numerical results will be given in Section~\ref{numerics} only.

%
%

\section{The case of strictly diagonally dominant matrices}

\label{sect_diag_dom}

Now we shall show that the condition of {\em strict diagonal dominance by rows}
of the coefficient matrix $A$, that is,
\begin{equation}
\|B_{J}\|_\infty = \|L+U\|_\infty < 1,
\label{str-diag-dom-row}
\end{equation}
which is obviously sufficient for the convergence of the Jacobi method and is also
known to be sufficient for the convergence of the Gauss-Seidel methods, implies indeed
the convergence of all the splitting methods of the type (\ref{general_method}).

\begin{proposition}
Let ${\cal B}_{(d)}=\{B_{1},\ldots ,B_{d}\}$ be a splitting of order $d$ of the
Jacobi iteration matrix $B_{J}$.
Then the corresponding iteration matrix (\ref{general_compact2}) is such that:
\begin{itemize}
\item
$\|B_{J}\|_\infty \le 1 \ \ \ \Longrightarrow \ \ \
\|{\cal B}_{(d)}\|_\infty = \|B_{J}\|_\infty$;

\vskip 0.2cm

\item
$\|B_{J}\|_\infty > 1 \ \ \ \Longrightarrow \ \ \
\|B_{J}\|_\infty \le \|{\cal B}_{(d)}\|_\infty \le \|B_{J}\|_\infty^d$.
\end{itemize}
\label{norma_inf}
\end{proposition}

\vskip 0.2cm

\begin{proof}
The first n rows of the ($dn \times dn$)-matrix ${\cal B}_{(d)}$ are given by the
($n \times dn$)-matrix $[B_{1} \ldots  B_{d}]$ and, hence,
$$
\|{\cal B}_{(d)}\|_\infty \ge \|[B_{1} \ldots  B_{d}]\|_\infty.
$$
On the other hand, the properties of the splittings (see
Proposition~\ref{proposition_splitting}) clearly imply
\begin{equation}
\|[B_{1} \ldots  B_{d}]\|_\infty = \|B_{J}\|_\infty
\label{uguaglianza}
\end{equation}
and, therefore, that
\begin{equation}
\|{\cal B}_{(d)}\|_\infty \ge \|B_{J}\|_\infty.
\label{disequazione}
\end{equation}

Now consider the transformation
$$
y := [y_1^T, \ldots ,y_d^T]^T := {\cal B}_{(d)}x \ \ \ {\rm with} \ \ \
x:=[x_1^T, \ldots ,x_d^T]^T,
$$ 
where $y_1, \ldots ,y_d,x_1, \ldots ,x_d \in \re^n$, and the iterative scheme
(\ref{general_method}), which suggests the equivalent system of equalities
\begin{equation}
\left\{
\begin{array}{ccc}
y_1 & := & B_{1}x_1 + B_{2}x_2 + \cdots + B_{d-1}x_{d-1} + B_{d}x_d, \\
y_2 & := & B_{1}y_1 + B_{2}x_2 + \cdots + B_{d-1}x_{d-1} + B_{d}x_d, \\
\cdot & \cdot & \cdots \\
\cdot & \cdot & \cdots \\
\cdot & \cdot & \cdots \\
y_d & := & B_{1}y_1 + B_{2}y_2 + \cdots + B_{d-1}y_{d-1} + B_{d}x_d.
\end{array}
\right.
\label{equivalent_system}
\end{equation}
Observe that (\ref{uguaglianza}) and the first of the above equalities ensure that
\begin{equation}
\|y_1\|_\infty \le \|B_{J}\|_\infty \, \|x\|_\infty.
\label{inequality_for_1}
\end{equation}

We start with the case $\|B_{J}\|_\infty \le 1$.

Let $1\le k\le d-1$ and assume by inductive hypothesis that
\begin{equation}
\|y_i\|_\infty \le \|B_{J}\|_\infty \, \|x\|_\infty \, \le \|x\|_\infty,
\ \ \ \ \ i=1, \ldots ,k,
\label{inequality_for_k-1}
\end{equation}
which is true for $k=1$ thanks to (\ref{inequality_for_1}).
Then (\ref{uguaglianza}) and the $(k+1)$-st equality of (\ref{equivalent_system})
clearly imply
$$
\|y_{k+1}\|_\infty \le \|B_{J}\|_\infty \, \|x\|_\infty,
$$
too.
Therefore, (\ref{inequality_for_k-1}) holds up to $k=d$, that is
$$
\|y\|_\infty = \|{\cal B}_{(d)} x\|_\infty \le \|B_{J}\|_\infty \, \|x\|_\infty,
$$
and, consequently,
$$
\|{\cal B}_{(d)}\|_\infty = \|B_{J}\|_\infty.
$$

Then we consider the case $\|B_{J}\|_\infty > 1$. 

Like before, let $1\le k\le d-1$ and assume by inductive hypothesis that
\begin{equation}
\|y_i\|_\infty \le \|B_{J}\|_\infty^k \, \|x\|_\infty,
\ \ \ \ \ i=1, \ldots ,k,
\label{inequality_for_k-2}
\end{equation}
which, again, is true for $k=1$ thanks to (\ref{inequality_for_1}).
Then (\ref{uguaglianza}) and the $(k+1)$-st equality of (\ref{equivalent_system})
clearly imply
$$
\|y_{k+1}\|_\infty \le \|B_{J}\|_\infty \,
\max \{\|B_{J}\|_\infty^k \, \|x\|_\infty, \|x\|_\infty \} =
\|B_{J}\|_\infty^{k+1} \, \|x\|_\infty,
$$
too. 
Therefore, (\ref{inequality_for_k-2}) also holds for $k=d$ and, similarly as before,
we can conclude that
$$
\|B_{J}\|_\infty \le \|{\cal B}_{(d)}\|_\infty \le \|B_{J}\|_\infty^d.
$$
\end{proof}

\vskip 0.2cm

\begin{corollary}
Let ${\cal B}_{(d)}=\{B_{1},\ldots ,B_{d}\}$ be a splitting of order $d$ of the
Jacobi iteration matrix $B_{J}$.
If the matrix $A$ is strictly diagonally dominant by rows, i.e., if
(\ref{str-diag-dom-row}) holds, then the corresponding iterative scheme
(\ref{general_method}) is such that $\rho({\cal B}_{(d)}) < 1$ and, therefore,
is convergent.
\label{convergence_Bd}
\end{corollary}

\vskip 0.2cm

It is interesting to remark that a result analogous to Proposition~\ref{norma_inf}
does not hold for the $1$-norm, as is shown by the following simple counterexample.
\begin{example}
Consider the $2 \times 2$ Jacobi iteration matrix
$$
B_J=\left[ \begin{array}{cc}
0 & 1 \\ 
1 & 0
\end{array} \right]
$$
for which $\|B_J\|_\infty = \|B_J\|_1 = 1$.

The iteration matrix of the corresponding T$_U$-method (\ref{Tagliaf1}) is
$$
{\cal B}_T=\left[ \begin{array}{cccc}
0 & 1 & 0 & 0 \\ 
0 & 0 & 1 & 0 \\
0 & 0 & 1 & 0 \\
0 & 0 & 1 & 0
\end{array} \right]
$$
for which $\|{\cal B}_T\|_\infty = 1$, as assured by Proposition~\ref{norma_inf}.
On the contrary, regarding the $1$-norm, we have that $\|{\cal B}_T\|_1 = 3$.
\label{esempio-2}
\end{example}

\vskip 0.2cm

However, in spite of the foregoing counterexample, an analogous result to
Corollary~\ref{convergence_Bd} does hold also for the {\em strict diagonal dominance
by columns}, i.e., when
\begin{equation}
\|B_{J}\|_1 = \|L+U\|_1 < 1.
\label{str-diag-dom-col}
\end{equation}

Along the proof of the next Proposition~\ref{convergence_Bd_norm1}, we shall use
the notation (\ref{absval}) and the fact that
\begin{equation}
\|\, |B|\, \|_1 = \|B\|_1
\label{abs(B)_norm}
\end{equation}
for all matrices $B$. In fact, this will allow us to apply some results from the
previous Section~\ref{sect_L-matrices}.

\begin{proposition}
Let ${\cal B}_{(d)}=\{B_{1},\ldots ,B_{d}\}$ be a splitting of order $d$ of the
Jacobi iteration matrix $B_{J}$.
If the matrix $A$ is strictly diagonally dominant by columns, i.e., if
(\ref{str-diag-dom-col}) holds, then the corresponding iterative scheme
(\ref{general_method}) is such that $\rho({\cal B}_{(d)}) < 1$ and, therefore,
is convergent.
\label{convergence_Bd_norm1}
\end{proposition}

\begin{proof}
We consider the splitting ${\cal B}_{(d)}'$ defined by (\ref{B_d'}),
for which (\ref{rhoB_d<=rhoB_d'}) holds.

On the other hand, by applying Corollary~\ref{not_slower-3} to $|B_J|$ and
${\cal B}_{(d)}'$, from (\ref{str-diag-dom-col}) and (\ref{abs(B)_norm}) we get
$$
\rho({\cal B}_{(d)}') \le \rho(|B_J|) \le \|B_J\|_1 < 1.
$$
Thus the proof is complete.
\end{proof}

%
%

\section{The AMKS-methods}

\label{IEEE}

In this section we consider the class of AMKS-methods introduced by
Ahmadi et al.~\cite{AMKS} and, for the sake of clarity, we briefly recall their
formulation.

\begin{definition}
We say that the $d$-tuple of matrices
${\cal P}_{(d)}=\{P_{1},\ldots ,P_{d}\}$, $P_{p} \in \re^{n\times n}$, $1\le d\le n$,
is a {\em decomposition of the identity (DoI) in} $\re^{n}$ if the following three
conditions hold:
\begin{itemize}
\item
$P_{p} \ne O$ \ {\rm for all} \ $p=1,\ldots ,d$;
\item
the $P_{p}$'s are {\em logical matrices}, i.e.,
$(P_{p})_{i,j} \in \{0,1\}$ \ {\rm for all} \ $p=1,\ldots ,d$ {\rm and} $i,j=1,\ldots ,n$;
\item
$\sum_{p=1}^d P_{p} = I$, where $I$ is the identity matrix.
\end{itemize}
\label{DoI}
\end{definition}

\vskip 0.2cm

\begin{remark}
Observe that the matrices of the DoI ${\cal P}_{(d)}$ clearly satisfy the condition
$$
P_{p} P_{q} = \delta_{pq}I,
$$
where $\delta_{pq}$ is the Kronecker delta.
\label{remark-DoI}
\end{remark}

\vskip 0.2cm

The corresponding iterative AMKS-method is given by
\begin{equation}
\left\{
\begin{array}{ccl}
x^{(k+1)}_0 & = & \bar{x}^{(k)}, \\
x^{(k+1)}_p & = & P_{p} (B_{J}x^{(k+1)}_{p-1} + c) + (I-P_{p})x^{(k+1)}_{p-1},
\ \ \ p=1,\ldots ,d, \\
\bar{x}^{(k+1)} & = & x^{(k+1)}_{d},
\end{array}
\right.
\label{AMKS}
\end{equation}
and, after elimination of the intermediate vectors $x^{(k+1)}_p$'s, becomes

\begin{equation}
\bar{x}^{(k+1)} =  {\cal P}_{(d)} \bar{x}^{(k)} + {\cal Q}_{(d)} c,
\label{AMKS-2}
\end{equation}
where the iteration matrix
\begin{equation}
{\cal P}_{(d)} = \prod_{p=d}^1 (P_{p}B_{J}+I-P_{p})
\label{AMKS-3}
\end{equation}
is denoted by the same symbol used for the corresponding DoI and
$$
{\cal Q}_{(d)} = \sum_{q=1}^d \left[ \prod_{p=d}^{q+1} (P_{p}B_{J}+I-P_{p}) \right] P_{q}.
$$

We remark that, as anticipated in Section~\ref{introduction}, the matrices $P_{p}B_{J}$
involved in (\ref{AMKS}) and (\ref{AMKS-2}) are made up by some rows of the Jacobi
iteration matrix $B_{J}$.

Now we will show that, substantially, the above class of methods is part of the more
general class of methods introduced by Definitions~\ref{pre-splitting} and \ref{splitting}.
More precisely, we have the following result.

\begin{proposition}
The AMKS-method (\ref{AMKS}) corresponding to the DoI
${\cal P}_{(d)}=\{P_{1},\ldots ,P_{d}\}$ and the method (\ref{general_compact1})
corresponding to the splitting mask ${\cal M}_{(d)}=\{M_{1},\ldots ,M_{d}\}$, where
\begin{equation}
M_{p} = P_{p}E, \ \ \ p=1,\ldots ,d,
\label{AMKS_mask}
\end{equation}
are spectrum-equivalent.
\label{AMKS_equivalence}
\end{proposition}
\begin{proof}
The splitting corresponding to the mask ${\cal M}$ is given by
$$
{\cal B}_{(d)}=\{P_{1}B_J,\ldots ,P_{d}B_J\}.
$$

First we assume that $\lambda \ne 0$ be an eigenvalue of the iteration matrix
${\cal P}_{(d)}$ of the AMKS-method (\ref{AMKS}).
Then there exists $x \in \re^{n}$, $x\ne 0$, such that
\begin{equation}
{\cal P}_{(d)} x = \lambda x.
\label{AMKS_eigen}
\end{equation}

Now we define the $n$-vector
$$
\alpha_1 := x
$$
and the further $d-1$ $n$-vectors
$$
\alpha_p := \lambda(P_{1}x+ \cdots +P_{p-1}x) + P_{p}x+ \cdots +P_{d}x,
\ \ \ \ \ p=2, \ldots, d.
$$

By using the properties of the DoI ${\cal P}_{(d)}$ (see Definition~\ref{DoI} and
Remark~\ref{remark-DoI}), it is easy to verify that the above definitions all
together are equivalent to
\begin{equation}
\lambda P_{p}\alpha_1 = \cdots = \lambda P_{p}\alpha_{p} =
P_{p}\alpha_{p+1} = \cdots = P_{p}\alpha_{d}, \ \ \ \ \ p=1, \ldots, d.
\label{P_alpha}
\end{equation}
Moreover, since $\lambda \ne 0$, it is clear that
\begin{equation}
P_{p}\alpha_1 = \cdots = P_{p}\alpha_{p}, \ \ \ \ \ p=2, \ldots, d.
\label{P_alpha_x}
\end{equation}

Consequently, if we apply the matrices $P_{p}$ to (\ref{AMKS_eigen}), each of them
separately for $p=1, \ldots, d$, still by using the properties of the DoI ${\cal P}_{(d)}$,
long but easy calculations allow us to conclude that
\begin{equation}
P_{p}B_J \alpha_{p} = \lambda P_{p}\alpha_{1}, \ \ \ \ \ p=1, \ldots, d.
\label{P_alpha_B}
\end{equation}

Finally, with $B_{p}=P_{p}B_J$, $p=1, \ldots, d$, it is immediate to see that
(\ref{P_alpha}), (\ref{P_alpha_x}) and (\ref{P_alpha_B}) imply the validity of
(\ref{eigen1}), meaning that $\lambda$ is an eigenvalue of the iteration matrix
${\cal B}_{(d)}$ as well.

Vice versa, now assume that $\lambda \ne 0$ be an eigenvalue of the iteration matrix
${\cal B}_{(d)}$. Then there exists a related eigenvector
$\big[\alpha_1^T, \ldots, \alpha_d^T\big]^T \in \re^{dn}$
with $\alpha_p \ne 0$, $p=1, \ldots, d$, (see Proposition~\ref{nonzero1}) and
the eigenvalue equation (\ref{eigen1}) holds true with $B_{p}=P_{p}B_J$, $p=1, \ldots, d$.

Now, by applying all the matrices $P_{p}$ to (\ref{eigen1}), each of them separately
for $p=1, \ldots, d$, similar calculations as before allow us to prove the validity
of the equalities (\ref{P_alpha}), (\ref{P_alpha_x}) and (\ref{P_alpha_B}).

In turn, by using such sets of equalities in the iteration matrix ${\cal P}_{(d)}$ of
the AMKS-method (\ref{AMKS}) starting from the last factor $P_{1}B_{J}+I-P_{1}$ and
going backward up to the first factor $P_{d}B_{J}+I-P_{d}$, we can easily prove that
$$
{\cal P}_{(d)} \alpha_{1} = \left[ \prod_{p=d}^q (P_{p}B_{J}+I-P_{p}) \right] \alpha_{q},
\ \ \ \ \ q=2, \ldots, d.
$$
In particular, for $q=d$ we get
$$
{\cal P}_{(d)} \alpha_{1} = (P_{d}B_{J}+I-P_{d}) \, \alpha_{d}.
$$
The final step leads us to find that
$$
{\cal P}_{(d)} \alpha_{1} = \lambda \alpha_{1}
$$
and hence, since $\alpha_{1} \ne 0$, to conclude that $\lambda$ is an eigenvalue of
the iteration matrix ${\cal P}_{(d)}$ as well.
\end{proof}

\vskip 0.2cm

A first interesting consequence of the foregoing result is that, substantially, even the
forward Gauss-Seidel method (\ref{f_Gauss-Seidel}) is included in the class of methods 
(\ref{general_compact1}).
In fact, Proposition~5 in \cite{AMKS} shows that it may be obtained as the AMKS-method
corresponding to the DoI ${\cal P}_{(n)}=\{P_{1},\ldots ,P_{n}\}$, where $(P_{p})_{p,p}=1$
for all $p=1,\ldots ,n$ and all the other elements vanish and, consequently, it
corresponds to the splitting mask ${\cal M}_{(n)}=\{M_{n},\ldots ,M_{1}\}$, where
$M_{p}$ is given by (\ref{AMKS_mask}).

Conversely, we remark that, since any splitting method (\ref{general_compact1}) that
is spectrum-equivalent to an AMKS-method (\ref{AMKS}) is a refinement of the Jacobi method
(if different from it), Theorem~\ref{faster} implies, for the case of irreducible
nonnegative Jacobi iteration matrices $B_J$, the result of Proposition~10 in \cite{AMKS}
as a particular instance.

We conclude this section by observing that an analogous class of methods could be
considered, where the Jacobi matrix $B_J$ is split into sets of columns (instead of rows).
This might be done by using splitting masks ${\cal M}_{(d)}=\{M_{1},\ldots ,M_{d}\}$
associated to a DoI ${\cal P}_{(d)}=\{P_{1},\ldots ,P_{d}\}$,
where $M_{p} = EP_{p}$, $p=1,\ldots ,d$.
However, we are not going to propose and study these methods in this paper.

%
%

\section{The T$_U$-method and its refinements}

\label{T-method}

In this section we reconsider the T$_U$-method (\ref{Tagliaf1}), which is associated to
the splitting ${\cal B}_{T_U}=\{U,L\}$, and the method associated to its cyclic shift
${\cal S}({\cal B}_{T_U}) =\{L,U\}$ (see Definition~\ref{def_shift}), which we call
the {\em T$_L$-method}.
Moreover, we introduce and analyze some of their refinements of particular interest.

From now on we shall often refer to the T$_U$ and the T$_L$-method and to their
refinements as the {\em triangular methods}.

Note that, in most cases, a triangular method and a splitting method equivalent to
some AMKS-method are not comparable between them in the partial ordering introduced by
Definitions~\ref{order1} and \ref{order2}.

\begin{lemma}
Let the Jacobi iteration matrix $B_J$ be such that $\rho(B_J)>0$ and
let $L\ne O$ and $U\ne O$.
Then it necessarily holds that $UL\ne O$ and $LU\ne O$, too.
\label{LU_zero}
\end{lemma}
\begin{proof}
We assume by contradiction that $UL=O$.
Then an easy calculation shows that
$$
(L+U)^{2n} = \sum_{h=0}^{2n} L^h U^{2n-h}
$$
and thus, since $L^n=U^n=O$ (recall that $L,U$ are triangular $n\times n$-matrices),
we get $(L+U)^{2n} =O$. Consequently, we obtain the absurde equality
$$
\rho(B_J) = \rho(L+U) = 0.
$$

We arrive at the same conclusion if we assume that $LU=O$.
\end{proof}

\vskip 0.2cm

In view of Definition~\ref{order1} and, in particular, of condition (III),
we can immediately conclude with the following result.
\begin{proposition}
Let the Jacobi iteration matrix $B_J$ be such that $\rho(B_J)>0$ and
let $L\ne O$ and $U\ne O$.
Then the T$_U$-method and the T$_L$-method are more refined than the Jacobi method, i.e.,
the corresponding splittings are such that
$$
{\cal B}_{T_U}=\{U,L\} \succ \{B_J\} \ \ \ {\rm and} \ \ \ {\cal B}_{T_L}=\{L,U\} \succ \{B_J\}.
$$
\label{T_refined}
\end{proposition}

%
%

\subsection{Upper and lower triangular column methods}

\label{LUC methods}

We define the {\em upper triangular column (UTC) methods} by considering the
T$_U$-method associated to ${\cal B}_{T_U}=\{U,L\}$ and by splitting further the upper triangle
$U$ only into some subsets of columns.
This idea is the basis for the possible generation of many methods.


\begin{definition}
Given the Jacobi iteration matrix $B_J=L+U \in \re^{n\times n}$ and $j\in \{2,\ldots ,n\}$, the matrix $U^{(j)}_{c}$ defined by
$$
(U^{(j)}_{c})_{i,h}:=\left\{
\begin{array}{ll}
u_{i,j} & {\rm if} \ h=j, \ i \le j-1, \\
0 & {\rm otherwise},
\end{array}
\right.
$$
is called the {\em $j$-th upper triangular column matrix} of $B_J$.
\label{upper_T_C_M}
\end{definition}

\vskip 0.2cm

With reference to the discussion of Section~\ref{introduction} related to
Definitions~\ref{splitting_mask}, \ref{pre-splitting} and \ref{splitting} and to
Proposition~\ref{proposition_splitting}, for the sake of simplicity (and without loss
of generality) we assume that $L\ne O$ and that all the upper triangular column matrices
$U^{(j)}_{c}\ne O$, $j=2, \ldots ,n$.

The proof of the next result is straightforward.

\begin{proposition}
Given the Jacobi iteration matrix $B_J=L+U \in \re^{n\times n}$, it holds that
$$
U^{(j)}_{c}U^{(k)}_{c}=O \ \ \ \ \ {\rm for \ all} \ \ \ 2\le k \le j \le n.
$$
\label{column-product-zero}
\end{proposition}

\begin{remark}
In the light of Proposition~\ref{prodotto_zero_raggio}, the foregoing result tells us
that, in order to get an effective improvement by splitting the upper triangle $U$,
its columns must be selected in decreasing order with respect to the column index.

In particular, as a limit case, the splitting
$\{U_{c}^{(2)}, \ldots, U_{c}^{(n)},L\}$ determines a method which is equivalent to
the T$_U$-method itself and, therefore, no improvement is obtained at all.
\label{increasing-order}
\end{remark}

\vskip 0.2cm

On the contrary, as a limit case in the opposite direction, we consider the following
{\em full upper triangular column (FUTC)} splitting
$$
{\cal B}_{FUTC} := \{U_{c}^{(n)}, \ldots, U_{c}^{(2)},L\}
$$
and the corresponding iterative scheme
\begin{equation}
\left\{
\begin{array}{ccl}
x^{(k+1)}_{1} & = & U_{c}^{(n)}x^{(k)}_1 + U_{c}^{(n-1)}x^{(k)}_2 + \cdots +
U_{c}^{(2)}x^{(k)}_{n-1} +  Lx^{(k)}_n+ c, \\
x^{(k+1)}_2 & = & U_{c}^{(n)}x^{(k+1)}_1 + U_{c}^{(n-1)}x^{(k)}_2 + \cdots +
U_{c}^{(2)}x^{(k)}_{n-1} +  Lx^{(k)}_n+ c, \\
\cdot & \cdot & \ \ \ \ \ \ \ \ \ \ \ \ \ \ \ \ \ \ \ \ \ \ \ \ \ \ \ \ \ \ \ \ \ \ \ \ \cdots \\
\cdot & \cdot & \ \ \ \ \ \ \ \ \ \ \ \ \ \ \ \ \ \ \ \ \ \ \ \ \ \ \ \ \ \ \ \ \ \ \ \ \cdots \\
\cdot & \cdot & \ \ \ \ \ \ \ \ \ \ \ \ \ \ \ \ \ \ \ \ \ \ \ \ \ \ \ \ \ \ \ \ \ \ \ \ \cdots \\
x^{(k+1)}_n & = & U_{c}^{(n)}x^{(k+1)}_1 + U_{c}^{(n-1)}x^{(k+1)}_2 + \cdots +
U_{c}^{(2)}x^{(k+1)}_{n-1} +  Lx^{(k)}_n+ c.
\end{array}
\right.
\label{FUTC_method}
\end{equation}

\begin{theorem}
The FUTC-method (\ref{FUTC_method}) and the backward Gauss-Seidel method are
spectrum-equivalent.
\label{Gauss-Seidel_equivalence}
\end{theorem}
\begin{proof}
First assume that $\lambda \ne 0$ is an eigenvalue of ${\cal B}_{FUTC}$ and that
$$
\big[\alpha_1^T, \ldots, \alpha_{n-1}^T, \beta^T\big]^T
$$
is a related eigenvector. Then by (\ref{eigen1}) we get
\begin{equation}
\left\{
\begin{array}{ccc}
U_{c}^{(n)}\alpha_1 + U_{c}^{(n-1)}\alpha_2 + \cdots + U_{c}^{(2)}\alpha_{n-1} +
 L\beta & = & \lambda \alpha_1, \\
\lambda U_{c}^{(n)}\alpha_1 + U_{c}^{(n-1)}\alpha_2 + \cdots + U_{c}^{(2)}\alpha_{n-1} +
 L\beta & = & \lambda \alpha_2, \\
\cdots & \cdot & \cdot \\
\cdots & \cdot & \cdot \\
\cdots & \cdot & \cdot \\
\lambda U_{c}^{(n)}\alpha_1 + \lambda U_{c}^{(n-1)}\alpha_2 + \cdots +
U_{c}^{(2)}\alpha_{n-1} +  L\beta & = & \lambda \alpha_{n-1}, \\
\lambda U_{c}^{(n)}\alpha_1 + \lambda U_{c}^{(n-1)}\alpha_2 + \cdots +
\lambda U_{c}^{(2)}\alpha_{n-1} +  L\beta & = & \lambda \beta.
\end{array}
\right.
\label{FUTC_eigen}
\end{equation}

Multiplying the first $n-1$ equalities by $U_{c}^{(n)}, \ldots, U_{c}^{(2)}$, respectively,
and using Proposition~\ref{column-product-zero} yield

\begin{equation}
\left\{
\begin{array}{lcc}
U_{c}^{(n)} L\beta & = & \lambda U_{c}^{(n)} \alpha_1, \\
U_{c}^{(n-1)} \big(\lambda U_{c}^{(n)}\alpha_1 + L\beta \big) & = & 
\lambda U_{c}^{(n-1)}\alpha_2, \\
U_{c}^{(n-2)} \big(\lambda U_{c}^{(n)}\alpha_1 + \lambda U_{c}^{(n-1)}\alpha_2 +
L\beta \big) & = & \lambda U_{c}^{(n-2)}\alpha_3, \\
\ \ \ \ \ \ \ \ \ \ \ \ \ \ \ \ \ \ \ \ \cdots & \cdot & \cdot \\
\ \ \ \ \ \ \ \ \ \ \ \ \ \ \ \ \ \ \ \ \cdots & \cdot & \cdot \\
\ \ \ \ \ \ \ \ \ \ \ \ \ \ \ \ \ \ \ \ \cdots & \cdot & \cdot \\
U_{c}^{(2)} \big(\lambda U_{c}^{(n)}\alpha_1 + \lambda U_{c}^{(n-1)} \alpha_2 +
\cdots + \lambda U_{c}^{(3)}\alpha_{n-2} + L\beta \big) & = &
\lambda U_{c}^{(2)} \alpha_{n-1}.
\end{array}
\right.
\label{FUTC_eigen-2}
\end{equation}
Analogously, multiplying the last equality of (\ref{FUTC_eigen}) by
$U_{c}^{(n)}, \ldots, U_{c}^{(2)}$ in sequence one after the other yields
\begin{equation}
\left\{
\begin{array}{lcc}
U_{c}^{(n)} L\beta & = & \lambda U_{c}^{(n)} \beta, \\
U_{c}^{(n-1)} \big(\lambda U_{c}^{(n)}\alpha_1 + L\beta \big) & = & 
\lambda U_{c}^{(n-1)}\beta, \\
U_{c}^{(n-2)} \big(\lambda U_{c}^{(n)}\alpha_1 + \lambda U_{c}^{(n-1)}\alpha_2 +
L\beta \big) & = & \lambda U_{c}^{(n-2)}\beta, \\
\ \ \ \ \ \ \ \ \ \ \ \ \ \ \ \ \ \ \ \ \cdots & \cdot & \cdot \\
\ \ \ \ \ \ \ \ \ \ \ \ \ \ \ \ \ \ \ \ \cdots & \cdot & \cdot \\
\ \ \ \ \ \ \ \ \ \ \ \ \ \ \ \ \ \ \ \ \cdots & \cdot & \cdot \\
U_{c}^{(2)} \big(\lambda U_{c}^{(n)}\alpha_1 + \lambda U_{c}^{(n-1)} \alpha_2 +
\cdots + \lambda U_{c}^{(3)}\alpha_{n-2} + L\beta \big) & = &
\lambda U_{c}^{(2)} \beta.
\end{array}
\right.
\label{FUTC_eigen-3}
\end{equation}
Therefore, since $\lambda \ne 0$, comparing (\ref{FUTC_eigen-2}) to (\ref{FUTC_eigen-3})
leads to
\begin{equation}
U_{c}^{(n-i+1)} \alpha_{i} = U_{c}^{(n-i+1)} \beta, \ \ \ \ \ i=1, \ldots, n-1.
\label{alpha_i-beta}
\end{equation}
Consequently, substitution into the last equality of (\ref{FUTC_eigen}) yields
\begin{equation}
\lambda \big( U_{c}^{(2)} + \cdots + U_{c}^{(n)} \big) \beta + L \beta =
\lambda U \beta + L \beta = \lambda \beta,
\label{eigen_B_{bGS}}
\end{equation}
which is equivalent to
$$
B_{bGS} \beta = (I-U)^{-1}L \beta = \lambda \beta.
$$
Since $\beta \ne 0$ (see Proposition~\ref{nonzero1}), we can conclude that $\lambda$
is an eigenvalue of the iteration matrix $B_{bGS}$.

Viceversa, let $\lambda \ne 0$ be an eigenvalue of $B_{bGS}$.
Then (\ref{eigen_B_{bGS}}) holds for some $\beta \ne 0$.

With $\sigma$ given by (\ref{sigma-2}), we set
\begin{equation}
\alpha_i := \beta + \sigma \big( U_{c}^{(n-i+1)} + \cdots +
U_{c}^{(2)} \big) \beta, \ \ \ \ \ i=1, \ldots, n-1,
\label{def_alpha_i}
\end{equation}
so that, by multiplying each equality by
$U_{c}^{(n-i+1)}$, respectively, and by using Proposition~\ref{column-product-zero},
we get (\ref{alpha_i-beta}) as well.

Now, it is easy to verify that (\ref{alpha_i-beta}), (\ref{eigen_B_{bGS}}) and
(\ref{def_alpha_i}) all together imply all the equalities of (\ref{FUTC_eigen}).

Thus $\lambda$ is an eigenvalue of ${\cal B}_{FUTC}$ and
$\big[\alpha_1^T, \ldots, \alpha_{n-1}^T, \beta^T\big]^T$ is the corresponding
eigenvector.
\end{proof}

\vskip 0.2cm

\begin{remark}
Since $U$ is strictly upper triangular, $\beta \ne 0$ and
$(1-\lambda)/\lambda~\ne~0$, in (\ref{def_alpha_i}) it must be $\alpha_i \ne 0$
for all $i=1, \ldots, n-1$, in perfect agreement with Proposition~\ref{nonzero1}.
\label{inutile}
\end{remark}

\vskip 0.2cm

The equivalence between the FUTC-method (\ref{FUTC_method})
and the backward Gauss-Seidel method is not limited to having the same spectrum of
the corresponding iteration matrix. Indeed, also the sequences of the respective
approximate solutions are the same.

\begin{theorem}
Let $\{y^{(k)}\}_{k\ge 1}$ be the sequence of approximations produced by the backward
Gauss-Seidel method starting from the initial value $y^{(0)}$, i.e.,
$$
y^{(k+1)} = B_{bGS}y^{(k)}+c_{bGS} = (I-U)^{-1} \big( Ly^{(k)}+c \big), \ \ \ \ \ k\ge 0.
$$
Then the sequence of approximations
$\big\{\big[{x^{(k)}_1}^T, \ldots, {x^{(k)}_{n-1}}^T, {x^{(k)}_n}^T \big]^T \big\}_{k\ge 1}$
produced by the FUTC-method (\ref{FUTC_method})
with initial value
$\big[{x^{(0)}_1}^T, \ldots, {x^{(0)}_{n-1}}^T, {y^{(0)}}^T \big]^T$
is such that
\begin{equation}
x^{(k)}_n = y^{(k)} \ \ \ \ \ {\rm for \ all} \ \ \ k\ge 0
\label{same-sequence}
\end{equation}
independently of the first $n-1$ initial vector-components $x^{(0)}_i$, $i=1, \ldots,n-1$.
\label{Gauss-Seidel_equivalence-2}
\end{theorem}
\begin{proof}
The equality (\ref{same-sequence}) is true for $k=0$ by hypothesis and we assume, by
induction, that it holds for a given $k\ge 1$.

Then, like in the proof of Theorem~\ref{Gauss-Seidel_equivalence}, we multiply
the first $n-1$ equalities of (\ref{FUTC_method}) by $U_{c}^{(n)}, \ldots, U_{c}^{(2)}$,
respectively, and the last again by $U_{c}^{(n)}, \ldots, U_{c}^{(2)}$ in sequence
one after the other, so that to obtain
\begin{equation}
U_{c}^{(n-i+1)} x^{(k+1)}_{i} = U_{c}^{(n-i+1)} x^{(k+1)}_{n}, \ \ \ \ \ i=1, \ldots, n-1.
\label{x_i-y}
\end{equation}
Therefore, by using again the last equality of (\ref{FUTC_method}), we can conclude that
$$
x^{(k+1)}_{n} = \big( U_{c}^{(2)} + \cdots + U_{c}^{(n)} \big) x^{(k+1)}_{n} + Ly^{(k)}+c=
Ux^{(k+1)}_{n} + Ly^{(k)}+c,
$$
which is equivalent to
$$
x^{(k+1)}_{n} = (I-U)^{-1} \big( Ly^{(k)}+c \big) = y^{(k+1)}.
$$
Hence, the induction works and the proof is complete.
\end{proof}

\vskip 0.2cm

Now we briefly treat the symmetric case of the {\em lower triangular column
(LTC) methods} by considering the T$_L$-method associated to ${\cal B}_{T_L}=\{L,U\}$ and
by splitting further the lower triangle $L$ only into some subsets of columns.

Everything is completely analogous to the previous case of the UTC-methods and, hence,
no proof needs to be repetead.

\begin{definition}
Given the Jacobi iteration matrix $B_J=L+U \in \re^{n\times n}$ and
$j\in \{1,\ldots ,n-1\}$, the matrix $L^{(j)}_{c}$ defined by
$$
(L^{(j)}_{c})_{i,h}:=\left\{
\begin{array}{ll}
l_{i,j} & {\rm if} \ h=j, \ i \ge j+1, \\
0 & {\rm otherwise},
\end{array}
\right.
$$
is called the {\em $j$-th lower triangular column matrix} of $B_J$.
\label{lower_T_C_M}
\end{definition}

\vskip 0.2cm

Again, for the sake of simplicity (and without loss of generality) we assume that
$U\ne O$ and that all the lower triangular column matrices
$L^{(j)}_{c}\ne O$, $j=1, \ldots ,n-1$.
\begin{proposition}
Given the Jacobi iteration matrix $B_J=L+U \in \re^{n\times n}$, it holds that
$$
L^{(j)}_{c}L^{(k)}_{c}=O \ \ \ \ \ {\rm for \ all} \ \ \ 1\le j \le k \le n-1.
$$
\label{column-product-zero-2}
\end{proposition}

\begin{remark}
In the light of Proposition~\ref{prodotto_zero_raggio}, the foregoing result tells us
that, in order to get an effective improvement by splitting the lower triangle $L$,
its columns must be selected in increasing order with respect to the column index.

In particular, as a limit case, the splitting
$\{L_{c}^{(n-1)}, \ldots, L_{c}^{(1)},U\}$ determines a method which is equivalent to
the T$_L$-method itself and, therefore, no improvement is obtained at all.
\label{increasing-order-2}
\end{remark}

\vskip 0.2cm

On the contrary, as a limit case in the opposite direction, it is easy to verify that
the symmetric {\em full lower triangular column (FLTC)} splitting
$$
{\cal B}_{FLTC} := \{L_{c}^{(1)}, \ldots, L_{c}^{(n-1)},U\}
$$
gives rise to an iterative scheme (the {\em FLTC-method}) which is spectrum-equivalent
to the forward Gauss-Seidel method and that the analogous results to
Theorems~\ref{Gauss-Seidel_equivalence} and \ref{Gauss-Seidel_equivalence-2} hold.

It is worth remarking that, in the light of the results of this section, any
UTC-method lies between of the T$_U$- and the backward Gauss-Seidel method with
respect to the partial order relation of refinement ''$\succ$'' introduced on the corresponding
splittings by Definitions~\ref{order1} and \ref{order2}.
Consequently, a generic UTC-method often (but not always) performs faster then the
T$_U$-method and slower than the backward Gauss-Seidel method.
In particular, such a behaviour is assured if the methods are applied to linear
systems characterized by an irreducible nonnegative Jacobi iteration matrix
(see Section~\ref{sect_L-matrices}).

%
%

\subsection{Upper and lower triangular row methods}

\label{LUR methods}
Similarly to what we did in the previous Section~\ref{LUC methods}, we can define the
{\em upper triangular row (UTR) methods} by considering the T$_U$-method associated
to ${\cal B}_{T_U}=\{U,L\}$ and by splitting further the upper triangle $U$ only into
some subsets of rows.

\begin{definition}
Given the Jacobi iteration matrix $B_J=L+U \in \re^{n\times n}$ and
$i\in \{1,\ldots ,n-1\}$, the matrix $U^{(i)}_{r}$ defined by
$$
(U^{(i)}_{r})_{h,j}:=\left\{
\begin{array}{ll}
u_{i,j} & {\rm if} \ h=i, \ j \ge i+1, \\
0 & {\rm otherwise},
\end{array}
\right.
$$
is called the {\em $i$-th upper triangular row matrix} of $B_J$.
\label{upper_T_R_M}
\end{definition}

\vskip 0.2cm

Again, for the sake of simplicity (and without loss of generality) we assume that $L\ne O$
and that all the upper triangular row matrices $U^{(i)}_{r}\ne O$, $i=1, \ldots ,n-1$.

The proofs of all the next results are completely analogous to those of
the corresponding ones in the previous Section~\ref{LUC methods}.

\begin{proposition}
Given the Jacobi iteration matrix $B_J=L+U \in \re^{n\times n}$, it holds that
$$
U^{(i)}_{r}U^{(l)}_{r}=O \ \ \ \ \ {\rm for \ all} \ \ \ 1\le l \le i \le n-1.
$$
\label{column-product-zero-3}
\end{proposition}

\begin{remark}
In the light of Proposition~\ref{prodotto_zero_raggio}, the foregoing result tells us
that, in order to get an effective improvement by splitting the upper triangle $U$,
its rows must be selected in decreasing order with respect to the row index.

In particular, as a limit case, the splitting
$\{U_{r}^{(1)}, \ldots, U_{r}^{(n-1)},L\}$ determines a method which is
spectrum-equivalent to the T$_U$-method itself and, therefore, no improvement is
obtained at all.
\label{increasing-order-3}
\end{remark}

\vskip 0.2cm

On the contrary, as a limit case in the opposite direction, we consider the following
{\em full upper triangular row (FUTR)} splitting
$$
{\cal B}_{FUTR} := \{U_{r}^{(n-1)}, \ldots, U_{r}^{(1)},L\}
$$
and the corresponding iterative scheme
\begin{equation}
\left\{
\begin{array}{ccl}
x^{(k+1)}_{1} & = & U_{r}^{(n-1)}x^{(k)}_1 + U_{r}^{(n-2)}x^{(k)}_2 + \cdots +
U_{r}^{(1)}x^{(k)}_{n-1} +  Lx^{(k)}_n+ c, \\
x^{(k+1)}_2 & = & U_{r}^{(n-1)}x^{(k+1)}_1 + U_{r}^{(n-2)}x^{(k)}_2 + \cdots +
U_{r}^{(1)}x^{(k)}_{n-1} +  Lx^{(k)}_n+ c, \\
\cdot & \cdot & \ \ \ \ \ \ \ \ \ \ \ \ \ \ \ \ \ \ \ \ \ \ \ \ \ \ \ \ \ \ \ \ \ \ \ \ \cdots \\
\cdot & \cdot & \ \ \ \ \ \ \ \ \ \ \ \ \ \ \ \ \ \ \ \ \ \ \ \ \ \ \ \ \ \ \ \ \ \ \ \ \cdots \\
\cdot & \cdot & \ \ \ \ \ \ \ \ \ \ \ \ \ \ \ \ \ \ \ \ \ \ \ \ \ \ \ \ \ \ \ \ \ \ \ \ \cdots \\
x^{(k+1)}_n & = & U_{r}^{(n-1)}x^{(k+1)}_1 + U_{r}^{(n-2)}x^{(k+1)}_2 + \cdots +
U_{r}^{(1)}x^{(k+1)}_{n-1} +  Lx^{(k)}_n+ c.
\end{array}
\right.
\label{FUTR_method}
\end{equation}

Then we have the following equivalence result with the backward Gauss-Seidel method.

\begin{theorem}
The FUTR-method (\ref{FUTR_method}) and the backward Gauss-Seidel method are
spectrum-equivalent.
\label{Gauss-Seidel_equivalence-3}
\end{theorem}

\vskip 0.2cm

Also in this case, the equivalence between the FUTR-method (\ref{FUTR_method})
and the backward Gauss-Seidel method is not limited to having the same spectrum of
the corresponding iteration matrix.

\begin{theorem}
Let $\{y^{(k)}\}_{k\ge 1}$ be the sequence of approximations produced by the backward
Gauss-Seidel method starting from the initial value $y^{(0)}$, i.e.,
$$
y^{(k+1)} = B_{bGS}y^{(k)}+c_{bGS} = (I-U)^{-1} \big( Ly^{(k)}+c \big), \ \ \ \ \ k\ge 0.
$$
Then the sequence of approximations
$\big\{\big[{x^{(k)}_1}^T, \ldots, {x^{(k)}_{n-1}}^T, {x^{(k)}_n}^T \big]^T \big\}_{k\ge 0}$
produced by the FUTR-method (\ref{FUTR_method})
with initial value
$\big[{x^{(0)}_1}^T, \ldots, {x^{(0)}_{n-1}}^T, {y^{(0)}}^T \big]^T$
is such that
\begin{equation}
x^{(k)}_n = y^{(k)} \ \ \ \ \ {\rm for \ all} \ \ \ k\ge 0
\label{same-sequence-2}
\end{equation}
independently of the first $n-1$ initial vector-components $x^{(0)}_i$, $i=1, \ldots,n-1$.
\label{Gauss-Seidel_equivalence-4}
\end{theorem}

\vskip 0.2cm

We conclude this section with the symmetric case of the {\em lower triangular row
(LTR) method} by considering the T$_L$-method associated to ${\cal B}_{T_L}=\{L,U\}$ and
by splitting further the lower triangle $L$ only into some subsets of rows.

\begin{definition}
Given the Jacobi iteration matrix $B_J=L+U \in \re^{n\times n}$ and
$i\in \{2,\ldots ,n\}$, the matrix $L^{(i)}_{r}$ defined by
$$
(L^{(i)}_{r})_{h,j}:=\left\{
\begin{array}{ll}
l_{i,j} & {\rm if} \ h=i, \ j \le i-1, \\
0 & {\rm otherwise},
\end{array}
\right.
$$
is called the {\em $i$-th lower triangular row matrix} of $B_J$.
\label{lower_T_R_M}
\end{definition}

\vskip 0.2cm

Again, for the sake of simplicity (and without loss of generality) we assume that
$U\ne O$ and that all the lower triangular row matrices
$L^{(i)}_{r}\ne O$, $i=2, \ldots ,n$.
\begin{proposition}
Given the Jacobi iteration matrix $B_J=L+U \in \re^{n\times n}$, it holds that
$$
L^{(i)}_{r}L^{(l)}_{r}=O \ \ \ \ \ {\rm for \ all} \ \ \ 2\le i \le l \le n.
$$
\label{column-product-zero-4}
\end{proposition}

\begin{remark}
In the light of Proposition~\ref{prodotto_zero_raggio}, the foregoing result tells us
that, in order to get an effective improvement by splitting the lower triangle $L$,
its columns must be selected in increasing order with respect to the row index.

In particular, as a limit case, the splitting
$\{L_{r}^{(n)}, \ldots, L_{r}^{(2)},U\}$ determines a method which is spectrum-equivalent
to the T$_L$-method itself and, therefore, no improvement is obtained at all.
\label{increasing-order-4}
\end{remark}

\vskip 0.2cm

On the contrary, as a limit case in the opposite direction, it is easy to verify that
the symmetric {\em full lower triangular row (FLTR)} splitting
$$
{\cal B}_{FLTR} := \{L_{r}^{(2)}, \ldots, L_{r}^{(n)},U\}
$$
gives rise to an iterative scheme (the {\em FLTR-method}) which is spectrum-equivalent
to the forward Gauss-Seidel method and that the analogous results to
Theorems~\ref{Gauss-Seidel_equivalence-3} and \ref{Gauss-Seidel_equivalence-4} hold.

%
%

\subsection{The triangular column and row methods}

\label{TC-TR methods}

We define the {\em triangular column (TC) methods} by considering the T$_L$-method
associated to ${\cal B}_{T_L}=\{L,U\}$ and by splitting further both the lower and the
upper triangles $L$ and $U$, the former before the latter separately, into some subsets
of columns.

Analogously, we define the {\em triangular row (TR) methods} by splitting
further both the lower and the upper triangles $L$ and $U$, the former before the
latter separately, into some subsets of rows.

Now we propose one of the simplest possible choices of TC- and TR-methods, but already
somehow significant as we shall see in Section~\ref{numerics} by means of some numerical
experiments.

In order to simplify the subsequent notation, we set
\begin{equation}
\nu := \left\{
\begin{array}{ll}
\frac{n}{2}-1 & {\rm if} \ n \ {\rm is \ even}, \\
\frac{n-1}{2} & {\rm if} \ n \ {\rm is \ odd}.
\end{array}
\right.
\label{nu}
\end{equation}
Then, with reference to Definitions~\ref{upper_T_C_M} and \ref{lower_T_C_M}, we consider
the splitting
$$
{\cal B}_{TC(2,2)} := \{L_{c}^{(1,\nu)}, L_{c}^{(\nu+1,n-1)},
U_{c}^{(n-\nu+1,n)}, U_{c}^{(2,n-\nu)}\},
$$
where
$$
\begin{array}{l}
L_{c}^{(1,\nu)}:=L_{c}^{(1)}+ \cdots +L_{c}^{(\nu)}, \\
L_{c}^{(\nu+1,n-1)}:=L_{c}^{(\nu+1)}+ \cdots +L_{c}^{(n-1)}, \\
U_{c}^{(2,n-\nu)}:=U_{c}^{(2)}+ \cdots +U_{c}^{(n-\nu)}, \\
U_{c}^{(n-\nu+1,n)}:=U_{c}^{(n-\nu+1)}+ \cdots +U_{c}^{(n)},
\end{array}
$$
which gives rise to the iterative scheme
\begin{equation}
\left\{
\begin{array}{l}
x^{(k+1)}_{1} = L_{c}^{(1,\nu)}x^{(k)}_1 + L_{c}^{(\nu+1,n-1)}x^{(k)}_2 +
U_{c}^{(n-\nu+1,n)}x^{(k)}_{3} + U_{c}^{(2,n-\nu)}x^{(k)}_{4} + c, \\
x^{(k+1)}_{2} = L_{c}^{(1,\nu)}x^{(k+1)}_1 + L_{c}^{(\nu+1,n-1)}x^{(k)}_2 +
U_{c}^{(n-\nu+1,n)}x^{(k)}_{3} + U_{c}^{(2,n-\nu)}x^{(k)}_{4} + c, \\
x^{(k+1)}_{3} = L_{c}^{(1,\nu)}x^{(k+1)}_1 + L_{c}^{(\nu+1,n-1)}x^{(k+1)}_2 +
U_{c}^{(n-\nu+1,n)}x^{(k)}_{3} + U_{c}^{(2,n-\nu)}x^{(k)}_{4} + c, \\
x^{(k+1)}_{4} = L_{c}^{(1,\nu)}x^{(k+1)}_1 + L_{c}^{(\nu+1,n-1)}x^{(k+1)}_2 +
U_{c}^{(n-\nu+1,n)}x^{(k+1)}_{3} + U_{c}^{(2,n-\nu)}x^{(k)}_{4} + c.
\end{array}
\right.
\label{TC(2,2)_method}
\end{equation}

Note that, in view of Propositions~\ref{column-product-zero} and \ref{column-product-zero-2},
we get
$$
L_{c}^{(1,\nu)} L_{c}^{(\nu+1,n-1)} = U_{c}^{(n-\nu+1,n)} U_{c}^{(2,n-\nu)} = O
$$
and, therefore, the order of the elements in the splitting ${\cal B}_{TC(2,2)}$ can not
be changed.

The specification $(2,2)$ in ${\cal B}_{TC(2,2)}$ obviously indicates the
number of divisions of $L$ and $U$.

Analogously, with reference to Definitions~\ref{upper_T_R_M} and \ref{lower_T_R_M},
we consider the splitting
$$
{\cal B}_{TR(2,2)} := \{L_{r}^{(2,n-\nu)}, L_{r}^{(n-\nu+1,n)},
U_{r}^{(\nu+1,n-1)}, U_{r}^{(1,\nu)}\},
$$
where
$$
\begin{array}{l}
L_{r}^{(2,n-\nu)}:=L_{r}^{(2)}+ \cdots +L_{r}^{(n-\nu)}, \\
L_{r}^{(n-\nu+1,n)}:=L_{r}^{(n-\nu+1)}+ \cdots +L_{r}^{(n)}, \\
U_{r}^{(1,\nu)}:=U_{r}^{(1)}+ \cdots +U_{r}^{(\nu)}, \\
U_{r}^{(\nu+1,n-1)}:=U_{r}^{(\nu+1)}+ \cdots +U_{r}^{(n-1)},
\end{array}
$$
which gives rise to the iterative scheme
\begin{equation}
\left\{
\begin{array}{l}
x^{(k+1)}_{1} = L_{r}^{(2,n-\nu)}x^{(k)}_1 + L_{r}^{(n-\nu+1,n)}x^{(k)}_2 +
U_{r}^{(\nu+1,n-1)}x^{(k)}_{3} + U_{r}^{(1,\nu)}x^{(k)}_{4} + c, \\
x^{(k+1)}_{2} = L_{r}^{(2,n-\nu)}x^{(k+1)}_1 + L_{r}^{(n-\nu+1,n)}x^{(k)}_2 +
U_{r}^{(\nu+1,n-1)}x^{(k)}_{3} + U_{r}^{(1,\nu)}x^{(k)}_{4} + c, \\
x^{(k+1)}_{3} = L_{r}^{(2,n-\nu)}x^{(k+1)}_1 + L_{r}^{(n-\nu+1,n)}x^{(k+1)}_2 +
U_{r}^{(\nu+1,n-1)}x^{(k)}_{3} + U_{r}^{(1,\nu)}x^{(k)}_{4} + c, \\
x^{(k+1)}_{4} = L_{r}^{(2,n-\nu)}x^{(k+1)}_1 + L_{r}^{(n-\nu+1,n)}x^{(k+1)}_2 +
U_{r}^{(\nu+1,n-1)}x^{(k+1)}_{3} + U_{r}^{(1,\nu)}x^{(k)}_{4} + c.
\end{array}
\right.
\label{TR(2,2)_method}
\end{equation}

Like before, Propositions~\ref{column-product-zero-3} and \ref{column-product-zero-4}
imply
$$
L_{r}^{(2,n-\nu)} L_{r}^{(n-\nu+1,n)} = U_{r}^{(\nu+1,n-1)} U_{r}^{(1,\nu)} = O
$$
and, therefore, the order of the elements in the splitting ${\cal B}_{TR(2,2)}$ can not
be changed either.

In the light of their spectrum-equivalence with the full triangular upper/lower
column/row methods (proved in the foregoing Sections~\ref{LUC methods} and
\ref{LUR methods}), the natural idea to get a simultaneous improvement of both the
backward and forward Gauss-Seidel methods is that of finding a common refinement of their
corresponding splittings, hopefully maximal in the partial order relation
''$\succ$'' introduced by Definitions~\ref{order1} and \ref{order2}.

We start by observing that, if we split an $i$-th lower triangular column (row) matrix
$L_{c}^{(i)}$ ($L_{r}^{(i)}$) into a pair $\{L', L'' \}$ in whatever way, it always holds
that the product
$$
L' L'' = O.
$$
Therefore, the splittings of the lower triangular matrix $L$ employed by the FLTC- and
the FLTR-methods are maximal in the partial order relation ''$\succ$''.

Analogously, if we split an $i$-th lower triangular column (row) matrix
$U_{c}^{(i)}$ ($U_{r}^{(i)}$) into a pair $\{U', U'' \}$ in whatever way, it always holds
that the product
$$
U' U'' = O
$$
and, hence, also the splittings of the upper triangular matrix $U$
employed by the FUTC- and the FUTR-methods are maximal in the partial order relation
''$\succ$''.

In conclusion, it is easy to see that the following {\em full triangular column (FTC)}
splitting
$$
{\cal B}_{FTC} := \{L_{c}^{(1)}, \ldots, L_{c}^{(n-1)},U_{c}^{(n)}, \ldots, U_{c}^{(2)}\}
$$ 
is a maximal refinement of the full lower triangular column splitting ${\cal B}_{FLTC}$.

The corresponding iterative scheme is
\begin{equation}
\left\{
\begin{array}{ccl}
x^{(k+1)}_{1} & = & L_{c}^{(1)}x^{(k)}_1 + L_{c}^{(2)}x^{(k)}_2 + \cdots +
L_{c}^{(n-2)}x^{(k)}_{n-2} + L_{c}^{(n-1)}x^{(k)}_{n-1} + \\
 & & + U_{c}^{(n)}x^{(k)}_{n} + U_{c}^{(n-1)}x^{(k)}_{n+1} + \cdots +
U_{c}^{(3)}x^{(k)}_{2n-3} + U_{c}^{(2)}x^{(k)}_{2n-2} + c, \\
x^{(k+1)}_2 & = & L_{c}^{(1)}x^{(k+1)}_1 + L_{c}^{(2)}x^{(k)}_2 + \cdots +
L_{c}^{(n-2)}x^{(k)}_{n-2} + L_{c}^{(n-1)}x^{(k)}_{n-1} + \\
 & & + U_{c}^{(n)}x^{(k)}_{n} + U_{c}^{(n-1)}x^{(k)}_{n+1} + \cdots +
U_{c}^{(3)}x^{(k)}_{2n-3} + U_{c}^{(2)}x^{(k)}_{2n-2} + c, \\
\cdot & \cdot & \ \ \ \ \ \ \ \ \ \ \ \ \ \ \ \ \ \ \ \ \ \ \ \ \ \ \ \ \ \ \cdots \\
\cdot & \cdot & \ \ \ \ \ \ \ \ \ \ \ \ \ \ \ \ \ \ \ \ \ \ \ \ \ \ \ \ \ \ \cdots \\
\cdot & \cdot & \ \ \ \ \ \ \ \ \ \ \ \ \ \ \ \ \ \ \ \ \ \ \ \ \ \ \ \ \ \ \cdots \\
x^{(k+1)}_{n-1} & = & L_{c}^{(1)}x^{(k+1)}_1 + L_{c}^{(2)}x^{(k+1)}_2 + \cdots +
L_{c}^{(n-2)}x^{(k+1)}_{n-2} + L_{c}^{(n-1)}x^{(k)}_{n-1} + \\
 & & + U_{c}^{(n)}x^{(k)}_{n} + U_{c}^{(n-1)}x^{(k)}_{n+1} + \cdots +
U_{c}^{(3)}x^{(k)}_{2n-3} + U_{c}^{(2)}x^{(k)}_{2n-2} + c, \\
x^{(k+1)}_{n} & = & L_{c}^{(1)}x^{(k+1)}_1 + L_{c}^{(2)}x^{(k+1)}_2 + \cdots +
L_{c}^{(n-2)}x^{(k+1)}_{n-2} + L_{c}^{(n-1)}x^{(k+1)}_{n-1} + \\
 & & + U_{c}^{(n)}x^{(k)}_{n} + U_{c}^{(n-1)}x^{(k)}_{n+1} + \cdots +
U_{c}^{(3)}x^{(k)}_{2n-3} + U_{c}^{(2)}x^{(k)}_{2n-2} + c, \\
x^{(k+1)}_{n+1} & = & L_{c}^{(1)}x^{(k+1)}_1 + L_{c}^{(2)}x^{(k+1)}_2 + \cdots +
L_{c}^{(n-2)}x^{(k+1)}_{n-2} + L_{c}^{(n-1)}x^{(k+1)}_{n-1} + \\
 & & + U_{c}^{(n)}x^{(k+1)}_{n} + U_{c}^{(n-1)}x^{(k)}_{n+1} + \cdots +
U_{c}^{(3)}x^{(k)}_{2n-3} + U_{c}^{(2)}x^{(k)}_{2n-2} + c, \\
\cdot & \cdot & \ \ \ \ \ \ \ \ \ \ \ \ \ \ \ \ \ \ \ \ \ \ \ \ \ \ \ \ \ \ \cdots \\
\cdot & \cdot & \ \ \ \ \ \ \ \ \ \ \ \ \ \ \ \ \ \ \ \ \ \ \ \ \ \ \ \ \ \ \cdots \\
\cdot & \cdot & \ \ \ \ \ \ \ \ \ \ \ \ \ \ \ \ \ \ \ \ \ \ \ \ \ \ \ \ \ \ \cdots \\
x^{(k+1)}_{2n-2} & = & L_{c}^{(1)}x^{(k+1)}_1 + L_{c}^{(2)}x^{(k+1)}_2 + \cdots +
L_{c}^{(n-2)}x^{(k+1)}_{n-2} + L_{c}^{(n-1)}x^{(k+1)}_{n-1} + \\
 & & + U_{c}^{(n)}x^{(k+1)}_{n} + U_{c}^{(n-1)}x^{(k+1)}_{n+1} + \cdots +
U_{c}^{(3)}x^{(k+1)}_{2n-3} + U_{c}^{(2)}x^{(k)}_{2n-2} + c.
\end{array}
\right.
\label{FTC_method}
\end{equation}

Furthermore, since the splitting  $\{L, U_{c}^{(n)}, \ldots, U_{c}^{(2)}\}$ is clearly
a cyclic permutation of the full upper triangular column splitting ${\cal B}_{FUTC}$,
the splitting ${\cal B}_{FTC}$ is also a maximal refinement of it modulo a cyclic
permutation which, substantially, does not modify the corresponding method
(see Corollary~\ref{cyclic_spectrum} and the subsequent discussion).

Similarly, the following {\em full triangular row (FTR)} splitting
$$
{\cal B}_{FTR} := \{L_{r}^{(2)}, \ldots, L_{r}^{(n)},U_{r}^{(n-1)}, \ldots, U_{r}^{(1)}\}
$$
is a maximal refinement of the full lower triangular row splitting ${\cal B}_{FLTR}$
and, modulo a cyclic permutation, of the full upper triangular row splitting
${\cal B}_{FUTR}$ as well.

The corresponding iterative scheme is
\begin{equation}
\left\{
\begin{array}{ccl}
x^{(k+1)}_{1} & = & L_{r}^{(2)}x^{(k)}_1 + L_{r}^{(3)}x^{(k)}_2 + \cdots +
L_{r}^{(n-1)}x^{(k)}_{n-2} + L_{r}^{(n)}x^{(k)}_{n-1} + \\
 & & + U_{r}^{(n-1)}x^{(k)}_{n} + U_{r}^{(n-2)}x^{(k)}_{n+1} + \cdots +
U_{r}^{(2)}x^{(k)}_{2n-3} + U_{r}^{(1)}x^{(k)}_{2n-2} + c, \\
x^{(k+1)}_2 & = & L_{r}^{(1)}x^{(k+1)}_1 + L_{r}^{(2)}x^{(k)}_2 + \cdots +
L_{r}^{(n-2)}x^{(k)}_{n-2} + L_{r}^{(n-1)}x^{(k)}_{n-1} + \\
 & & + U_{r}^{(n-1)}x^{(k)}_{n} + U_{r}^{(n-2)}x^{(k)}_{n+1} + \cdots +
U_{r}^{(2)}x^{(k)}_{2n-3} + U_{r}^{(1)}x^{(k)}_{2n-2} + c, \\
\cdot & \cdot & \ \ \ \ \ \ \ \ \ \ \ \ \ \ \ \ \ \ \ \ \ \ \ \ \ \ \ \ \ \ \cdots \\
\cdot & \cdot & \ \ \ \ \ \ \ \ \ \ \ \ \ \ \ \ \ \ \ \ \ \ \ \ \ \ \ \ \ \ \cdots \\
\cdot & \cdot & \ \ \ \ \ \ \ \ \ \ \ \ \ \ \ \ \ \ \ \ \ \ \ \ \ \ \ \ \ \ \cdots \\
x^{(k+1)}_{n-1} & = & L_{r}^{(1)}x^{(k+1)}_1 + L_{r}^{(2)}x^{(k+1)}_2 + \cdots +
L_{r}^{(n-2)}x^{(k+1)}_{n-2} + L_{r}^{(n-1)}x^{(k)}_{n-1} + \\
 & & + U_{r}^{(n-1)}x^{(k)}_{n} + U_{r}^{(n-2)}x^{(k)}_{n+1} + \cdots +
U_{r}^{(2)}x^{(k)}_{2n-3} + U_{r}^{(1)}x^{(k)}_{2n-2} + c, \\
x^{(k+1)}_{n} & = & L_{r}^{(1)}x^{(k+1)}_1 + L_{r}^{(2)}x^{(k+1)}_2 + \cdots +
L_{r}^{(n-2)}x^{(k+1)}_{n-2} + L_{r}^{(n-1)}x^{(k+1)}_{n-1} + \\
 & & + U_{r}^{(n-1)}x^{(k)}_{n} + U_{r}^{(n-2)}x^{(k)}_{n+1} + \cdots +
U_{r}^{(2)}x^{(k)}_{2n-3} + U_{r}^{(1)}x^{(k)}_{2n-2} + c, \\
x^{(k+1)}_{n+1} & = & L_{r}^{(1)}x^{(k+1)}_1 + L_{r}^{(2)}x^{(k+1)}_2 + \cdots +
L_{r}^{(n-2)}x^{(k+1)}_{n-2} + L_{r}^{(n-1)}x^{(k+1)}_{n-1} + \\
 & & + U_{r}^{(n-1)}x^{(k+1)}_{n} + U_{r}^{(n-2)}x^{(k)}_{n+1} + \cdots +
U_{r}^{(2)}x^{(k)}_{2n-3} + U_{r}^{(1)}x^{(k)}_{2n-2} + c, \\
\cdot & \cdot & \ \ \ \ \ \ \ \ \ \ \ \ \ \ \ \ \ \ \ \ \ \ \ \ \ \ \ \ \ \ \cdots \\
\cdot & \cdot & \ \ \ \ \ \ \ \ \ \ \ \ \ \ \ \ \ \ \ \ \ \ \ \ \ \ \ \ \ \ \cdots \\
\cdot & \cdot & \ \ \ \ \ \ \ \ \ \ \ \ \ \ \ \ \ \ \ \ \ \ \ \ \ \ \ \ \ \ \cdots \\
x^{(k+1)}_{2n-2} & = & L_{r}^{(1)}x^{(k+1)}_1 + L_{r}^{(2)}x^{(k+1)}_2 + \cdots +
L_{r}^{(n-2)}x^{(k+1)}_{n-2} + L_{r}^{(n-1)}x^{(k+1)}_{n-1} + \\
 & & + U_{r}^{(n-1)}x^{(k+1)}_{n} + U_{r}^{(n-2)}x^{(k+1)}_{n+1} + \cdots +
U_{r}^{(2)}x^{(k+1)}_{2n-3} + U_{r}^{(1)}x^{(k)}_{2n-2} + c.
\end{array}
\right.
\label{FTR_method}
\end{equation}

Indeed, now we shall see that the corresponding {\em FTC-method} and {\em FTR-method}
are themselves substantially equivalent to each other and, also, to the symmetric
Gauss-Seidel method (see (\ref{sym_Gauss-Seidel}) for the general case of a nonsingular
diagonal matrix $D$).
\begin{theorem}
The FTC-method, the FTR-method and the symmetric Gauss-Seidel method are all
spectrum-equivalent.
\label{FTC-FTR-sGS_equivalence}
\end{theorem}
\begin{proof}
Writing (\ref{eigen1}) for the splitting ${\cal B}_{FTC}$ (with respect to the
eigenvalue $\lambda$ and its related eigenvector
$[\alpha_{1}^T, \ldots, \alpha_{n-1}^T,\beta_{1}^T, \ldots, \beta_{n-1}^T]^T$) leads to
\begin{equation}
\left\{
\begin{array}{ccc}
L_{c}^{(1)}\alpha_1 + L_{c}^{(2)}\alpha_2 + \cdots +
L_{c}^{(n-2)}\alpha_{n-2} + L_{c}^{(n-1)}\alpha_{n-1} + \\
+ U_{c}^{(n)}\beta_1 + U_{c}^{(n-1)}\beta_2 + \cdots +
U_{c}^{(3)}\beta_{n-2} + U_{c}^{(2)}\beta_{n-1} & = & \lambda \alpha_{1}, \\

\lambda L_{c}^{(1)}\alpha_1 + L_{c}^{(2)}\alpha_2 + \cdots +
L_{c}^{(n-2)}\alpha_{n-2} + L_{c}^{(n-1)}\alpha_{n-1} + \\
+ U_{c}^{(n)}\beta_1 + U_{c}^{(n-1)}\beta_2 + \cdots +
U_{c}^{(3)}\beta_{n-2} + U_{c}^{(2)}\beta_{n-1} & = & \lambda \alpha_{2}, \\
\cdots & \cdot & \cdot \\
\cdots & \cdot & \cdot \\
\cdots & \cdot & \cdot \\
\lambda L_{c}^{(1)}\alpha_1 + \lambda L_{c}^{(2)}\alpha_2 + \cdots +
\lambda L_{c}^{(n-2)}\alpha_{n-2} + L_{c}^{(n-1)}\alpha_{n-1} + \\
+ U_{c}^{(n)}\beta_1 + U_{c}^{(n-1)}\beta_2 + \cdots +
U_{c}^{(3)}\beta_{n-2} + U_{c}^{(2)}\beta_{n-1} & = & \lambda \alpha_{n-1}, \\
\lambda L_{c}^{(1)}\alpha_1 + \lambda L_{c}^{(2)}\alpha_2 + \cdots +
\lambda L_{c}^{(n-2)}\alpha_{n-2} + \lambda L_{c}^{(n-1)}\alpha_{n-1} + \\
+ U_{c}^{(n)}\beta_1 + U_{c}^{(n-1)}\beta_2 + \cdots +
U_{c}^{(3)}\beta_{n-2} + U_{c}^{(2)}\beta_{n-1} & = & \lambda \beta_{1}, \\
\lambda L_{c}^{(1)}\alpha_1 + \lambda L_{c}^{(2)}\alpha_2 + \cdots +
\lambda L_{c}^{(n-2)}\alpha_{n-2} + \lambda L_{c}^{(n-1)}\alpha_{n-1} + \\
+ \lambda U_{c}^{(n)}\beta_1 + U_{c}^{(n-1)}\beta_2 + \cdots +
U_{c}^{(3)}\beta_{n-2} + U_{c}^{(2)}\beta_{n-1} & = & \lambda \beta_{2}, \\
\cdots & \cdot & \cdot \\
\cdots & \cdot & \cdot \\
\cdots & \cdot & \cdot \\
\lambda L_{c}^{(1)}\alpha_1 + \lambda L_{c}^{(2)}\alpha_2 + \cdots +
\lambda L_{c}^{(n-2)}\alpha_{n-2} + \lambda L_{c}^{(n-1)}\alpha_{n-1} + \\
+ \lambda U_{c}^{(n)}\beta_1 + \lambda U_{c}^{(n-1)}\beta_2 + \cdots +
\lambda U_{c}^{(3)}\beta_{n-2} + U_{c}^{(2)}\beta_{n-1} & = & \lambda \beta_{n-1}.
\end{array}
\right.
\label{FTC_eigen}
\end{equation}

Since $\lambda \ne 0$, we can consider the constant $\sigma$ given by (\ref{sigma-2}).
Therefore, by taking the difference of each pair of consecutive equalities in 
(\ref{FTC_eigen}) we obtain
\begin{equation}
\alpha_{p+1} = \big(I-\sigma L_{c}^{(p)}\big)\alpha_p, \ \ \ p=1,\ldots ,n-2,
\label{differences-3}
\end{equation}
\begin{equation}
\beta_{1} = \big(I-\sigma L_{c}^{(n-1)}\big)\alpha_{n-1},
\label{differences-4}
\end{equation}
\begin{equation}
\beta_{p+1} = \big(I-\sigma U_{c}^{(n+1-p)}\big)\beta_p, \ \ \ p=1,\ldots ,n-2.
\label{differences-5}
\end{equation}

By repeatedly and recursively using the above equalities (\ref{differences-3}),
(\ref{differences-4}) and (\ref{differences-5}) along with
Propositions~\ref{column-product-zero} and \ref{column-product-zero-2}, some
tedious easy calculations (which we do not report here
for the sake of brevity) allow us to express all the $2n-3$ vector-components
$\alpha_{p}$, $p=2,\ldots ,n-1$, and $\beta_{p}$, $p=1,\ldots ,n-1$, in terms of
$\alpha_{1}$ only.

Furthermore, since $L^n=U^n=O$, it holds that
\begin{equation}
I + \sum_{h=1}^{n-1} (-\sigma)^h L^h = (I+ \sigma L)^{-1} \ \ \ \ \ {\rm and}
\ \ \ \ \ I + \sum_{h=1}^{n-1} (-\sigma)^h U^h = (I+ \sigma U)^{-1}.
\label{expansion}
\end{equation}
Therefore, using again Propositions~\ref{column-product-zero} and
\ref{column-product-zero-2}, a subsequent substitution into the first equality of
(\ref{FTC_eigen}) yields

\begin{equation}
\big[ L + \big( I+ \sigma U\big )^{-1} U \big] \big( I+ \sigma L \big) ^{-1} \alpha_1
= \lambda \alpha_1,
\label{FTC-eigen-equivalent}
\end{equation}
where $\alpha_{1} \ne 0$ (see Proposition~\ref{nonzero1}).

Vice versa, we assume that (\ref{FTC-eigen-equivalent}) holds for some $\lambda \ne 0$ and
$\alpha_{1} \ne 0$.

Then, by imposing the equalities (\ref{differences-3}), (\ref{differences-4}) and
(\ref{differences-5}), we define the vectors $\alpha_{p}$, $p=2,\ldots ,n-1$,
and $\beta_{p}$, $p=1,\ldots ,n-1$ and, starting from (\ref{FTC-eigen-equivalent})
and using (\ref{expansion}), we can reverse the flow of the previous calculations so
as to arrive at the system (\ref{FTC_eigen}), proving its equivalence to
(\ref{FTC-eigen-equivalent}).

It is immediate to see that, given the identical form of the schemes (\ref{FTC_method})
and (\ref{FTR_method}) (corresponding to the splittings ${\cal B}_{FTC}$ and
${\cal B}_{FTR}$, respectively) and the identical computation rules
imposed by Propositions~\ref{column-product-zero}, \ref{column-product-zero-2},
\ref{column-product-zero-3} and \ref{column-product-zero-4}, completely analogous
calculations lead to the same equivalent eigenvalue equation (\ref{FTC-eigen-equivalent})
also for the iteration matrix ${\cal B}_{FTR}$.

To finish with, we have to prove that (\ref{FTC-eigen-equivalent}) is also equivalent to
$$
B_{sGS} \, \alpha_1 = \lambda \alpha_1,
$$
where
\begin{equation}
B_{sGS}=(I-U)^{-1}L(I-L)^{-1}U.
\label{sGS_matrix}
\end{equation}

To this purpose, we multiply both sides of (\ref{FTC-eigen-equivalent}) by the nonsingular
matrix $(I+\sigma L)(I+\sigma U)$ and, after some more tedious easy calculations,
we arrive at the further equivalent equality
$$
LU \alpha_1 = \lambda (I-L) (I-U) \alpha_1.
$$
Finally, multipliplying both sides by $(I-U)^{-1} (I-L)^{-1}$ immediately concludes the
proof.
\end{proof}

\vskip 0.2cm

The equivalence among the FTC- and FTR-methods (\ref{FTC_method}) and (\ref{FTR_method})
and the symmetric Gauss-Seidel method is not limited to having the same spectrum of
the corresponding iteration matrix.
Indeed, we shall see that the same sequence of approximations $\{y^{(k)}\}_{k\ge 1}$
produced by the symmetric Gauss-Seidel method can be always obtained also by means of
the FTC- or the FTR-method.

In order to prove this fact, we first need to introduce an alternative formulation of
the symmetric Gauss-Seidel iteration which, in a sequential computation environment, has
the additional advantage of halving the computational cost of each step with respect to
a naive standard implementation.

\begin{lemma}
Given the initial value $y^{(0)}$, consider the symmetric Gauss-Seidel iterative scheme
\begin{equation}
y^{(k+1)} = B_{sGS} \, y^{(k)} + c_{sGS}, \ \ \ \ \ k\ge 0,
\label{iteration_sGS}
\end{equation}
where $B_{sGS}$ is given by (\ref{sGS_matrix}) and
\begin{equation}
c_{sGS}=(I-U)^{-1} (I-L)^{-1}c.
\label{c_sGS}
\end{equation}
Then, by setting
\begin{equation}
z^{(k)} := Uy^{(k)}, \ \ \ \ \ k\ge 0,
\label{smart_sGS-1}
\end{equation}
we get
\begin{equation}
y^{(k+1)} = (I-U)^{-1} [(I-L)^{-1} - I] z^{(k)} + c_{sGS}, \ \ \ \ \ k\ge 0,
\label{smart_sGS-2}
\end{equation}
and the {\em modified symmetric Gauss-Seidel} iteration
\begin{equation}
z^{(k+1)} = B_{sGS}^* \, z^{(k)} + Uc_{sGS}, \ \ \ \ \ k\ge 0,
\label{smart_sGS-3}
\end{equation}
where
\begin{equation}
B_{sGS}^* := [(I-U)^{-1} - I] \cdot [(I-L)^{-1} - I]
\label{smart_sGS-4}
\end{equation}
\label{sym-Gauss-Seidel_smart}
is the {\em modified symmetric Gauss-Seidel} iteration matrix.
\end{lemma}
\begin{proof}
The proof is trivial and follows from the obvious equalities
$$
(I-L)^{-1} = I+L(I-L)^{-1} \ \ \ \ \ {\rm and} \ \ \ \ \ (I-U)^{-1} = I+U(I-U)^{-1}.
$$
\end{proof}

\vskip 0.2cm

\begin{remark}
The naive implementation of the method (\ref{iteration_sGS}) using formula
(\ref{sGS_matrix}) for the computation of the approximation $y^{(h)}$ (for a given
$h\ge 1$) costs about twice as much $h$ steps of the Jacobi (or forward or backward
Gauss-Seidel) method.

Instead, it is evident that starting with (\ref{smart_sGS-1}) for $k=0$,
going on with the iteration scheme (\ref{smart_sGS-3}) up to the next to last
approximation $z^{(h-1)}$ and concluding with (\ref{smart_sGS-2}) for $k=h-1$ almost
halves the total computational cost and, thus, reduces it to about the same cost of
$h$ steps of the Jacobi (or forward or backward Gauss-Seidel) method.
\label{sym-Gauss-Seidel_remark}
\end{remark}

\vskip 0.2cm

\begin{theorem}
Let $\{z^{(k)}\}_{k\ge 1}$ be the sequence of approximations produced by the modified
symmetric Gauss-Seidel method (\ref{smart_sGS-3}) starting from the initial value
$z^{(0)}:=Uy^{(0)}$.
Then the sequence of approximations
$$
\big\{X^{(k)}=\big[{x^{(k)}_1}^T, \ldots, {x^{(k)}_{n-1}}^T, {x^{(k)}_n}^T, \ldots,
{x^{(k)}_{2n-2}}^T \big]^T \big\}_{k\ge 0}
$$
produced by the FTC-method (\ref{FTC_method}) with initial value
$$
X^{(0)}:=\big[{x^{(0)}_1}^T, \ldots, {x^{(0)}_{n-1}}^T, {y^{(0)}}^T, \ldots,
{y^{(0)}}^T \big]^T
$$
is such that
\begin{equation}
U_{c}^{(n)}x^{(k)}_{n} + \cdots + U_{c}^{(2)}x^{(k)}_{2n-2} = z^{(k)}
 \ \ \ \ \ {\rm for \ all} \ \ \ k\ge 0,
\label{same-sequence-sGS}
\end{equation}
independently of the first $n-1$ initial vector-components $x^{(0)}_i$, $i=1, \ldots,n-1$.

Analogously, for the FTR-method (\ref{FTR_method}) it holds that
\begin{equation}
U_{r}^{(n-1)}x^{(k)}_{n} + \cdots + U_{r}^{(1)}x^{(k)}_{2n-2} = z^{(k)}
 \ \ \ \ \ {\rm for \ all} \ \ \ k\ge 0.
\label{same-sequence-sGS-2}
\end{equation}
\label{sym-Gauss-Seidel_equivalence}
\end{theorem}
\begin{proof}
The proof is carried out by induction and, since it is rather technical but not difficult,
for the sake of brevity we give here an outline only.

Since $U_{c}^{(n)} + \cdots + U_{c}^{(2)} = U$
and since the chosen initial value satisfies
$$
x^{(0)}_{n} = \cdots = x^{(0)}_{2n-2} = y^{(0)},
$$
it is obvious that (\ref{same-sequence-sGS}) holds for $k=0$.

Now we assume that it holds for $k$ and prove that it is true for $k+1$, too.

To this purpose, by using the iterative scheme (\ref{FTC_method}) and the obvious equality 
$$
I + \sum_{h=1}^{n-1} L^h = (I-L)^{-1}
$$
and by taking Propositions~\ref{column-product-zero} and \ref{column-product-zero-2}
into account, long and tedious calculations allow us to express the second part
$\big[{x^{(k+1)}_n}^T, \ldots, {x^{(k+1)}_{2n-2}}^T \big]^T$ of $X^{(k+1)}$
in terms of the second part $\big[{x^{(k)}_n}^T, \ldots, {x^{(k)}_{2n-2}}^T \big]^T$
only of $X^{(k)}$ and on the vector $c$.
In particular, as an interesting sample of the new vector components $x^{(k+1)}_{i}$,
$i=n, \ldots, 2n-2$, we get
$$
x^{(k+1)}_{n} = (I-L)^{-1} \big( z^{(k)}+c \big).
$$

Therefore, we have just proved that the second part of all the iterates $X^{(k)}$ is
independent of the first $n-1$ initial vector-components $x^{(0)}_i$, $i=1, \ldots,n-1$.

Finally, using also the obvious equality
$$
I + \sum_{h=1}^{n-1} U^h = (I-U)^{-1},
$$
some further computations show that
$$
U_{c}^{(n)}x^{(k+1)}_{n} + \cdots + U_{c}^{(2)}x^{(k+1)}_{2n-2} =
B_{sGS}^* \, z^{(k)} + Uc_{sGS}.
$$
Hence, in view of (\ref{smart_sGS-3}), the induction works. 

The proof of (\ref{same-sequence-sGS-2}) is completely analogous.
\end{proof}

\vskip 0.2cm

All the results displayed for the FTC- and FTR-methods may be repeated in the same way
for the cyclic permutations
$$
{\cal B}'_{FTC} := \{U_{c}^{(n)}, \ldots, U_{c}^{(2)},L_{c}^{(1)}, \ldots, L_{c}^{(n-1)}\}
$$
and
$$
{\cal B}'_{FTR} := \{U_{r}^{(n-1)}, \ldots, U_{r}^{(1)},L_{r}^{(2)}, \ldots, L_{r}^{(n)}\}
$$
of the splittings ${\cal B}_{FTC}$ and ${\cal B}_{FTR}$, respectively, and for
the cyclic modification (\ref{sym_Gauss-Seidel-2}) of the symmetric Gauss-Seidel method
(\ref{sym_Gauss-Seidel}).

%
%

\section{Some further proposals of splitting methods}

\label{section_altrimetodi}

In this section we propose a few particular splitting methods, which have not been
already considered in Section~\ref{T-method} and have a good potential
for possibly attaining faster convergence than the symmetric Gauss-Seidel method
(i.e., the FTC- or FTR-method).

Such methods are not comparable in the partial order relation ''$\succ$''
introduced by Definitions~\ref{order1} and \ref{order2},
either among them or with the methods introduced in Section~\ref{T-method}, including
the T$_U$- and the T$_L$-method.
Therefore, it is not possible to establish a priori which of them all is the fastest,
even if applied to a nonnegative Jacobi iteration matrix.

%
%

\subsection{The alternate triangular column and row methods}

\label{ATC-ATR_methods}

We define the {\em alternate triangular column (ATC) methods} by considering the
T$_L$-method associated to ${\cal B}_{T_L}=\{L,U\}$ (or the T$_U$-method associated to
${\cal B}_{T_U}=\{U,L\}$) and by splitting further both the
lower and the upper triangles $L$ and $U$ separately into the same number of subsets
of columns, which are ordered in an alternate way, one from $L$ and one from $U$.

Analogously, we define the {\em alternate triangular row (ATR) methods} by splitting
further both the lower and the upper triangles $L$ and $U$ separately into the same
number of subsets of rows, which are ordered in an alternate way, one from $L$
and one from $U$.

In particular, as the most promising choice, we consider the lower and upper triangular
column/row matrices introduced by Definitions~\ref{upper_T_C_M}, \ref{lower_T_C_M},
\ref{upper_T_R_M} and \ref{lower_T_R_M} and the splittings
$$
{\cal B}_{AFTC_L} := \{L_{c}^{(1)},U_{c}^{(n)}, \ldots, L_{c}^{(n-1)},U_{c}^{(2)} \},
$$
$$
{\cal B}_{AFTC_U} := \{U_{c}^{(n)},L_{c}^{(1)}, \ldots, U_{c}^{(2)},L_{c}^{(n-1)} \},
$$
$$
{\cal B}_{AFTR_L} := \{L_{r}^{(2)},U_{r}^{(n-1)}, \ldots, L_{r}^{(n)},U_{r}^{(1)} \},
$$
$$
{\cal B}_{AFTR_U} := \{U_{r}^{(n-1)},L_{r}^{(2)}, \ldots, U_{r}^{(1)},L_{r}^{(n)} \},
$$
which are obtained by alternating one column (row) of $L$ to one column (row) of $U$
in the same order used for the FTC- and FTR-methods.
We call them the {\em alternate full triangular column splitting starting from $L$
(AFTC$_L$)} or {\em starting from $U$ (AFTC$_U$)} and the
{\em alternate full triangular row splitting starting from $L$
(AFTR$_L$)} or {\em starting from $U$ (AFTR$_U$)}, respectively.

However, all the above splittings are not essential (see Definition~\ref{essentiality}).
In fact, it is immediate to prove the following result.
\begin{proposition}
Given the Jacobi iteration matrix $B_J=L+U \in \re^{n\times n}$, it holds that
$$
U^{(j)}_{c}L^{(k)}_{c}=O \ \ \ \ \ {\rm for \ all} \ \ \ 2\le j \le k \le n-1,
$$
$$
L^{(j)}_{c}U^{(k)}_{c}=O \ \ \ \ \ {\rm for \ all} \ \ \ 2\le k \le j \le n-1,
$$
$$
U^{(i)}_{r}L^{(l)}_{r}=O \ \ \ \ \ {\rm for \ all} \ \ \ 2\le l \le i \le n-1,
$$
$$
L^{(i)}_{r}U^{(l)}_{r}=O \ \ \ \ \ {\rm for \ all} \ \ \ 2\le i \le l \le n-1.
$$
\label{UL=LU=O}
\end{proposition}

Therefore, by using Proposition~\ref{prodotto_zero_raggio}, easy calculations allow us to
conclude that the iteration matrices corresponding to the above splittings have the same
nonzero eigenvalues of those corresponding to
$$
{\cal B}'_{AFTC_L} := \{L_{c}^{(1)},U_{c}^{(n)}, \ldots, L_{c}^{(n-\nu-1)},U_{c}^{(\nu+2)},
L_{c}^{(n-\nu)}+U_{c}^{(\nu+1)}, \ldots, L_{c}^{(n-1)}+U_{c}^{(2)}\},
$$
$$
{\cal B}'_{AFTC_U} := \{U_{c}^{(n)},L_{c}^{(1)}, \ldots, U_{c}^{(\nu+2)},L_{c}^{(n-\nu-1)},
L_{c}^{(n-\nu)}+U_{c}^{(\nu+1)}, \ldots, L_{c}^{(n-1)}+U_{c}^{(2)}\},
$$
$$
{\cal B}'_{AFTR_L} := \{L_{r}^{(2)}+U_{r}^{(n-1)}, \ldots, L_{r}^{(\nu+1)}+U_{r}^{(n-\nu)},
L_{r}^{(\nu+2)},U_{r}^{(n-\nu-1)}, \ldots, L_{r}^{(n)},U_{r}^{(1)}\},
$$
$$
{\cal B}'_{AFTR_U} := \{L_{r}^{(2)}+U_{r}^{(n-1)}, \ldots, L_{r}^{(\nu+1)}+U_{r}^{(n-\nu)},
U_{r}^{(n-\nu-1)},L_{r}^{(\nu+2)}, \ldots, U_{r}^{(1)},L_{r}^{(n)}\},
$$
respectively, which have $2n- 2 -\nu= O(3n/2)$ elements, where $\nu$ is given by (\ref{nu}).

These splittings turn out to be maximal in the partial order relation ''$\succ$'' introduced by
Definitions~\ref{order1} and \ref{order2}.

Moreover, differently from the case of the FTC- and FTR-methods, all the four
AFTC$_L$-, AFTC$_U$-, AFTR$_L$- and AFTR$_U$-methods generally have different speeds
of convergence.

Some numerical results will be given in Section~\ref{numerics}.

%
%

\section{Numerical examples}

\label{numerics}

In this section we give some numerical examples based on three classes of matrices $A$,
say, {\em Class}~1, {\em Class}~2 and {\em Class}~3, which do not necessarily play a
particular role in applications, but are suitable enough to illustrate the developed
theory.

The matrices of {\em Class}~1 are constructed by choosing {\em uniformly distributed random elements} such that
$$
a_{i,j}\in [-1,1] \ \ \ \ \ \forall \ (i,j) \ \ \ {\rm with} \ \ \ i\ne j.
$$
Moreover, the diagonal elements are computed by using the rule
\begin{equation}
a_{i,i} := \frac{1}{\phi} \sum_{j=1,\, j\ne i}^n |a_{i,j}| \, \ \ \ \ \ \forall \ i,
\label{diag_elements}
\end{equation}
where $\phi$ is a given parameter.
\begin{remark}
For the corresponding Jacobi iteration matrix it always holds that
$$
\|B_J\|_{\infty} = \phi.
$$

In particular, if $\phi < 1$, then the matrix $A$ is strictly diagonally dominant.
\label{rowfact}
\end{remark}

\vskip 0.2cm

For each matrix $A$ of {\em Class}~1 we define the corresponding matrix $A'$ of
{\em Class}~2 by setting
$$
a_{i,j}' := -|a_{i,j}| \ \ \ \ \ \forall \ (i,j) \ \ \ {\rm with} \ \ \ i\ne j
$$
and
$$
a_{i,i}' := a_{i,i} \ \ \ \ \ \forall \ i
$$
and the corresponding matrix $A''$ of {\em Class}~3 by setting
$$
a_{i,j}'' := |a_{i,j}| \ \ \ \ \ \forall \ (i,j) \ \ \ {\rm with} \ \ \ i\ne j
$$
and, again,
$$
a_{i,i}'' := a_{i,i} \ \ \ \ \ \forall \ i.
$$

\begin{remark}
The matrices belonging to Class~2 are such that $B_J \ge O$ elementwise and those
belonging to Class~3 are such that $B_J \le O$ elementwise. Therefore, the theory
developed in Section~\ref{sect_L-matrices} applies.

Furthermore, in both cases we have that
$$
\rho(B_J) = \|B_J\|_{\infty} = \phi.
$$
\label{neg-pos}
\end{remark}

The elements of the matrices being randomly chosen, we repeat each experiment for
$N=100$ times in order to give it a reasonable statistical meaning.

For each matrix $A$, $A'$ and $A''$ we compute the spectral radius $\rho(B_J)$ of the
corresponding Jacobi iteration matrix $B_J$ along with the spectral radius
$\rho({\cal B})$ of the iteration matrix ${\cal B}$ for all the other considered methods.

Furthermore, we compute the {\em speedup factors}
\begin{equation}
sp({\cal B}) := log(\rho({\cal B}))/log(\rho(B_J))
\label{speedup}
\end{equation}
of the various methods with respect to Jacobi.

However, since we make many experiments using different matrices
$A^{(r)}$, $r=1,\ldots ,N$, for each method we only show the
{\em average values}
\begin{equation}
\overline{\rho}(B_J) := \frac{1}{N} \sum_{r=1}^N \rho(B_J^{(r)}) \, , \ \ \ \ \
\overline{\rho}({\cal B}) := \frac{1}{N} \sum_{r=1}^N \rho({\cal B}^{(r)})
\label{average_rho}
\end{equation}
and
\begin{equation}
\overline{sp}({\cal B}) := \frac{1}{N} \sum_{r=1}^N sp({\cal B}^{(r)})
\label{average_speedup}
\end{equation}
of the spectral radii and of the speedup factors, respectively, together with the
corresponding {\em standard deviations}
\begin{equation}
s.d.{\rho}(B_J) := \sqrt{\frac{1}{N} \sum_{r=1}^N \Big( \rho(B_J^{(r)}) -
\overline{\rho}(B_J) \Big)^2} \, , \ \ \ \ \
s.d.{\rho}({\cal B}) := \sqrt{\frac{1}{N} \sum_{r=1}^N \Big( \rho({\cal B}^{(r)}) -
\overline{\rho}({\cal B}) \Big)^2}
\label{s.d._rho}
\end{equation}
and
\begin{equation}
s.d.{sp}({\cal B}) := \sqrt{\frac{1}{N} \sum_{r=1}^N \Big( sp({\cal B}^{(r)}) -
\overline{sp}({\cal B}) \Big)^2} \, .
\label{s.d._speedup}
\end{equation}

The numerical values of the average spectral radii and standard deviations are given
in the fixed point notation, truncated to the 5th significant digit, whereas the
standard deviations are given in the floating point notation, truncated to the 3rd
significant digit.

Besides Jacobi we consider, in the following order, the methods T$_U$, fGS, bGS, TC(2,2),
TR(2,2), sGS, AFTC$_L$, AFTC$_U$, AFTR$_L$, AFTR$_U$ in relation to the value
$$
\phi = 0.9
$$
in the definition (\ref{diag_elements}) of the diagonal elements.
The chosen dimension of the matrices is $n=100$.

\vskip 1cm


\begin{center}
\begin{tabular}{||c||}
\hline
Experiments with matrices $A$ of Class~1 for $\phi=0.9$ \\
\hline
\end{tabular}
\end{center}
\begin{center}
\begin{tabular}{||c|c|c|c|c||}
\hline
\ \ method \ \ & \ \ $\overline{\rho}({\cal B})$ \ \ & \ \ $s.d.{\rho}({\cal B})$ \ \ &
$\ \ \overline{sp}({\cal B})$ \ \ & \ \ $s.d.{sp}({\cal B})$ \ \ \\
\hline\hline
Jacobi & 0.10962 \ \ & 3.85\,E-03 & ===== & ===== \\
\hline
T$_U$ & 0.057121 \ & 1.90\,E-03 & 1.2950 & 2.12\,E-02 \\
\hline
fGS & 0.042714 \ & 1.56\,E-03 & 1.4296 & 2.84\,E-02 \\
\hline
bGS & 0.042434 \ & 1.54\,E-03 & 1.4296 & 2.84\,E-02 \\
\hline
TC(2,2) & 0.043724 \ & 1.61\,E-03 & 1.4160 & 2.80\,E-02 \\
\hline
TR(2,2) & 0.043949 \ & 1.58\,E-03 & 1.4137 & 2.76\,E-02 \\
\hline
sGS & 0.0075707 & 5.93\,E-04 & 2.2103 & 5.04\,E-02 \\
\hline
AFTC$_L$ & 0.032672 \ & 1.26\,E-03 & 1.5478 & 2.77\,E-02 \\
\hline
AFTC$_U$ & 0.032815 \ & 1.24\,E-03 & 1.5458 & 2.59\,E-02 \\
\hline
AFTR$_L$ & 0.032552 \ & 1.34\,E-03 & 1.5496 & 3.03\,E-02 \\
\hline
AFTR$_U$ & 0.032762 \ & 1.26\,E-03 & 1.5466 & 2.98\,E-02 \\
\hline
\end{tabular}
\end{center}

\newpage


\begin{center}
\begin{tabular}{||c||}
\hline
Experiments with matrices $A'$ of Class~2 for $\phi=0.9$ \\
\hline
\end{tabular}
\end{center}
\begin{center}
\begin{tabular}{||c|c|c|c|c||}
\hline
\ \ method \ \ & \ \ $\overline{\rho}({\cal B})$ \ \ & \ \ $s.d.{\rho}({\cal B})$ \ \ &
$\ \ \overline{sp}({\cal B})$ \ \ & \ \ $s.d.{sp}({\cal B})$ \ \ \\
\hline\hline
Jacobi & 0.90000 & 0 & ===== & ===== \\
\hline
T$_U$ & 0.85418 & 1.25\,E-04 & 1.4960 & 1.39\,E-03 \\
\hline
fGS & 0.81286 & 5.82\,E-04 & 1.9670 & 6.77\,E-03 \\
\hline
bGS & 0.81282 & 5.79\,E-04 & 1.9670 & 6.77\,E-03 \\
\hline
TC(2,2) & 0.82388 & 3.22\,E-04 & 1.8387 & 3.71\,E-03 \\
\hline
TR(2,2) & 0.82385 & 3.53\,E-04 & 1.8391 & 4.07\,E-03 \\
\hline
sGS & 0.73472 & 3.50\,E-04 & 2.9259 & 5.08\,E-03 \\
\hline
AFTC$_L$ & 0.78174 & 5.80\,E-04 & 2.3370 & 7.04\,E-03 \\
\hline
AFTC$_U$ & 0.78179 & 5.88\,E-04 & 2.3365 & 7.14\,E-03 \\
\hline
AFTR$_L$ & 0.78167 & 5.43\,E-04 & 2.3379 & 6.59\,E-03 \\
\hline
AFTR$_U$ & 0.78162 & 5.45\,E-04 & 2.3385 & 6.62\,E-03 \\
\hline
\end{tabular}
\end{center}

\vskip 1cm


\begin{center}
\begin{tabular}{||c||}
\hline
Experiments with matrices $A''$ of Class~3 for $\phi=0.9$ \\
\hline
\end{tabular}
\end{center}
\begin{center}
\begin{tabular}{||c|c|c|c|c||}
\hline
\ \ method \ \ & \ \ $\overline{\rho}({\cal B})$ \ \ & \ \ $s.d.{\rho}({\cal B})$ \ \ &
$\ \ \overline{sp}({\cal B})$ \ \ & \ \ $s.d.{sp}({\cal B})$ \ \ \\
\hline\hline
Jacobi & 0.90000 \ & 0 & ===== & ===== \\
\hline
T$_U$ & 0.40932 \ & 8.36\,E-04 & \ \ \ 8.4782 & 1.94\,E-02 \\
\hline
fGS & 0.19544 \ & 2.23\,E-03 & 15.498 & 1.07\,E-01 \\
\hline
bGS & 0.19537 \ & 2.20\,E-03 & 15.498 & 1.07\,E-01 \\
\hline
TC(2,2) & 0.22573 \ & 1.87\,E-03 & 14.127 & 7.88\,E-02 \\
\hline
TR(2,2) & 0.22555 \ & 1.94\,E-03 & 14.135 & 8.19\,E-02 \\
\hline
sGS & 0.17146 \ & 9.75\,E-04 & 16.737 & 5.39\,E-02 \\
\hline
AFTC$_L$ & 0.098689 & 1.78\,E-03 & 21.981 & 1.72\,E-01 \\
\hline
AFTC$_U$ & 0.098751 & 1.76\,E-03 & 21.975 & 1.70\,E-01 \\
\hline
AFTR$_L$ & 0.098327 & 1.72\,E-03 & 22.016 & 1.66\,E-01 \\
\hline
AFTR$_U$ & 0.098319 & 1.70\,E-03 & 22.017 & 1.65\,E-01 \\
\hline
\end{tabular}
\end{center}

\vskip 1cm

The first general observation we feel to make is that the standard deviations always
are small enough in comparison to the average values of the various spectral radii and
speedup factors, meaning that the above three classes of matrices are substantially
stable with respect to these two evaluation indices of an iterative method.

Looking at the experiments on the above three classes of matrices, the second general
observation is that the heuristic idea for which the more the splitting is refined,
the faster is the convergence rate, seems to be substantially confirmed.
Nevertheless, a fixed score among the various methods cannot be established.

It is clear enough that Jacobi is the slowest method, followed by the T$_U$-method
with a good improvement, especially when applied to matrices of Class~3.
However, these two methods present an evident advantage with respect to the more refined
ones if they are used in a parallel computation environment.

Then it looks like that the TC(2,2) and TR(2,2) methods could be a good alternative
to forward and backward Gauss-Seidel since, although being a little bit slower, they
also present a clear advantage if used in a parallel computation environment.

It also seems to be quite frequent that symmetric Gauss-Seidel be the fastest method,
closely followed by the four variants of AFTC- and AFTR- methods (see the experiments
with matrices of Class~1 and Class~2).
The second ones perform better on matrices of Class~3.
In any case, the general better performance of these two types of methods is consistent
to the fact that, among those we have considered, they are the only ones that are
potentially optimal (in the sense of Definition~\ref{optimality}).

Finally, we remark that the experiments we made within Class~2 clearly confirm the
theoretical speed score implied by Theorem~\ref{faster} and that, given the mutual
relation among the matrices of Class~2 and Class~3, the expected inequality
(\ref{rhoB_d<=rhoB_d'}) is largely verified in the strict form for all the considered
methods except, of course, Jacobi (see Proposition~\ref{Jacobi-nonpos}).

In connection to the observations we made at the end of Section~\ref{subsect_nonpos},
now we present some results that illustrate the evident superiority of all the considered
splitting methods with respect to Jacobi when they are applied to matrices of Class~3.
Again, we consider matrices of dimension $n=100$.

Starting from $\phi=1$, we consider a grid of values for the parameter $\phi$ with fixed 
increment $\delta=0.5$ and, for the sake of space, we report only the results related
to those instances of $\phi$ that determine the exceeding of the critical value
$\rho({\cal B})=1$ for one or more of the methods under consideration.
We highlight in boldface all the values $\rho({\cal B})>1$.

The experiments are based on a particular random matrix $A''$ of Class~3 corresponding
to the value $\phi=1$ and, keeping fixed all of its off-diagonal elements,
by using formula (\ref{diag_elements}) to pass from $\phi$ to the next value
$\phi + \delta$. Thus we only change the diagonal elements and,
consequently, $\rho(B_J)$ (see Remark~\ref{neg-pos}).

This time we do not need to make any statistics because our goal is to analyze
the spectral radii $\rho({\cal B})$ as functions of the sole parameter $\phi$.

\vskip 1cm


\begin{center}
\begin{tabular}{||c||}
\hline
$\rho({\cal B})$ as function of $\phi=\rho(B_J)$ for a matrix $A''$ of Class~3 \\
\hline
\end{tabular}
\end{center}
\begin{center}
\begin{tabular}{||c|c|c|c|c|c||}
\hline
\ method \ & \ \ $\phi=2.0$ \ & \ \ $\phi=3.5$ \ &
\ $\phi=14.0$ \ & \ $\phi=15.0$ \ & \ $\phi=15.5$ \ \\
\hline\hline
T$_U$ & {\bf 1.0175} & {\bf 1.8977} & {\bf 67.415} & {\bf 78.376} & {\bf 84.160} \\
\hline
fGS & 0.40715 & 0.60249 & 0.95948 & {\bf 1.0509} & {\bf 1.1194} \\
\hline
bGS & 0.39091 & 0.57263 & {\bf 1.0043} & {\bf 1.0621} & {\bf 1.1080} \\
\hline
TC(2,2) & 0.52197 & {\bf 1.2888} & {\bf 387.61} & {\bf 530.03} & {\bf 614.41} \\
\hline
TR(2,2) & 0.51804 & {\bf 1.2774} & {\bf 378.11} & {\bf 517.14} & {\bf 599.52} \\
\hline
sGS & 0.43217 & 0.65676 & 0.95397 & 0.94779 & {\bf 1.0011} \\
\hline
AFTC$_L$ & 0.16267 & 0.45879 & 0.92226 & {\bf 1.0288} & {\bf 1.0993} \\
\hline
AFTC$_U$ & 0.17824 & 0.45429 & 0.92915 & {\bf 1.0644} & {\bf 1.1364} \\
\hline
AFTR$_L$ & 0.18127 & 0.44260 & 0.92957 & {\bf 1.0577} & {\bf 1.1346} \\
\hline
AFTR$_U$ & 0.17969 & 0.44696 & 0.92392 & {\bf 1.0539} & {\bf 1.1319} \\
\hline
\end{tabular}
\end{center}

\vskip 1cm

Our experiments suggest that, as the parameter $\phi$ increases, the resistance of the
methods to exceed the critical value $\rho({\cal B})=1$ is generally stronger for those
that are more refined and that the better performance of the four variants of AFTC- and
AFTR- methods tends to disappear in favour of symmetric Gauss-Seidel.

Finally, we give some numerical results referred to an application
to B-spline approximation (Wang et al. [WHV]).

The underlying matrix $A$, still of dimension $n=100$, is a particular instance of
matrices of Class~3. In fact, it is symmetric and 9-diagonal with elements
$$
[\ 0\ \cdots \ 0\ \ 1\ \ 4\ \ 1\ \ 4\ \ 16\ \ 4\ \ 1\ \ 4\ \ 1\ \ 0\ \cdots \ 0\ ]
$$
in each row, where the diagonal one is $a_{i,i}=16$ (apart from the first and the
last four rows, which are modified in an obvious way).

The particular form of $A$ makes many of the spectral radii of the iteration matrices
of the considered splitting methods to coincide among them. More precisely,
\begin{eqnarray*}
\rho({\cal B}_{fGS}) & = & \rho({\cal B}_{bGS}), \\
\rho({\cal B}_{TC(2,2)}) & = & \rho({\cal B}_{TR(2,2)}), \\
\rho({\cal B}_{AFTC_L}) = \rho({\cal B}_{AFTC_U}) & = & \rho({\cal B}_{AFTR_L})
= \rho({\cal B}_{AFTR_U}).
\end{eqnarray*}

The matrix $A$ is not diagonally dominant and
$$
\rho(B_J)=1.2464... > 1,
$$
so that Jacobi's method does not converge, whereas all the other methods do.
The spectral radii of the respective iteration matrices are given in the following table.

\vskip 1cm


\begin{center}
\begin{tabular}{||c|c|c|c|c||}
\hline
\ \ method \ \ & \ \ $\rho({\cal B})$ \ \ \\
\hline\hline
T$_U$ & 0.68383 \\
\hline
fGS & 0.56821 \\
\hline
TC(2,2) & 0.68087 \\
\hline
sGS & 0.35876 \\
\hline
AFTC$_L$ & 0.38260 \\
\hline
\end{tabular}
\end{center}

\vskip 1cm

We remark that, on this particular form of matrices of Class~3, symmetric Gauss-Seidel
seems to perform slightly better than the AFTC-methods.

%
%

\section{Conclusions}

\label{conclusions}

In this paper, taking inspiration from the T$_U$-method defined in \cite{Ta}, we have
defined a general class of {\em splitting} methods which, starting from Jacobi,
include, among others, all the well-known methods of Gauss-Seidel (forward, backward
and symmetric) and have the property to share all the same cost in a sequential
computation environment.

A theoretical ranking of the convergence properties of all these methods is given
on the basis of a certain {\em refinement} partial order relation which definitely works
when they are applied to matrices $A$ such that the corresponding Jacobi iteration matrix
$B_J$ is nonnegative (often called L-matrices).

Two particular proposals of new methods have been done, namely the TC(2,2)- and
TR(2,2)-method and the four variants of AFTC- and AFTR-methods, that seem to be
promising enough.

It is clear that further work could still be done to improve the potential of the
proposed splitting methods such as, for example, the use of a relaxation
parameter $\omega$.

Finally, also an accurate comparison of the performances of the various splitting methods
would be worth doing assuming to work in parallel computation environments.


\section*{Aknowledgement}
This work was partially supported by GNCS-INdAM and FRA-University of Trieste.
Paolo Novati is member of the INdAM research group GNCS.


\newcommand{\etalchar}[1]{$^{#1}$}

\end{document}